\documentclass[a4paper,12pt]{amsart}
\usepackage{amsmath,amsfonts,amssymb,amsthm}
\usepackage{graphics}
\usepackage{epsfig}
\usepackage{esint}
\usepackage{shuffle}
\usepackage{istgame}
\usepackage{vmargin}
\setpapersize{A4}
\setmargins{2.5cm}{1.5cm}{16.5cm}{23.42cm}{10pt}{1cm}{0pt}{2cm}

\usepackage{enumerate}

\newcommand{\ignore}[1]{}

\usepackage[hidelinks,hyperindex,breaklinks]{hyperref}

\newtheorem{proposition}{Proposition}
\newtheorem{theorem}{Theorem}

\newtheorem{lemma}{Lemma}
\newtheorem{corollary}{Corollary}

\title{On minimizing curves in a Brownian potential}
\author{Felix Otto, Matteo Palmieri, and Christian Wagner}

\begin{document}

\begin{abstract}

	We study a $(1+1)$-dimensional semi-discrete random variational problem 
	that can be interpreted as the geometrically linearized version of 
	the critical $2$-dimensional random field Ising model.~The scaling
	of the correlation length of the latter was recently characterized
	in \cite{dw} and \cite[Section 5]{DingarXiv}; our analysis is reminiscent of 
	the multi-scale approach of the latter work and of
	\cite{LeightonShor}. We show that at
	every dyadic scale from the system size down to the lattice spacing
	the minimizer contains at most order-one Dirichlet energy per unit length.
	We also establish a quenched homogenization result in the sense that
	the leading order of the minimal energy becomes deterministic
	as the ratio system size / lattice spacing diverges.
	To this purpose we adapt arguments from \cite{Peled} on the $(d+1)$-dimensional version
 	our the model, with a Brownian replacing the white noise potential,
	to obtain the initial large-scale bounds. Based on our estimate of the $(p=3)$-Dirichlet
	energy, we give an informal justification of the geometric linearization.
        Our bounds, which are oblivious to the microscopic cut-off scale provided by 
	the lattice spacing, yield tightness of the law of minimizers 
	in the space of continuous functions as the lattice spacing is sent to zero.

\end{abstract}

\maketitle

\section{Setting and main results}

\subsection{Setting}
Consider the one-dimensional lattice of size $L$,
\begin{align*}
\{0,\cdots,L\}\ni x;
\end{align*}
our configuration space is given by all height functions $h$
\begin{align}\label{ao49}
\{0,\cdots,L\}\ni x\mapsto h(x)\in\mathbb{R}\quad\mbox{with}\quad h(0)=h(L)=0
\end{align}
that vanish at the two boundary points.
Introducing the (discrete) Dirichlet energy
\begin{align}\label{ao51}
D(h):=\frac{1}{2}\sum_{x=1}^L(h(x)-h(x-1))^2
\end{align}
and a field term of the form  
\begin{align}\label{ao53}
W(h):=\sum_{x=1}^{L-1} W(x,h(x)),
\end{align}
we are interested in the variational problem of
\begin{align}\label{ao48}
\mbox{minimize}\;\;E(h):=D(h)-W(h)\;\;\mbox{among all $h$ of the form (\ref{ao49})}.
\end{align}
Interpreting $x$ as a time variable, $E=D-W$ can be viewed as the Lagrangian of a 
particle in (time-dependent) potential.
The model (\ref{ao48}) can also be interpreted as a directed polymer $h$ at zero
thermal temperature in a potential. Like for directed polymers, we choose
the potential $W$ to be random. In \cite{Peled}, Dembin, Elboim, Hadas, and Peled 
consider lattices also of higher dimension so that (\ref{ao48}) 
models\footnote{with the surface energy
replaced by $D$, its quadratic approximation} a membrane in a potential.
Compared to \cite{Peled} we consider a different ensemble for the potential, namely
\begin{align}\label{ao50}
\{W(x,\cdot)\}_{x=1,\cdots,L-1}\quad
\mbox{are independent two-sided Brownian motions}.
\end{align}
Alternatively, we can write $W(h)$ $=\sum_x\int_0^{h(x)}\xi(x, dy)$
where $\xi$ is semi-discrete white noise. 
As will be explained in Subsection \ref{SS:der},
we interpret (\ref{ao48}) as a (1+1)-dimensional version of the 2-dimensional
random field Ising model, cf.~Subsection \ref{SS:RFIM} for its definition; 
(\ref{ao48}) zooms in on a flatter region of the domain boundary,
represents the boundary by a height function $h$, and replaces the boundary length
by its quadratic approximation, the Dirichlet energy. This approximation amounts
to the geometric linearization of the length
\begin{align}\label{ao90}
\sum_{x=1}^{L}|(x,h(x))-(x-1,h(x-1))|\approx L+D(h).
\end{align}

\medskip

As in \cite[Lemma 5.4]{Peled} we have
\begin{lemma}\label{L:7}
Almost surely, the variational problem \eqref{ao48} admits a unique minimizer.
\end{lemma}

\ignore{
Indeed\footnote{The notation $A\gg_{\alpha_1, \alpha_2, \ldots} B$, used input
statements, means that there exists a constant $C = C(\alpha_1, \alpha_2,
\ldots)$ such that the statement holds if $A \geq CB$. The
notation $A \lesssim_{\beta_1, \beta_2, \ldots} B$, used in output statements,
means that there exist $C' = C'(\beta_1, \beta_2, \ldots)$ such that $A \leq
C'B$. In absence of subscripts, the constants $C, C'$ must be thought as
universal. Furthermore, we shall write $ A \sim B $ to mean that 
both $ A \lesssim B $ and $ B \lesssim A $ hold.},
by the law of the iterated logarithm\footnote{The subscript $\omega$ indicates that
the constant depends on the realization of the Brownian motions (\ref{ao50}).} 
\begin{align*}
	W(h) & \lesssim_{\omega, L} \sum_{x = 1}^L \sqrt{|h(x)|\,\ln\ln(e+|h(x)|)} 
\lesssim_{L} D(h)^{1/4} \ln \ln (e+D(h)),
\end{align*}
$E$ is (almost surely) coercive, that is,
\begin{align*}
E(h) = D(h) - W(h) \geq \frac{1}{2}D(h)\quad\mbox{whenever}\;D(h) \gg_{\omega, L} 1.
\end{align*}
%
}
which allows us to use the notation
\begin{align}\label{ao84}
h_*={\rm argmin} E.
\end{align}

\medskip

After the first arXiv post of this work appeared, Dembin, Elboim, and Peled posted \cite{Peled2}, which extends \cite{Peled} to \eqref{ao48}. They allow for general dimension and codimension, which in our case are both equal to one. Furthermore, they cover fractional Brownian motion, whilst we only consider a special case as motivated in Section \ref{S:Mot}. Our main result focus on different aspects of the model. We refer to Section \ref{S:strategy_of_proof} for a more in-depth comparison with the results of \cite{Peled2}.


\subsection{Scale invariances}\label{S:scal}

Thanks to the continuous nature of the configuration space (\ref{ao49}) in the
vertical/height variable there is an exact scale invariance: Under the change of variables
%
\begin{align}\label{ao58}
h=H\hat h
\end{align}
the functionals (\ref{ao51}) and (\ref{ao53}) transform as
\begin{align}\label{ao54}
D=H^2\hat D\quad\mbox{and}\quad W=H^\frac{1}{2}\hat W\;\;\mbox{in law},
\end{align}
where the latter follows from the scale invariance in law of Brownian motion, cf.~(\ref{ao50}).
We learn from (\ref{ao54}) that there was no loss of generality in ignoring prefactors
in $D-W$, see also Section \ref{SS:der}. 
Since for $H\gg 1$, we have $H^2\gg H^\frac{1}{2}$, the scaling (\ref{ao54})
reflects the fact that the $D-W$ is (almost surely) bounded from below.

\medskip

Due to the discrete nature of the configuration space (\ref{ao49}) in
the horizontal variable $x$, the problem setting comes with two input scales: 
the microscopic scale induced by the lattice (which we set to unity), 
and the macroscopic scale given by the system size $L$. 
We momentarily send the microscopic scale to zero, by replacing the configuration space 
(\ref{ao49}) by the fully continuum one $x\in[0,L]\mapsto h(x)\in\mathbb{R}$, 
holding on to the boundary conditions, and replace (\ref{ao51}) \& (\ref{ao53}) by
\begin{align}\label{ao47}
D_{cont}(h)=\int_0^L dx(\frac{dh}{dx})^2
\quad\mbox{and}\quad W_{cont}(h)=\int_0^L\int_0^{h}d\xi_{cont}
\end{align}
and where $\xi_{cont}$ is (continuum) white noise. 
This fully continuum version allows to complement (\ref{ao58}) by a rescaling of $x$
\begin{align}\label{ao140}
h=l\hat h\quad\mbox{and}\quad
x=l\hat x\;\;\mbox{and thus}\;\;\hat L=\frac{L}{l}.
\end{align}
where $l$ is a (mesoscopic) length scale under which
(\ref{ao47}) transform in an identical way:
\begin{align}\label{ao60}
D_{cont}=l\hat D_{cont}
\quad\mbox{and}\quad W_{cont}=l\hat W_{cont}\;\;\mbox{in law},
\end{align}
where the last statement follows from the scale invariance of white noise $\xi_{cont}$.
In particular, (\ref{ao60}) suggests that the dependence on the macroscopic scale $L$
is determined by scaling:
Ignoring the small-scale divergence when minimizing $D_{cont}-W_{cont}$, 
see Subsection \ref{SS:Tal}, (\ref{ao60}) suggests that the minimizer $h_{cont,*}$ satisfies
\begin{align}\label{ao61}
\mbox{the law of}\quad
\frac{D_{cont}(h_{cont,*})}{L}
\quad\mbox{is independent of $L$},
\end{align}
which would amount to what is called an extensive scaling, i.~e.~proportionality
(in law) of the energy contributions to the system ``volume'' $L$.


\subsection{A maximum of Gaussians}\label{SS:Tal}
We reformulate our variational problem (\ref{ao48}) as
\begin{align}\label{ao56}
\frac{1}{L}{\rm inf}_h E(h)
={\rm inf}_{\hat D\in[0,\infty)}\big(\hat D-\sup_{h:\frac{D(h)}{L}\le \hat D}\frac{W(h)}{L}\big),
\end{align}
where the normalization by $L$ is motivated by the extensive scaling (\ref{ao61}).
Appealing to (\ref{ao58}) with $H=\sqrt{\hat D}$ we obtain from (\ref{ao54}) for the inner supremum
\begin{align}\label{ao63}
\sup_{h:\frac{D(h)}{L}\le \hat D}\frac{W(h)}{L}=\hat D^\frac{1}{4}
\sup_{\hat h:\frac{D(\hat h)}{L}\le 1}\frac{W(\hat h)}{L}\quad\mbox{in law}.
\end{align}
Hence there is a close connection between our variational problem (\ref{ao48})
and the maximum (\ref{ao63}) over a family of centered Gaussian variables $\frac{W(\hat h)}{L}$.
From (\ref{ao50}), we obtain for the variance of such a variable
\begin{align}\label{ao65}
\mathbb{E}(\frac{W(\hat h)}{L})^2=\frac{1}{L}\sum_{x=1}^{L-1}|\frac{\hat h(x)}{L}|
\le(\frac{D(\hat h)}{L})^\frac{1}{2},
\end{align}
where the inequality can be seen as a discrete version of a Poincar\'e inequality
\footnote{which is very easy to establish in our one-dimensional setting:
$\frac{1}{L^2}\sum_{x}|\hat h(x)|$ $\le\frac{1}{L}\max_x|\hat h(x)|$ 
$\le\frac{1}{L}\sum_x|\hat h(x)-\hat h(x-1)|$
$\le(\frac{1}{L}\sum_x(\hat h(x)-\hat h(x-1))^2)^\frac{1}{2}$}. In particular,
in the r.~h.~s.~supremum in (\ref{ao63}) all Gaussians have variance $\le 1$. Hence
we obtain from concentration, cf.~(\ref{c2}) in Lemma \ref{L:4}, and the obvious\footnote{
select a ''macroscopic'' $\hat h$ with $\frac{1}{L}\sum_{x=1}^{L-1}|\frac{\hat h(x)}{L}|\sim 1$ 
but $\frac{D(\hat h)}{L}\le 1$, appeal to (\ref{ao65}), and to $\mathbb{E}X^2$ 
$\le e\|X\|_s^2$ as a consequence of (\ref{ao24})}
$\big\|\sup_{\hat h:\frac{D(\hat h)}{L}\le 1}\frac{W(\hat h)}{L}\big\|_2$ $\gtrsim 1$,
that\footnote{here and in the following by $A \lesssim B$ we mean that there
exists a universal constant $C$ such that $A \leq CB$, while $A \sim B$ is used
if both $A \lesssim B$ and $B \lesssim A$ hold.}
\begin{align}\label{ao81}
\big\|\sup_{\hat h:\frac{D(\hat h)}{L}\le 1}\frac{W(\hat h)}{L}\big\|_2
\sim1+\mathbb{E}\sup_{\hat h:\frac{D(\hat h)}{L}\le 1}\frac{W(\hat h)}{L},
\end{align}
where we introduced
a norm\footnote{Note that $\|\cdot\|_s$ is homogeneous by construction and convex
as a consequence of the convexity of
$X\mapsto\exp(|X|^s)$ -- and thus indeed satisfies the triangle inequality which will
be important. Note that the r.~h.~s.~$e$ in (\ref{ao24}) is chosen such that
$\|X\|_{ s } =|X|$ for a deterministic/constant $X$, which by Jensen's inequality provides a normalization and yields
$|\mathbb{E}X|$ $\le\|X\|_s$.}
on random variables $X$ that captures stretched exponential tails
of exponent $s\ge 1$ by the Orlicz-type definition
\begin{align}\label{ao24}
\|X\|_s:=\inf\{\,\nu>0\,|\,\mathbb{E}\exp(|\frac{X}{\nu}|^s)\le e\,\}.
\end{align}

\medskip

There is a closer relation between the supremum of Gaussians (\ref{ao63}) 
and our variational problem (\ref{ao49}):
In terms of scaling and on the level of the norms (\ref{ao24}) we have

\begin{lemma}\label{L:5} For $s\ge 1$
\begin{align}\label{ao83}
\big\|\frac{D(h_*)}{L}\big\|_{s}
\sim
\big\|\sup_{\hat h:\frac{D(\hat h)}{L}\le 1}\frac{W(\hat h)}{L}\big\|_{\frac{4}{3}s}^{\frac{4}{3}}.
\end{align}
\end{lemma}

For the readers' convenience, we display a proof of Lemma \ref{L:5} in the appendix.

\medskip

Here and in the following we will assume for simplicity that
$L\in2^{\mathbb{N}}$. Boiling down an explicit scale-by-scale martingale construction of an $h$ by Ding and Wirth 
\cite[Section 2.2 and 3]{dw}
from a geometrically nonlinear 2-dimensional setting, see Subsection \ref{SS:RFIM}, 
to our simpler (1+1)-dimensional setting (and in terms of the decomposition (\ref{ao67})
below), we obtain the lower bound

\begin{lemma}[{\cite[Section 2.2]{dw}}]\label{L:6}
\begin{align*}
\mathbb{E}\sup_{h:\frac{D(h)}{L}\le \ln
(e+L)}\frac{W(h)}{L}\gtrsim\ln L.
\end{align*}
\end{lemma}

Note that by (\ref{ao63}) and (\ref{ao81}) , Lemma \ref{L:6} yields in particular the lower bounds
\begin{align}\label{ao108}
\big\|\sup_{\hat h:\frac{D(\hat h)}{L}\le 1}\frac{W(\hat h)}{L}\big\|_2 
\gtrsim\ln^\frac{3}{4}L
\quad\mbox{and thus by Lemma \ref{L:5}}\quad
\|\frac{D(h_*)}{L}\|_{\frac{3}{2}} \gtrsim\ln L.
\end{align}
This shows that the reasoning that leads to (\ref{ao61}) was too naive.
As a side-effect of our main result Theorem \ref{T}, 
we obtain the matching upper bounds, see Corollary \ref{C:4}.
Our proof of Corollary \ref{C:4} is logically independent from
-- but related to in terms of the ingredients of concentration, covering/binning, and counting --
a general result by Talagrand \cite{Talagrandbook1,Talagrandbook2}.

\medskip

This relation has been discovered by Ding and Wirth in \cite[Section 2.4]{dw},
where Talagrand's reasoning is used to establish a lower bound on the correlation length 
in the random field Ising model
in terms of the field amplitude $\epsilon\ll 1$, see (\ref{ao91}).
They are led to study the related, geometrically nonlinear problem
\begin{align*}
	\sup_{P\subset B(0,L)}\frac{\xi(P)}{\text{Per}(P)}
	\stackrel{[12]}{\sim}\ln^{\frac{3}{4}}L ,
\end{align*}
where $\xi$ is the two-dimensional white noise and $P$ is a polygon of side length at least 1.
Previously, a similar quantity had appeared in the Leighton and Shor's grid matching
problem \cite[Theorem 2]{LeightonShor}, with the centered empirical measure of
independent points in place of the white noise, and had been
revisited by Talagrand via his chaining method, see \cite[Section 3.4]{Talagrandbook1}. The original
construction \cite[Theorem 1]{LeightonShor} hints on how to build a competitor that attains
the lower bound, which later on appeared in \cite{dw} (see also Lemma~\ref{L:6} above): 
In particular the inequality in the polygonal
decomposition in \cite[Theorem 1]{LeightonShor} is sharp if the triangles are
similar isosceles (as in \cite{dw} construction) and with height much smaller
than the basis. This flatness condition suggests to look at our linearized model, on which we
comment more in Section~\ref{SS:Heu}.

\medskip

In giving a proof independent on Talagrand's chaining, our paper  
is similar to \cite{LeightonShor} and to the arXiv version \cite[Section 5]{DingarXiv} of \cite{dw}.
While we work on the level of $\frac{D(h_*)}{L}$ in a scale-by-scale way, 
\cite{DingarXiv} establishes strict negativity of the 2-dimensional analogue
of $\mathbb{E}\sup_{h}(\epsilon\frac{W(h)}{L}-\frac{L+D(h)}{L})$ (see (\ref{ao90}) for
a motivation of $L+D$),
whereas Talagrand directly tackles $\sup_{\hat h:\frac{D(\hat h)}{L}\le 1}\frac{W(\hat h)}{L}$
in great generality. Both \cite{DingarXiv} and this paper differ from
Talagrand's approach by working with non-centered random variables.
For the readers' convenience, we also include a proof of Lemma \ref{L:6} 
in the appendix.


\subsection{Decomposition by scales}
Talagrand's approach to estimating the maximum over families of Gaussian variables
$W(h)$ over an index set $h$ is based on the canonical metric 
\begin{align}\label{ao52}
\sqrt{\mathbb{E}(\frac{W(\tilde h)}{L}-\frac{W(h)}{L})^2}
\stackrel{(\ref{ao65})}{=}\sqrt{\frac{1}{L}\sum_{x=1}^{L-1}|\frac{\tilde h(x)-h(x)}{L}|}.
\end{align}
Hence the problem comes with two metrics, the Hilbertian one given by
the Dirichlet energy (\ref{ao51}) and the non-affine one given by (\ref{ao52}).
We indirectly capitalize on this by a decomposition of $h$ that is orthogonal w.~r.~t.~$D$.

\medskip

More precisely, we decompose a configuration $h$ into
components of a given dyadic scale $l=1,2,\cdots,\frac{L}{4},\frac{L}{2}$. 
We implement this by introducing the projection operator
\begin{align}\label{ao66}
h\mapsto h_{\ge l}\quad\mbox{where}\; h_{\ge l}\;\mbox{is the piecewise linear
interpolation of $h$ in $x\in l\mathbb{Z}$}.
\end{align}
\begin{figure}
\centering
\includegraphics[width=0.45\linewidth]{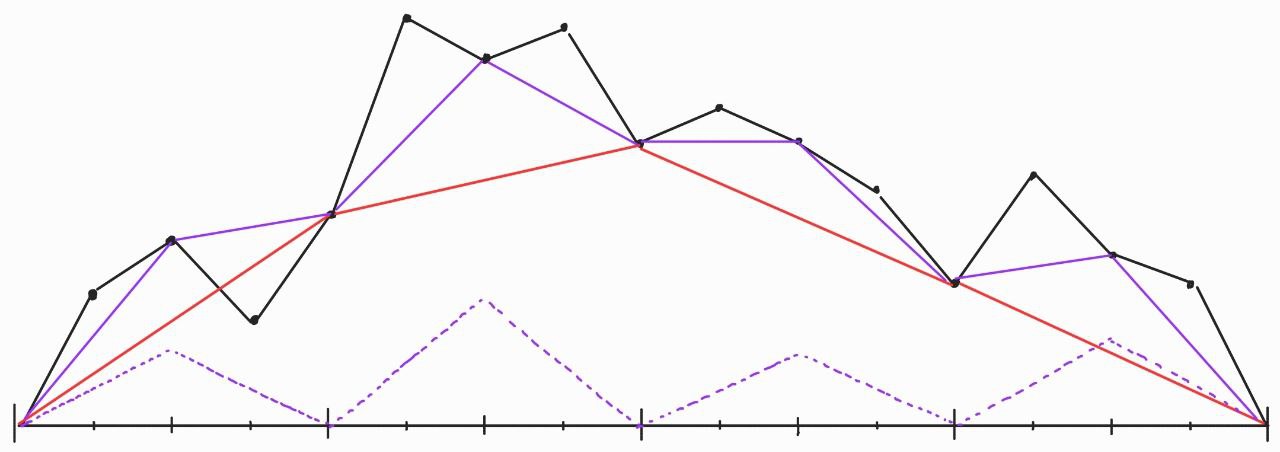}
\caption{Projection operator with $L = 16, l = 2$: $h$ in black, $h_{\ge l}$ in
violet, $h_{\ge 2l}$ in red, $h_l$ dotted in violet.}
\label{fig:f1}
\end{figure}
Since in particular $h_{\ge L}=0$ and $h_{\ge 1}=h$, this yields a decomposition: 
\begin{align}\label{ao67}
h_l:=h_{\ge l}-h_{\ge 2l}\quad\mbox{satisfies}\quad
h=\sum_{l} h_l.
\end{align}
An equivalent definition\footnote{used in the proof of Theorem \ref{T}} is obtained as follows:
Identify a configuration $h$ with its slope, the piecewise constant function $\frac{dh}{dx}$
given by the derivative of its piecewise linear interpolation; in view of the
boundary conditions in (\ref{ao49}), 
$\frac{dh}{dx}$ has spatial average zero and $h$ can be recovered
from it by integration. On this level, the definition (\ref{ao66}) \& (\ref{ao67}) amounts to
%
\begin{align}\label{ao116}
\frac{dh_l}{dx}\quad\mbox{}
&\mbox{is the $L^2(0,L)$-orthogonal projection of}\;\frac{dh}{dx}\;\mbox{on}\nonumber\\
&\mbox{the}\;\frac{L}{2l}-\mbox{dimensional span of the elements of the Haar 
basis of scale $2l$.}
\end{align}

\medskip

This decomposition is the (1+1)-dimensional version of a more adaptive
decomposition of polygons into polygonal curves in \cite[(44)]{DingarXiv}. 
Thanks to our geometrically linearized setting, 
we may in addition capitalize on (\ref{ao116}), of which we retain
\begin{align}\label{ao77}
h_{\ge l}=\sum_{\substack{\rho\;\text{dyadic} \\ l\le\rho<L}} h_\rho\quad\mbox{and}\quad
D(h_{\ge l})= \sum_{\substack{\rho\;\text{dyadic} \\ l\le\rho<L}}D(h_\rho).
\end{align}


\subsection{Main result} Our main result establishes that every scale $h_{*,l}$ contains (at most)
order one Dirichlet energy per length $\frac{D}{L}$, see (\ref{ao20}) for $p=2$.
In fact, we strengthen this control to $\frac{D_{p}}{L}$ with $ p < 3 $,
where the $p$-Dirichlet energy $D_p$ 
\begin{align}\label{ao106}
	D_p(h):=\frac{1}{2^{\frac{p}{2}}}\sum_{x=1}^L|h(x)-h(x-1)|^p\quad\mbox{so that}\quad D_2(h)=D(h)
\end{align}
captures improved spatial integrability of the (discrete) derivative. 
This strengthening to control of $D_p(h_*)$ for some $p>2$
will be crucial when informally justifying the
geometric linearization (\ref{ao90}) in Subsection \ref{SS:Heu}.
Note that as a consequence of Jensen's inequality
\begin{align*}
(\frac{D(h)}{L})^\frac{p}{2}\leq\frac{D_p(h)}{L}\quad\mbox{for}\;p\ge 2
\end{align*}
and thus as a consequence of the definition (\ref{ao24}) of $\|\cdot\|_s$:
\begin{align}\label{ao103}
\|\frac{D(h_{*,\ge l})}{L}\|_\frac{p}{2}^\frac{p}{2}
\le\|\frac{D_p(h_{*,\ge l})}{L}\|_1.
\end{align}
By (\ref{ao81}) it is natural to measure 
$\sup_{\hat h:\frac{D(\hat h)}{L}\le 1}\frac{W(\hat h)}{L}$ in the Gaussian norm $\|\cdot\|_2$;
hence in view of Lemma \ref{L:5}, it is natural to measure $\frac{D(h_*)}{L}$ 
in $\|\cdot\|_{\frac{3}{2}}$. 
Thus by (\ref{ao103}) we expect to measure $\frac{D_3(h_*)}{L}$ in $\|\cdot\|_1$. 
Up to an infinitesimal loss in integrability (both in space and probability), 
this is exactly what we obtain:

\begin{theorem}\label{T} 
For $2\le p<3$ we have\footnote{We add variables as subscripts below the symbol
$\lesssim$ (like $\lesssim_p$) if the implicit constant depends on
them.} for all $l=1,2,\cdots,\frac{L}{2}$
\begin{align}\label{ao20}
\big\|\frac{D_{p}(h_{*,l})}{L}\big\|_1\lesssim_p 1
\end{align}
and 
\begin{align}\label{ao80}
	\big\|\frac{D_p(h_{*,\ge l})}{L}\big\|_1\lesssim_p
	\ln^\frac{p}{2}\frac{L}{l}.
\end{align}
\end{theorem}

Note that by (\ref{ao103}), (\ref{ao80}) in particular implies
\begin{align*}
	\big\|\frac{D(h_{*,\ge
	l})}{L}\big\|_s\lesssim_s\ln\frac{L}{l}
\quad\mbox{for all}\;s<\frac{3}{2}.
\end{align*}
Theorem \ref{T} will be established by induction in the mesoscopic scale
$l=\frac{L}{2},\frac{L}{4},\cdots$, see Lemma \ref{L:3}, and thus from macroscopic to microscopic.
This order is reminiscent of the application of large-scale regularity theory, see \cite{AvellanedaLin},
in homogenization, as used for random environments in \cite{ArmstrongSmart}, 
also in a variational context.
This order is in contrast to the work \cite[Section 5.6]{DingarXiv} by Ding and Wirth,
where upper bounds on the geometrically nonlinear
analogue of $\mathbb{E}\sup_{h}(\epsilon\frac{W(h)}{L}-\frac{L+D(h)}{L})$ are established
by induction in the system size $L$, and thus from small to large scales.
We are able to run the induction in the opposite direction because of the 
macroscopic bound provided by Proposition \ref{P:1}, which is not available in the geometrically
nonlinear setting of \cite{DingarXiv}.
A benefit of our order is that we obtain incremental estimates like (\ref{ao20})
and estimates like (\ref{ao80}) that are oblivious to the small-scale cut-off 
given by the lattice spacing.
However, both inductions here and in \cite{DingarXiv}
can be seen as an induction from small to large complexity using a decomposition
by scales (multi-scale analysis).
The induction steps in both works rely on similar ingredients: 
concentration, binning, and counting.

\medskip

We close this section by drawing an easy conclusion from Theorem \ref{T}.
By (\ref{ao80}) for $p=2$ and $l=1$ we have $\|\frac{D(h_{*})}{L}\|_1$ $\lesssim
\ln L$;
by Lemma \ref{L:5} this translates to $\|\sup_{\hat h:\frac{D(\hat h)}{L}\le 1}
\frac{W(\hat h)}{L}\|_{\frac{4}{3}}^{\frac{4}{3}} $ $\lesssim\ln
L$. By Jensen's
inequality this yields $\mathbb{E}\sup_{\hat h:\frac{D(\hat h)}{L}\le 1}
\frac{W(\hat h)}{L}$ $\lesssim\ln^{\frac{3}{4}} L$, which by (\ref{ao81}) can be upgraded again to
$\|\sup_{\hat h:\frac{D(\hat h)}{L}\le 1}
\frac{W(\hat h)}{L}\|_{2}$ $\lesssim\ln^{\frac{3}{4}} L$. Appealing once more to Lemma \ref{L:5}
we learn

\begin{corollary}\label{C:4}
\begin{align}\label{ao82}
	\big\|\frac{D(h_{*})}{L}\big\|_\frac{3}{2}\lesssim\ln L\quad\mbox{next to}\quad
\big\|\sup_{\hat h:\frac{D(\hat h)}{L}\le 1}\frac{W(\hat h)}{L}\big\|_2
\lesssim\ln^\frac{3}{4} L.
\end{align}
\end{corollary}

We see Theorem \ref{T} 
as the main contribution of this paper: it provides finer scale-by-scale information
oblivious to the microscopic cut-off,
and thus in particular shows that the specific power 1 of the logarithm
appearing in the first item of (\ref{ao82}) is natural 
as it arises from counting dyadic scales; the power $\frac{3}{4}$
of the logarithm appearing in the second item of (\ref{ao82}) is then a plain consequence of
the scaling (\ref{ao54}).

\medskip

A superadditivity argument and the same techniques involved in the proof of
Theorem \ref{T} allow to prove the convergence of the rescaled minimal
energy to a deterministic constant, which amounts to a quenched
homogenization result.

\begin{theorem}\label{T:2}
	There exists a positive deterministic constant $\alpha$ such that
	\begin{align}
		&\lim_{L\to\infty}\frac{E(h_*)}{L\ln L}=-\alpha,\label{ao135}\\
		&\lim_{L\to\infty}\frac{D(h_*)}{L\ln
		L}= \frac{1}{3}  \alpha,\label{ao136}\\
		&\lim_{L\to\infty}\frac{1}{\ln^{\frac{3}{4}}L}\sup_{\frac{D(h)}{L}\le
		1}\frac{W(h)}{L}=\frac{4}{3^{\frac{3}{4}}}\alpha^{\frac{3}{4}}\label{ao137},
	\end{align}
	where the limits hold almost surely and with respect to the Orlicz norms
	$\|\cdot\|_{\frac{3}{2}}, \|\cdot\|_s $ (with $s<\frac{3}{2}$) and $\|\cdot\|_2$
	respectively.
\end{theorem}

The statement \eqref{ao137} can be seen as a first step in asymptotically characterizing the specific supremum of Gaussians.


\section{Motivation and Context}\label{S:Mot}


\subsection{The random field Ising model and its correlation length}\label{SS:RFIM}
Our motivation to study \eqref{ao48} comes from the random field Ising model,
which we introduce now. Given the $d$-dimensional lattice $\mathbb{Z}^d\ni z$
and the discrete spin space $\{-1,1\}\ni \sigma$, 
the configurations are given by $\mathbb{Z}^d\ni z\mapsto\sigma_z\in\{-1,1\}$.
The Hamiltonian is defined as
the sum of a perimeter term and a field term
\begin{equation}\label{i1}
	H(\sigma) = \sum_{|z-z'| = 1} |\sigma_z - \sigma_{z'}| + \epsilon \sum_{z} \xi_z
	\sigma_z,
\end{equation}
where $\epsilon>0$ is a parameter that we think of as being small, i.~e.~$\epsilon\ll 1$.
The field is given by independent Gaussians $(\xi_z)_{z\in\mathbb{Z}^d}$ 
of vanishing mean and unit variance.
In order for (\ref{i1}) to make sense, one considers a finite sub-lattice, 
e.~g.~$\mathbb{Z}^d\cap[-L,L]^d$,
and fixes boundary conditions outside, e.~g.~$\sigma_z=1$ for $z\not\in\mathbb{Z}^d\cap[-L,L]^d$,
so that both sums in (\ref{i1}) can be replaced by finite ones.

\medskip

The following Gedankenexperiment suggests that the critical dimension is $d=2$, 
the case on which we will focus: Flipping the spin in a regular domain of diameter $l$, 
say $\mathbb{Z}^d\cap[-l,l]^d$,
the perimeter term changes by an amount that scales like $l^{d-1}$, whereas the
field term changes by an amount that by the central limit theorem scales like $\sqrt{l^d}$; both expressions
are identical iff $d=2$.
It turns out that the critical dimension $d=2$ behaves like 
the lower dimension $d=1$ in the sense that
even at zero temperature, the case on which we will focus, 
there is no phase transition. By this one means that as $L\uparrow\infty$,
minimizers $\sigma_*$ will almost surely converge to a unique configuration. Hence the 
(quenched) randomness of the field suppresses the phase transition that is present
in $d=2$ even for (sufficiently small) positive temperature.
This prediction was first made by the physicists Imry and Ma \cite{IM75},
and later proven by Aizenman and Wehr \cite{AW90} by a purely qualitative argument
for any $\epsilon>0$. A more quantitative decay of spin-spin correlations 
$\mathbb{E}\sigma_z\sigma_0$ as $|z|\uparrow\infty$ has been
worked out in \cite{C19} (logarithmic decay), \cite{AP19} (power law),
and \cite{DingXia,AHP19} (exponential), again for any $\epsilon > 0$.

\medskip

A natural question is at which length scale, in terms of $\epsilon$, the influence of the boundary
conditions fades away: At which size $L$ of the system with boundary conditions $\equiv 1$
does the expectation $\mathbb{E}\sigma_*(0)$ of the spin $\sigma_*(0)$ at the origin $0$
drop below the value $\frac{1}{2}$?
This type of correlation length is expected to diverge as $\epsilon$ goes to zero;
at what rate does it diverge? This question has been addressed only more recently. 
In \cite{dw}, Ding and Wirth 
derive matching lower and upper bounds at zero temperature:
\begin{align}\label{ao91}
\ln(\mbox{correlation length})\sim\epsilon^{-4/3},
\end{align}
see \cite{BarNir} for suboptimal bounds for a more stringent
notion of correlation length.
Also the case of finite temperature has been tackled: The upper bound in (\ref{ao91})
was in fact established at all temperatures in \cite{dw}; the lower bound was later extended
to small temperatures in \cite{DingHuangXiaPert}. More recently, 
it was shown in \cite{DingHuangXia}
that (\ref{ao91}) holds up to the critical temperature.

\medskip

In \cite{rw} a continuum version of \eqref{i1} is considered: 
$\mathbb{Z}^2$ is replaced by $\mathbb{R}^2$ so that configurations $\sigma$ now are
given by $\mathbb{R}^2\ni z\mapsto\sigma_z \in \{-1, +1\}$;
likewise, the Hamiltonian is replaced by
\begin{equation}\label{i4}
H_{cont}(\sigma) = \mbox{perimeter of} \, \{\sigma=1\}+\frac{\epsilon}{2}\int\sigma d\xi
\end{equation}
where $\xi$ denotes white noise with a small-scale cut-off on scale one, 
for instance provided by a mollification. 
If the latter is radially symmetric,
the model (\ref{i4}) has the advantage over (\ref{i1}) of rotational invariance in law.
Incidentally, in \cite[Section 5.3]{DingarXiv}, the small-scale cut-off is introduced by
restricting $H_{\text{cont}}$ to $\sigma$'s such that the set $\{\sigma=1\}$ is a polygon with side length at least one.
Diffuse interface versions of
\eqref{i4} had already been considered in \cite{DirrOrlandi1,DirrOrlandi2}
extending the results by \cite{AW90}.
This continuum model allows to investigate a finer notion of correlation length,
namely the scale at which the interface ceases to be nearly flat.
Consider two points $z$ and $z'$ on the interface with $|z-z'|\gg1$ and 
let $\nu_z$ and $\nu_{z'}$ denote the corresponding normal.
Loosely speaking, it has been established in \cite{rw} that
\begin{equation*}
	|\nu_z - \nu_{z'}|^2 \lesssim \epsilon (\ln ( e + | z | ) + \ln ( e + |z'| ) )^{1/2}
(\ln|z-z'|)^{11/4}.
\end{equation*}
In particular,
\begin{equation}\label{i2}
\ln\max\{|z|,|z'|\}\ll \epsilon^{-4/13}\quad\Longrightarrow\quad
|\nu_z - \nu_{z'}| \ll 1,
\end{equation}
so that the logarithm of this finer correlation length is $\gtrsim \epsilon^{-4/13}$.


\subsection{Informal derivation of our (1+1) dimensional model}\label{SS:der}

\begin{figure}
\centering
\includegraphics[width=0.80\linewidth]{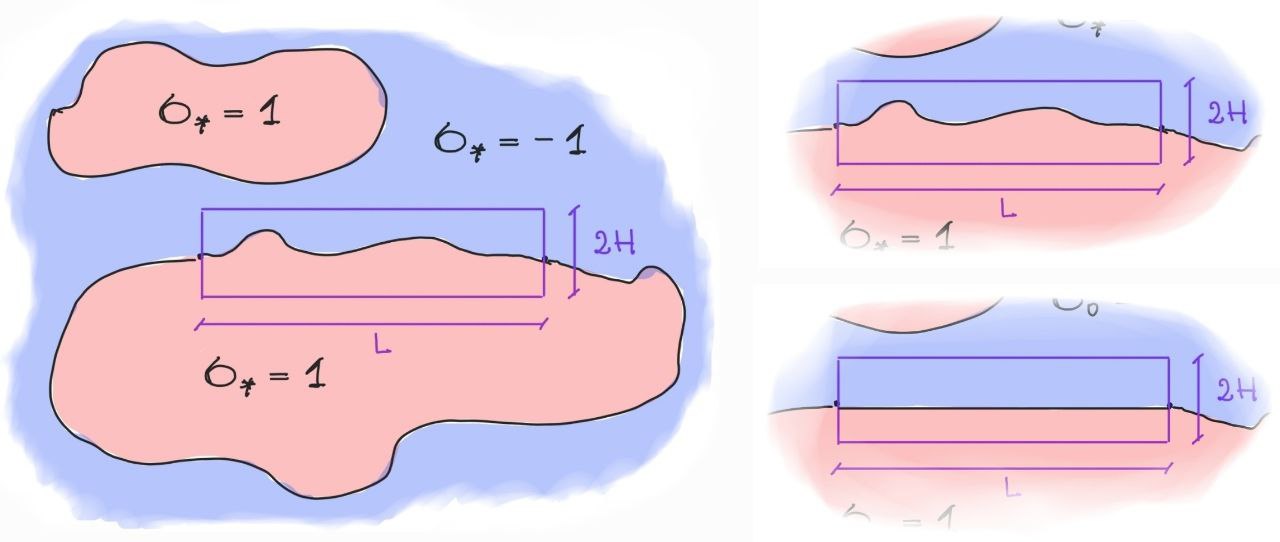}
\caption{Minimizing configuration $\sigma_*$ (left). Map from $\sigma_*$
	(up right) to the flat configuration $\sigma_0$ (down right) in the violet flat rectangle $[0,L]\times[-H,H]$.}
\label{fig:f3}
\end{figure}

We doubt that the exponent $\frac{4}{13}$ in (\ref{i2}), 
which is strictly smaller than the exponent $\frac{4}{3}$ in (\ref{ao91}), is optimal.
This motivates us to consider a reduced model, which we informally derive now,
in line with \cite[Section 6.2.2]{Peled}.
Consider a flat rectangle, which by rotational invariance we may take to be 
$[0,L]\times[-H,H]$ with $H\ll L$, in which the minimizer 
$\sigma_*$ is nearly flat (see Figure~\ref{fig:f3}), meaning that 
\begin{align*}
\{\sigma_*=1\}\cap([0,L]\times\{-H,H\})&=[0,L]\times\{-H\}\quad\mbox{and}\\
\{\sigma_*=1\}\cap(\{0,L\}\times[-H,H])&=\{0,L\}\times[-H,0].
\end{align*}
We compare $\sigma_*$ to the exactly flat configuration $\sigma_0$ given by 
\begin{align*}
\{\sigma_0=1\}\cap([0,L]\times[-H,H])=[0,L]\times[-H,0]
\end{align*}
and coinciding with $\sigma_*$ outside of $[0,L]\times[-H,H]$.
Then $\sigma_*$ minimizes $H_{cont}(\sigma)-H_{cont}(\sigma_0)$, 
which is well defined as a difference.
If we restrict to configurations $\sigma$ that are the subgraphs of a function
$[0,L]\ni x\mapsto h(x)\in[-H,H]$,  which in formulas means 
\begin{align*}
\{\sigma=1\}\cap([0,L]\times [ -H,H ] )=\{\,(x,y)\in[0,L]\times [ -H,H ] \,|\,y\le h(x)\,\},
\end{align*}
so that in particular $h(0)=h(L)=0$, then this difference of Hamiltonians can be written as
\begin{align}\label{ao94}
H_{cont}(\sigma)-H_{cont}(\sigma_0) 
=\int_0^L dx\big(\sqrt{1+(\frac{dh}{dx})^2}-1\big)
-\epsilon\int_0^Ldx\int_0^{h(x)}dy\xi(x, y).
\end{align}

\medskip

We now turn to the perimeter term in (\ref{ao94}). 
We note that the integrand $\sqrt{1+\zeta^2}-1$ as a function of
the slope $\zeta=\frac{dh}{dx}$
satisfies
\begin{align}\label{ao113}
\sqrt{1+\zeta^2}-1\approx\frac{1}{2}\zeta^2\quad\mbox{for}\;
|\zeta|\ll 1.
\end{align}
Hence if $\frac{dh}{dx}$ is small in an appropriate sense -- we will make this more
precise in Subsection \ref{SS:Heu} -- 
then the length term in (\ref{ao94}) is well-approximated by its 
geometric linearization\footnote{which amounts to a harmonic approximation}
\begin{align}\label{ao93}
\int_0^L dx(\sqrt{1 + (\frac{dh}{dx})^2}-1)\approx \frac{1}{2}\int_0^Ldx (\frac{dh}{dx})^2 ,
\end{align}
while the field term in (\ref{ao94}) can be rewritten as
\begin{align*}
\int_0^Ldx \int_0^{h(x)}dy \xi(x, y)=\int_0^Ldx W(x, h(x)) ,
\end{align*}
where $W(x,y)$ is the Gaussian field that arises from integrating our mollified white noise
in $y$, and thus acquires the character of two-sided Brownian motion in $y$ but
retains the short-range correlation in $x$. 
Incidentally, the approximation (\ref{ao93}) is present in 
\cite[Claim 5.20]{DingarXiv} in the form of a lower bound.
The last modification is to replace the small-scale cut-off provided by mollification of $\xi$ 
on scale one by a discretization of $[0,L]$ on scale one, that is, by
passing to the semi-discrete configuration space (\ref{ao49}) and to (\ref{ao50}):
\begin{equation}\label{i5}
\frac{1}{2}\sum_{x=1}^L(h(x)-h(x-1))^2-\epsilon\sum_{x=1}^LW(x,h(x))
\stackrel{(\ref{ao51}),(\ref{ao53})}{=}D(h) - \epsilon W(h).
\end{equation}
Finally, we exploit the scaling invariance (\ref{ao58}) \& (\ref{ao54}) in the specific form of
\begin{align}\label{ao95}
h=\epsilon^\frac{2}{3}\hat h\quad\mbox{and}\quad E=\epsilon^\frac{4}{3}\hat E,
\end{align}
where $\hat E(\hat h)$ $=D(\hat h)-W(\hat h)$ coincides with (\ref{ao48}). 
This is our motivation for studying (\ref{ao48}). 

\subsection{Heuristic justification of the geometric linearization}\label{SS:Heu}
We seek a heuristic identification of the regime in which the geometric linearization (\ref{ao93})
is legitimate. 
We start by noting that by the concavity of the square root we have
\begin{align*}
\sqrt{1+\zeta^2}-1\le\frac{1}{2} \zeta^2\quad\mbox{and thus}\quad
\int_0^Ldx(\sqrt{1+(\frac{dh}{dx})^2}-1)\le\frac{1}{2}\int_0^Ldx(\frac{dh}{dx})^2
\end{align*}
for any configuration $h$. Thanks to the $(p>2)$-version of Theorem \ref{T}, we
will show that the minimizer $h_*$ of the geometrically linearized 
problem (\ref{i5}) approximately satisfies the opposite 
inequality\footnote{we use the continuum notation which is exact if we interpolate
piecewise linearly}\textsuperscript{,} provided that
\begin{align}\label{ao110}
1 \ll \ln L\ll \epsilon^{-\frac{4}{3}},
\end{align}
in line with (\ref{ao91}).

\begin{proposition}\label{L:14}
If \eqref{ao110} holds then
\begin{align}
\int_0^Ldx(\sqrt{1+(\frac{dh_*}{dx})^2}-1)
= ( 1 - \eta ) \frac{1}{2}\int_0^Ldx(\frac{dh_*}{dx})^2 ,
\quad\mbox{with} ~ \eta \ll 1 ~ \mbox{in~probability}\footnotemark . \label{ao184}
\end{align}
\end{proposition}
\footnotetext{By this we mean that $ \eta \rightarrow 0 $ in probability as $ \ln L \rightarrow \infty $ and $ \epsilon^{ \frac{4}{3} } \ln L \rightarrow 0 $. Let us remark that with a bit more effort, the convergence can be quantified in some
Orlicz norm $\|\cdot\|_s$.}


\section{Strategy of proof}\label{S:strategy_of_proof}

The basis for Theorem \ref{T}, which characterizes the minimizer $h_*$ of
(\ref{ao48}) in terms of the Hilbertian metric given by the Dirichlet energy $D(h)$, 
is a bound in the other metric, cf.~Subsection \ref{SS:Tal}, namely
\begin{align}\label{ao79}
M(h):=\sum_{x=1}^{L-1}|h(x)|,
\end{align}
see (\ref{ao65}). We have

\begin{proposition}\label{P:1}
\begin{align*}
\big\|\frac{M(h_*)}{L^2}\big\|_3\lesssim 1.
\end{align*}
\end{proposition}
In establishing Proposition \ref{P:1}, we rely on two
ingredients: the characterization in law of the translation $E(\cdot+\phi)$ of
$E$, as crucially used in \cite[Proposition 2.1]{Peled} and the fluctuations of
the minimal energy proven in Lemma~\ref{L:4} below, see also \cite[Lemma 5.5]{Peled}.
However in our case this concentration of the minimal energy is not universal,
as the presence of $\bar M$ in $\eqref{c1}$ shows. We thus need a buckling argument.

\begin{lemma}\label{L:4}
\begin{align}
&\|E(h) - \mathbb{E}E(h)\|_2\lesssim \sqrt{M(h)},\label{ao117}\\
		& \|\min_{M(h) \leq \bar{M}} E(h) - \mathbb{E}\min_{M(h) \leq \bar{M}}
		E(h)\|_2 \lesssim \sqrt{\bar{M}} \qquad \text{for $\bar M<\infty$,}\label{c1}\\
		& \|\sup_{\frac{D(h)}{L} \leq 1} \frac{W(h)}{L} -
		\mathbb{E}\sup_{\frac{D(h)}{L} \leq 1} \frac{W(h)}{L}\|_2 \lesssim
		1\label{c2}.
	\end{align}
\end{lemma}

We will need multiple times the following standard equivalence of
stochastic norms, the proof of which can be found in the appendix.

\begin{lemma}\label{L:8}
For $s\ge1$, a random variable $X\ge 0$ and a constant $\nu_0<\infty$ we have
\begin{align}
& \forall\nu\ge\nu_0\;\;\ln\mathbb{P}(X\ge\nu)\le-\nu^s\quad\Longrightarrow\quad
\|X\|_{s}-\nu_0\lesssim 1;\label{ao39} \\
& \ln \mathbb{P}(X \geq
\nu)\lesssim-(\frac{\nu}{\|X\|_s})^s\qquad\mbox{for}\;\nu\gg\|X\|_s.\label{ao124}
\end{align}
\end{lemma}

\medskip

A simple post-processing of the proof of Proposition \ref{P:1} yields
the following pointwise-in-space bound:

\begin{corollary}\label{C:1} For all $x\in\{0,\cdots,L\}$
\begin{align}\label{ao21}
\big\|\frac{h_*(x)}{L}\big\|_{3}\lesssim 1.
\end{align}
\end{corollary}

Stronger versions of Lemma \ref{L:4} and Corollary \ref{C:1} have been established in the subsequent work \cite[Theorem 1.2]{Peled2}: In this work, \eqref{c1} is shown without the restriction $ M \leq \bar M $, cf.~Lemma \ref{L:10} below. Furthermore, Corollary \ref{C:1}, which is an annealed $ L^{ \infty } $-bound, is improved to the stronger quenched bound
\begin{align*}
	\big\| \sup_{ x } | \frac{ h_* ( x ) }{ L } | \big\|_{ 3 } \lesssim 1
\end{align*}
including also lower tail bounds.

\medskip

In order to obtain Theorem \ref{T}, we will apply Corollary \ref{C:1} on subsystems.
To do so, we need to extend Corollary \ref{C:1} to non-vanishing boundary values: 
For $h_0,h_1\in\mathbb{R}$ we consider 
\begin{align*}
h_{h_0,h_1,*}\quad\mbox{minimizes}\;E=D-W\;\mbox{under all}
\;h\;\mbox{with}\;h(0)=h_0\;\mbox{and}\;h(L)=h_1
\end{align*}
and monitor the analogue of (\ref{ao21}) for $x = \frac{L}{2}$
\begin{align}\label{ao114}
\frac{h_{h_0,h_1,*}(\frac{L}{2})-\frac{1}{2}(h_0+h_1)}{L}.
\end{align}
To estimate a bin of boundary data $(h_0,h_1)$ simultaneously, 
we rely on the following comparison principle\footnote{called order preservation
in \cite[Section 7.5]{Peled} but not used there}:

\begin{lemma}\label{L:MP}
We have
\begin{align}\label{ao125}
h_{*}\left\{\begin{array}{ccc}\leq& h_{h_0,h_1, *}&\mbox{for all}\;h_0,h_1\ge0\\
\geq& h_{h_0,h_1, *}&\mbox{for
all}\;h_0,h_1\le0\end{array}\right\}\quad\mbox{almost surely},
\end{align}
with the understanding that the null set does not depend on
$(h_0, h_1)$.
%
\end{lemma}

Indeed, Lemma \ref{L:MP} allows us to obtain control of (\ref{ao114})
with a locally\footnote{on scale $L$} quenched\footnote{by which we mean
that the supremum is inside the stochastic norm} uniformly in $(h_0,h_1)$. Loosely speaking, the
comparison principle provides a substitute for Dudley's bound in \cite[Section
5]{DingarXiv}, \cite[Corollary 3.9]{Peled}. 

\begin{corollary}\label{C:2} For all $\bar h_0,\bar h_1\in\mathbb{R}$ we consider the random
variable\footnote{in case of non-uniqueness the supremum in the r.~h.~s.~of \eqref{ao30} is taken among all
possible minimizers with the considered boundary conditions.}
\begin{align}\label{ao30}
X(\bar h_0,\bar h_1):=\sup_{h_0\in [\bar h_0-\frac{L}{2},\bar h_0+\frac{L}{2}],
h_1\in [\bar h_1-\frac{L}{2},\bar h_1+\frac{L}{2}]}
\frac{|h_{h_0,h_1,*}(\frac{L}{2})-\frac{1}{2}(h_0+h_1)|}{L}.
\end{align}
Then we have
\begin{align}
&\mbox{the law of}\;\;X(\bar h_0,\bar h_1)\;\;\mbox{is independent of}\;\;(\bar h_0,\bar h_1),
\label{ao31}\\
&\|X(\bar h_0,\bar h_1)\|_3\lesssim 1.\label{ao28}
\end{align}
\end{corollary}

\medskip

We will use Corollary \ref{C:2} on subsystems of size $2l$,
relying on the following connection between (\ref{ao114}) and $h_l$:
The definition (\ref{ao66}) \& (\ref{ao67}) of $h_l$ means that
\begin{align}\label{ao112}
\lefteqn{h_l\quad\mbox{is the piecewise linear interpolation of}\quad}\nonumber\\
&\left\{\begin{array}{cl}
0&\mbox{at}\;x=2\hat xl\\
h((2\hat x-1)l)-\frac{1}{2}\big(h(2(\hat x-1)l)+h(2\hat xl)\big)
&\mbox{at}\;x=(2\hat x-1)l
\end{array}\right\}\;\;\mbox{for}\;\hat x=0,\cdots,\frac{L}{2l}.
\end{align}
As a consequence, its $p$-Dirichlet energy (\ref{ao106}) can be expressed as
\begin{align}\label{ao27}
	D_p(h_l)=\frac{2l}{2^{\frac{p}{2}}}\sum_{\hat x=1}^{\frac{L}{2l}-1}
\Big|\frac{h((2\hat x-1)l)-\frac{1}{2}\big(h(2(\hat x-1)l)+h(2\hat xl)\big)}{l}\Big|^p.
\end{align}
If we use (\ref{ao27}) for the (global) minimizer $h_*$, then\footnote{up to a factor of 2} 
the argument of $|\cdot|^p$ is an instance of the expression (\ref{ao114}),
with the random boundary values $h_0=h_*(2(\hat x-1)l)$, $h_1=h_*(2\hat xl)$,
and with $L$ replaced by $2l$.

\medskip

In order to overcome the randomness of the boundary values, we introduce binning.
Given a fully discrete configuration on scale $2l$
\begin{align}\label{ao32}
\bar h\colon 2l\{0,\cdots,\frac{L}{2l}\}\rightarrow 2l\mathbb{Z}
\quad\mbox{with}\;\bar h(0)=\bar h(L)=0
\end{align}
we restrict to the event of
\begin{align}\label{ao29}
h_*(2\hat x l)\in(\bar h(2\hat x l)-l,\bar h(2\hat x l)+l].
\end{align}
Then we have by definition (\ref{ao30})
\begin{align}\label{ao46}
\Big|\frac{h_*((2\hat x-1)l)-\frac{1}{2}\big(h_*(2(\hat x-1)l)+h_*(2\hat xl)\big)}{2l}\Big|
\le X_{[2(\hat x-1)l,2\hat xl]}\big(\bar h(2(\hat x-1)l),\bar h(2\hat xl)\big),
\end{align}
where the subscript $[2(\hat x-1)l,2\hat xl]$ is to indicate that we replaced
$[0,L]$ by that interval (and the vertical variable is rescaled according to the
size of the interval). Incidentally, this combination of binning and coarse graining is also used
in \cite{LeightonShor} for the matching problem, and is implemented via coarser grids.\footnote{In that work, a coarse graining step is performed on a polygon by keeping only one vertex
every two, in strict analogy with our map $h_{\ge l}\mapsto h_{ \ge 2l } $.}

\medskip

The next important observation is that for given deterministic $\bar h$ of the form (\ref{ao32}),
the random variables on the r.~h.~s.~of (\ref{ao46}) are
independent (next to being identically distributed) so that the sum enjoys
good concentration properties for a large number $\frac{L}{l}$ of subsystems:

\begin{lemma}\label{L:1} For deterministic 
$\bar h$ of the form (\ref{ao32}) we have for all $\nu\gg 1$
\begin{align}\label{ao34}
\ln\mathbb{P}\big(\,\mbox{(\ref{ao29}) holds and}\;\;\frac{D_3(h_{*,l})}{L}\ge\nu\,\big)
\lesssim-\frac{L}{l}\nu.
\end{align}
\end{lemma}

\medskip

Another consequence of the localization argument provided by Corollary \ref{C:2} 
is an annealed modulus of continuity; using an iterative argument
similar to \cite[Lemma 3.8]{Peled} we obtain

\begin{proposition}\label{P:2} For all $1\le s<3$ and $x,y\in\{0,\cdots,L\}$
\begin{align*}
	\|h_*(y)-h_*(x)\|_{s}\lesssim|y-x|\big(1+\ln^\frac{4}{3}\frac{L}{|y-x|}\big).
\end{align*}
\end{proposition}

\medskip

In view of the stretched exponential integrability, by
Kolmogorov this annealed modulus translates
into a (slightly worse) quenched modulus of continuity. Together with the boundary condition
$h_*(0)=h_*(L)=0$, this yields tightness of the law of the
rescaled $\hat h_*(\hat x)=\frac{h_*(L\hat x)}{L}$ in $C^0([0,1])$ as $L\uparrow\infty$.
In view of Theorem \ref{T} one might expect
that the optimal exponent of the logarithm in Proposition \ref{P:2}
is $\frac{1}{2}$. Indeed, if we worked with periodic boundary conditions and
thus could capitalize on stationarity, we would obtain the exponent
$\frac{1}{2}$, but only for algebraic moments $p < 3$, i.~e.
$\mathbb{E}^{\frac{1}{p}}|h_*(x)-h_*(y)|^p\lesssim_p|y-x|( 1 +
\ln^{\frac{1}{2}}\frac{L}{|x-y|})$.

\medskip


The next ingredient for the proof of Theorem \ref{T} is a counting argument for
the bins, i.~e.~the functions $\bar h$ of the form (\ref{ao32}). 
In line with Talagrand's chaining, we think of a bin as capturing independence
structures. Thus there is no loss in applying a union bound over all bins, which
requires the following counting argument.
To this purpose, we identify $\bar h$ with a piecewise linear function on intervals of size $2l$, 
which allows us to define $\bar h_\rho$ for as in (\ref{ao66}) \& (\ref{ao67}).
We restrict the family of Dirichlet energies $D(\bar h_\rho)$ when counting the $\bar h$'s:

\begin{lemma}\label{L:2} We have for $\hat D\gg 1$
\begin{align}\label{ao69}
\ln\#\{\,\bar h\;\mbox{of form}\;(\ref{ao32})\,|\,
\max_\rho(\frac{l}{\rho})^2\frac{D(\bar h_\rho)}{L}<\hat D\,\}\lesssim \frac{L}{l}\ln \hat D,
\end{align}
where here and in the sequel, $\rho$ ranges over intermediate dyadic scales 
\begin{align}\label{ao71}
\rho\in\{2l,\cdots,\frac{L}{4},\frac{L}{2}\}.
\end{align}
\end{lemma}

As its proof reveals, (\ref{ao69}) amounts to counting lattice points
in the intersection of cylinders indexed by $\rho$, within the standard Hilbert
space. Working with the intersection of cylinders instead of one ball allows to
avoid double logarithms.
There is nothing special about how the prefactor $(\frac{l}{\rho})^2$
depends on the scale ratio $\frac{l}{\rho}\le 1$:
It could be exponentially small in $\frac{l}{\rho}$ without affecting
the scaling (\ref{ao69}), it needs to be logarithmically small for Lemma \ref{L:3}. 

\medskip


In order to combine 
the concentration of $D_3(h_*)$ within a bin indexed by $\bar h$, see Lemma \ref{L:1}, 
with the counting of $\bar h$'s in terms of $D(\bar h_\rho)$, see Lemma \ref{L:2},
to a statement on $D_3(h_*)$ and $D(h_*)$, we need to relate $D(\bar h_\rho)$ to $D(h_{*,\rho})$. 
Recall that $\bar h$ is related to $h_*$ via (\ref{ao29}) which is made such that
\begin{align*}
|\bar h(2l\hat x)-h_*(2l\hat x)|\le l
\quad\mbox{for}\;\hat x\in0,\cdots,\frac{L}{2l},
\end{align*}
By the linearity of the projection $h\mapsto h_\rho$ and in view of (\ref{ao112}) this yields
\begin{align*}
|\bar h_{\rho}((2\hat x-1)\rho)-h_{*,\rho}((2\hat x-1)\rho)|\le 2l
\quad\mbox{for}\;\hat x\in0,\cdots,\frac{L}{2\rho}
\end{align*}
for all $\rho$ satisfying (\ref{ao71}).
From the representation \eqref{ao112} \&
\eqref{ao27} of $D$ (for $p=2$) 
we therefore obtain by the triangle inequality
\begin{align*}
|\sqrt{\frac{D(\bar h_{\rho})}{L}}-\sqrt{\frac{D(h_{*,\rho})}{L}}|\le 1
\quad\mbox{and thus}\quad
\frac{D(\bar h_{\rho})}{L}\le 2(\frac{D(h_{*,\rho})}{L}+1),
\end{align*}
so that for $\hat D\ge 1$
\begin{align}\label{ao38}
\max_\rho(\frac{l}{\rho})^2\frac{D(h_{*,\rho})}{L}\le\hat D
\quad\Longrightarrow\quad
\max_\rho(\frac{l}{\rho})^2\frac{D(\bar h_{\rho})}{L}\le4\hat D.
\end{align}
This allows us to combine Lemmas \ref{L:1} and \ref{L:2} to

\begin{corollary}\label{C:3} Provided $\hat D\gg 1$ we have
\begin{align}\label{ao37}
\big\|I\big(\max_\rho(\frac{l}{\rho})^2\frac{D(h_{*,\rho})}{L}\le \hat D\big) 
\frac{D_3(h_{*,l})}{L}\big\|_{1}
\lesssim\ln\hat D,
\end{align}
where $\rho$ ranges over (\ref{ao71}).
\end{corollary}

\medskip


Corollary \ref{C:3} in turn allows us to bound $D_{3}(h_{*,l})$, with a slight
loss of spatial integrability to an exponent $p<3$. We have to give up spatial integrability 
in order to gain the stochastic integrability needed to cope with the
product of the indicator function $I$ and of $D_3(h_{*,l})$ in (\ref{ao37}). This
exchange of spatial and stochastic integrability relies on the observation that by
definition (\ref{ao106}) of $D_p$ and Jensen's inequality we have
\begin{align*}
(\frac{D_p(h)}{L})^\frac{1}{p}\quad\mbox{is monotone in}\quad p,
\end{align*}
so that by definition (\ref{ao24}) of $\|\cdot\|_s$
\begin{align}\label{ao99}
\big\|I\frac{D_p(h_{*,l})}{L}\big\|_\frac{3}{p}^\frac{3}{p}
\le\big\|I\frac{D_3(h_{*,l})}{L}\big\|_1\quad\mbox{for any $\{0,1\}$-valued random variable $I$}
\end{align}
cf.~(\ref{ao103}).

\begin{lemma}\label{L:3} For any $1\le p<3$ 
\begin{align}\label{ao40}
\big\|\frac{D_p(h_{*,l})}{L}\big\|_1\lesssim_{p} 
\ln^\frac{p}{3}(e+\max_\rho\big\|\frac{D(h_{*,\rho})}{L} \big\|_1) ,
\end{align}
where $\rho$ ranges over (\ref{ao71}).
\end{lemma}

While (\ref{ao20}) in Theorem \ref{T} is an easy consequence of Lemma \ref{L:3},
the aggregated (\ref{ao80}) requires a bit of harmonic analysis in case
of $p\not=2$; for $p=2$ it follows immediately from (\ref{ao77}).

\medskip

The main step in the proof of Theorem \ref{T:2} is establishing a
superadditivity property, up to lower order terms,
for the function\footnote{Recall that we are restricting to dyadic scale $ L $.}
\begin{align*}
	\ln L\mapsto-\frac{\mathbb{E}\min E}{L}.
\end{align*}
This is obtained by constructing a competitor at scale $L$ in terms of the
minimizers of the problems at scale $l$ and $\frac{L}{l}$. The idea is a
simple two-scale construction:
\begin{itemize}
	\item consider the minimizer of $E$ constraint to coarse-grained
		profiles, which are piecewise linear with interpolation points
		$[0,L]\cap l\mathbb{Z}$.
	\item on each subsystem of the form $[(\hat x-1)l, \hat xl]$ paste in the appropriately sheared
	minimizer on scale $l$ and microscopic scale 1.
\end{itemize}
In performing the first step, we construct Brownian motions
$\hat W(\hat x, \cdot)$ so that
\begin{align*}
\sum_{x=l(\hat x-1)+1}^{l\hat x}W(x,h(x))\approx l\hat W(\hat x, lh(\hat xl))
\qquad\mbox{for}~\hat x=1,\ldots,\frac{L}{l}.
\end{align*}
We use the Brownian motions on the new lattice, to construct an approximating
energy. The new problem has the advantage of having the same law (up to
rescaling) of the problem in \eqref{ao48} on the lattice $\hat
x\in\{0,\ldots,\frac{L}{l}\}$. 

\medskip

\begin{lemma}\label{L:9}
The limit
\begin{align*}
\lim_{L\to\infty}\frac{\mathbb{E}\min E}{L\ln L}\in[-\infty,0]
\end{align*}
exists.
\end{lemma}

Lemma~\ref{L:9} is indeed enough to prove the limit in \eqref{ao135} when
combined with the following concentration property\footnote{Incidentally, Lemma~\ref{L:10} is also used in the proof of
Lemma~\ref{L:9}.}.

\begin{lemma}\label{L:10}
\begin{align}\label{ao138}
\|\min E-\mathbb{E}\min E\|_{\frac{3}{2}}\lesssim L.
\end{align}
\end{lemma}

The proof of Lemma~\ref{L:10} is just a combination of the constrained
concentration in \eqref{c1} and the bound given by Proposition~\ref{P:1}.
The same concentration result as \eqref{ao138} is obtained by \cite[Theorem 1.2]{Peled2},
also proving lower tail bounds.

\medskip

In the proofs of Lemma~\ref{L:9}, one needs multiple times to provide upper bounds on the sum of
$\frac{L}{l}$ local errors $X_{\hat x}(\bar h)$ (which we further specify
in the proof) of order $1$, where $\bar h$
is a coarse-grained profile of the form
\begin{align}\label{ao170}
	\mbox{$\bar h$ piecewise linear with interpolation points $\bar h(\hat
	xl)\in l\mathbb{Z}$ for $\hat x\in\mathbb{Z}$.}
\end{align}
For a fixed $\bar h$, each variable $X_{\hat x}(\bar h)$ depends only on the noise in the strip $[(\hat
x-1)l,\hat xl]\times\mathbb{R}$ so that one can capitalize on independence
to get a concentration result in the spirit of Lemma~\ref{L:1}. A counting argument like
the one in Lemma~\ref{L:2} then allows to estimate the number of profiles $\bar h$. 
The results we need are gathered in the following lemma.

\begin{lemma}\label{L:12}
Assume that a set of random variables 
\begin{align*}
X_x(\bar h)\quad\mbox{for}\;
x=1,\ldots,L\;\mbox{and}\;\bar h:\{1,\ldots,L\}\to\mathbb{Z}
\end{align*}
is given such that for any $\bar h$
\begin{align}\label{ao189}
\mbox{$\{X_x(\bar h)\}_x$ are independent.}
\end{align}
If for some $s,p\in[1,\infty]$ with $\frac{1}{p}+\frac{1}{s}\le 1$ and
$\alpha\in[0,\frac{3}{2})$, we have\footnote{with the usual convention that for
$p=\infty$ we replace the l.h.s.~by $\sup_x\|X_x(\bar h)\|_s$.}
\begin{align}\label{f12}
\Big(\frac{1}{L}\sum_{x=1}^L\|X_x(\bar h)\|_s^p\Big)^{\frac{1}{p}}\le \big(\frac{D(\bar h)}{L}\big)^\alpha
\qquad\mbox{for every $\bar h$,}
\end{align}
then
\begin{align*}
\|\frac{1}{L}\sum_{x=1}^LX_x(\bar h_*)\|_t\lesssim\ln^\alpha L\qquad
\mbox{provided $\displaystyle\frac{1}{t}:=\frac{1}{s}+\frac{2}{3}\alpha$ and $t\ge 1$.}
\end{align*}
\end{lemma}

Finally, using the equivalent formulation of the problem in \eqref{ao56} and
exploiting the fluctuation bound in \eqref{c2}, the next lemma establishes a
concentration of the Dirichlet energy of the minimizer $h_*$ and shows the
equivalence of the existence of the two limits \eqref{ao136} and \eqref{ao137}.

\begin{lemma}\label{L:13}
	\begin{align}
		&\|\frac{D(h_*)}{L}-(\frac{1}{4}\mathbb{E}\sup_{\frac{D(h)}{L}\le
		1}\frac{W(h)}{L})^{\frac{4}{3}}\|_1\lesssim
		(\ln^{\frac{5}{8}}L)\ln^{\frac{1}{4}}\ln L, \label{ao171}\\
		&\|\frac{W(h_*)}{L}-(\frac{1}{4})^{\frac{1}{3}}
		(\mathbb{E}\sup_{\frac{D(h)}{L}\le 1}\frac{W(h)}{L})^{\frac{4}{3}}\|_1
	\lesssim(\ln^{\frac{5}{8}}L)\ln^{\frac{1}{4}}\ln L.\label{ao200}
	\end{align}
\end{lemma}

Lemma \ref{L:13}, together with the existence of the limit \eqref{ao135} of the
minimal energy, allows to prove the existence of limits \eqref{ao136} and \eqref{ao137} in
Theorem \ref{T:2}. Moreover, since we have upper and lower bounds for the
Dirichlet energy (cf.~\eqref{ao82} and \eqref{ao108}) and an explicit relation
between the three limits in Theorem \ref{T:2}, then they are not zero
nor infinity.

\section{Detailed proofs}

{\sc Proof of Proposition \ref{L:14}}.
We can assume that $h_*\neq0$, since otherwise the statement holds 
with $\eta=0$. Recall the definition of $\eta$ in the first item of
\eqref{ao184}. To show that $\eta$ vanishes, we will first argue that
for any deterministic parameter $ \nu $
\begin{align}\label{ao214}
\eta\le\frac{ \nu^2+\tilde{\eta} }{1+\nu^2} ,
\end{align}
where $ \tilde{\eta} $ is defined by
\begin{align}\label{ao217}
\frac{ 1 }{ 2 } \int_0^LdxI(|\frac{dh_*}{dx}|>\nu)(\frac{dh_*}{dx})^2=
\tilde{\eta}\frac{1}{2}\int_0^Ldx(\frac{dh_*}{dx})^2.
\end{align}
Using \eqref{ao214}, we we will conclude by first showing that for any positive $\nu$
\begin{align}\label{ao216}
\mbox{$\tilde{\eta}\to0$ in probability as $\ln L\to\infty$ and
$\epsilon^{\frac{4}{3}}\ln L\to 0$}
\end{align}
and then sending $\nu$ to zero.

\medskip

In order to get \eqref{ao214}, let us split the Dirichlet energy using $\nu$ and
appeal to $ \frac{ 1 }{ 2 } \zeta^2 I ( | \zeta | < \nu ) \le(1+\nu^2)(\sqrt{1+\zeta^2}-1)$
\begin{align*}
\frac{ 1 }{ 2 } \int_0^Ldx(\frac{dh_*}{dx})^2
&=\frac{ 1 }{ 2 }  \int_0^LdxI(|\frac{dh_*}{dx}|>\nu)(\frac{dh_*}{dx})^2+
\frac{ 1 }{ 2 }  \int_0^LdxI(|\frac{dh_*}{dx}|\le\nu)(\frac{dh_*}{dx})^2\\
&\stackrel{ \eqref{ao217} }{ \le } \tilde{\eta} \frac{ 1 }{ 2 }  \int_0^Ldx(\frac{dh_*}{dx})^2+(1+\nu^2)
\int_0^Ldx(\sqrt{1+(\frac{dh_*}{dx})^2}-1)\\
&\stackrel{ \eqref{ao184} }{ = } \tilde{\eta} \frac{ 1 }{ 2 }  \int_0^Ldx(\frac{dh_*}{dx})^2+
(1+\nu^2)(1-\eta) \frac{ 1 }{ 2 }  \int_0^Ldx(\frac{dh_*}{dx})^2.
\end{align*}
Cancelling $\frac{1}{2}\int_0^Ldx(\frac{dh_*}{dx} )^2$ and reordering the terms, this proves \eqref{ao214}.

\medskip

Here comes the proof of \eqref{ao216}. For any $p\ge2$ let us show
\begin{align}\label{ao215}
\tilde{\eta}\lesssim_p\big(\frac{\epsilon^\frac{2}{3}}{\nu}\big)^{p-2}\frac{D_p(\hat
h_*)}{D(\hat h_*)}\quad\mbox{where~we~recall~that}\quad\hat h_*=\epsilon^{-\frac{2}{3}}h_*
\end{align}
is the reduced variable. This comes from an application of H\"older's and Chebyshev's inequality
\begin{align*}
&\tilde{\eta} \frac{1}{2L}\int_0^Ldx(\frac{d\hat h_*}{dx})^2=
\frac{1}{2L}\int_0^Ldx I(|\frac{d\hat
h_*}{dx}|>\epsilon^{-\frac{2}{3}}\nu)(\frac{d\hat h_*}{dx})^2\\
&\le\big(\frac{{\mathcal L}(\{|\frac{d\hat h_*}{dx}|>\epsilon^{-\frac{2}{3}} \nu \})}{L}\big)^{1-\frac{2}{p}}
(\frac{D_p(\hat h_*)}{L})^\frac{2}{p}
\lesssim_{ p } \big(\frac{\epsilon^\frac{2}{3}}{\nu}\big)^{p-2}\frac{D_p(\hat h_*)}{L},
\end{align*}
where ${\mathcal L}$ denotes the one-dimensional Lebesgue measure.

\medskip

Inequality \eqref{ao215} can be rewritten as
\begin{align*}
&\tilde\eta\lesssim_p(\frac{\epsilon^{\frac{4}{3}}\ln
L}{\nu^2})^{\frac{p}{2}-1}(\frac{D_p(\hat h_*)}{L\ln^{\frac{p}{2}}L})(\frac{L\ln
L}{D(\hat h_*)})\quad\mbox{and thus}\\
&\tilde\eta\lesssim_p(\frac{\epsilon^{\frac{4}{3}}\ln
L}{\nu^2})^{\frac{p}{2}-1}(\frac{D_p(\hat h_*)}{L\ln^{\frac{p}{2}}L})
+I(\frac{D(\hat h_*)}{L\ln L}\le\frac{1}{2}\frac{4^{\frac{1}{3}}}{3}\alpha)
\end{align*}
where we used that $\tilde\eta$ is always less then 1 (cf. \eqref{ao217}). From
\eqref{ao80} in Theorem~\ref{T} applied with $l=1$, the first term converges to
zero in $\|\cdot\|_1$ (and thus in probability) for any $p\in(2,3)$ provided $\epsilon^{\frac{4}{3}}\ln
L$ vanishes. Finally, from the limit \eqref{ao136} in Theorem~\ref{T:2}, the indicator
function tends to zero in probability as $L$ diverges, so that \eqref{ao216} is established.
\qed

\medskip

{\sc Proof of Lemma \ref{L:4}}. Note that the Cameron-Martin space of (\ref{ao50}) is given by the direct sum
\begin{align*}
\bigoplus_{x=1}^{L-1}\{\,W\in \dot H^1(\mathbb{R})\,|\,W(0)=0\,\},
\end{align*}
where $\dot H^1(\mathbb{R})$ 
is to indicate that $\sqrt{\int_\mathbb{R}dh(\frac{dW(h)}{dh})^2}$ 
is the Hilbert norm of every component. In view of\footnote{which
is a consequence of Cauchy-Schwarz in $h$}
\begin{align*}
\big|W(x,h(x))\big|
\le\sqrt{|h(x)|}\sqrt{\int_\mathbb{R}dh(\frac{dW(x,h)}{dh})^2},
\end{align*}
$W\mapsto E(h)=D(h)-W(h)$ is Lipschitz continuous with
Lipschitz constant $\sqrt{\sum_{x=1}^{L-1}\sqrt{|h(x)|}^2}$ $=\sqrt{M(h)}$ in this geometry. 
As a consequence, also
\begin{align*}
W\mapsto\min_{M\le\bar M}E:=\min_{h:M(h)\le\bar M}E(h)
\end{align*}
is Lipschitz continuous with Lipschitz constant 
$\max\{\sqrt{M(h)}|M(h)\le\bar M\}$ $=\sqrt{\bar M}$, cf.~(\ref{ao79}). 
These Lipschitz bounds imply Gaussian concentration (see \cite[Theorem
4.5.7]{BogachevGM}), that is, for all $\nu\ge 0$
\begin{align*}
	\mathbb{P}\big(\min_{M\le\bar M}E-\mathbb{E}\min_{M\le\bar
	M}E\ge\sqrt{\bar M}\nu\big)
&\le\exp(-\frac{\nu^2}{2}),
\end{align*}
Appealing to Lemma~\ref{L:8} with $s=2$ and $\nu_0=0$, this implies
\begin{align*}
\big\|\frac{\max\{\min_{M \leq \bar{M}} E - \mathbb{E}\min_{M \leq \bar{M}}E,0\}}{\sqrt{\bar{M}}}
\big\|_2\lesssim 1 .
\end{align*}
The same argument applies to the negative part, which yields (\ref{c1})
by the triangle inequality.
The same reasoning implies (\ref{ao117}), 
and together with the Poincaré inequality in \eqref{ao65}, which bounds the
Lipschitz constant of $\frac{W(h)}{L}$ for all the $h$ under consideration, it also yields \eqref{c2}.
\qed

\medskip

{\sc Proof of Proposition \ref{P:1}}. It is convenient to introduce the Green function $\phi$ of the discrete
(one-dimensional) Laplacian $-\frac{d^2}{dx^2}$ arising from the Dirichlet energy $D$.
In terms of the related Dirichlet form 
\begin{align}\label{ao05}
D(h,\phi)=\sum_{x=1}^{L}(h(x)-h(x-1))(\phi(x)-\phi(x-1)),
\end{align}
and for every site $y\in\{1,\cdots,L-1\}$
the Green function is characterized by the representation
\begin{align}\label{ao01}
D(h,\phi_y)=h(y).
\end{align}
It is given by
\begin{align}\label{ao09}
\phi_y(x)=\min\{\frac{(L-y)x}{L},\frac{y(L-x)}{L}\}.
\end{align}
For later purposes we note
\begin{align}
M(\phi_y)\stackrel{(\ref{ao79}),(\ref{ao09})}{=}\frac{L}{2}\phi_y(y)
\stackrel{(\ref{ao09})}{=}\frac{(L-y)y}{2}
\le\frac{L^2}{8},\label{ao10}\\
D(\phi_y)
\stackrel{(\ref{ao51}),(\ref{ao05})}{=}\frac{1}{2}D(\phi_y,\phi_y)
\stackrel{(\ref{ao01})}{=}\frac{1}{2}\phi_y(y)
\stackrel{(\ref{ao09})}{=}\frac{(L-y)y}{2L}
\le\frac{L}{8}.\label{ao11}
\end{align}

\medskip

Fix a deterministic $\bar M\ge0$. We now derive the following inequality for every site $y=1,\cdots,L-1$ and $\lambda\in\mathbb{R}$
\begin{align}\label{ao02}
\lambda h_*(y)&\le(\min_{M\le\bar M}E-\mathbb{E}\min_{M\le\bar M}E)
-(\min_{M\le\bar M}E_{y, \lambda}-\mathbb{E}\min_{M\le\bar M}E_{y, \lambda})\nonumber\\
& - (E(\lambda\phi_y)-\mathbb{E}E(\lambda\phi_y)) +\frac{\lambda^2
L}{8}
\end{align}
provided
\begin{align}\label{ao02cont}
M(h_*)+\frac{|\lambda|L^2}{8}\le\bar M.
\end{align}
Here the functional $E_{y, \lambda}$ is defined through
\begin{align}\label{ao04}
	E_{y,\lambda}(h):= D(h) - W(h+\lambda\phi_y) + W(\lambda\phi_y) , 
\end{align}
where by the stationary increments of Brownian motions
$\{W(x,\cdot)\}_x$,
\begin{align}\label{ao16}
	E_{y, \lambda}\quad\mbox{has the same law as}\quad E
\end{align}
so that in particular
\begin{align}\label{ao07}
	\mathbb{E}\min_{M\le\bar M}E_{y,\lambda}=\mathbb{E}\min_{M\le\bar M}E.
\end{align}
Due to \eqref{ao02cont}
\begin{align*}
M(h_*-\lambda\phi_y)\stackrel{(\ref{ao79})}{\le}
M(h_*)+|\lambda|M(\phi_y)\stackrel{(\ref{ao10})}{\le}M(h_*)+\frac{|\lambda|L^2}{8}\le\bar M,
\end{align*}
which ensures
\begin{align}\label{ao08}
E(h_*)=\min_{M\le\bar M}E\quad\mbox{and}\quad
E_{y,\lambda}(h_*-\lambda\phi_y)\ge\min_{M\le\bar M}E_{y,\lambda}.
\end{align}

\medskip

By the centeredness of $W$ we have $\mathbb{E}E(\lambda\phi_y)$ $=D(\lambda\phi_y)$ 
so that by (\ref{ao11}) we obtain
\begin{align*}
 D (\lambda\phi_y) + W ( \lambda \phi_y ) \le - (E(\lambda\phi_y)-\mathbb{E}E(\lambda\phi_y)) + \frac{\lambda^2 L}{8}.
\end{align*}
As a consequence of (\ref{ao01}), (\ref{ao07}), and (\ref{ao08}),
the inequality (\ref{ao02}) thus reduces to the identity
\begin{align*}
\lambda D(h_*,\phi_y)=E(h_*)-E_{y,\lambda}(h_*-\lambda\phi_y)+ D (\lambda\phi_y) + W ( \lambda \phi_y ) ,
\end{align*}
which by the definitions (\ref{ao48}) and (\ref{ao04}) amounts to
\begin{align*}
\lambda D(h_*,\phi_y)=D(h_*)-D(h_*-\lambda\phi_y)+D(\lambda\phi_y)
\end{align*}
that is the polarization identity that relates the bilinear form
(\ref{ao05}) to its quadratic part (\ref{ao51}).

\medskip

Restricting to $\lambda\ge 0$ but then 
using (\ref{ao02}) for both $\pm\lambda$, we obtain
\begin{align}\label{ao89}
\lefteqn{\lambda|h_*(y)|}\nonumber\\
&\le|\min_{M\le\bar M}E-\mathbb{E}\min_{M\le\bar M}E|
+|\min_{M\le\bar M}E_{y, \lambda}-\mathbb{E}\min_{M\le\bar M}E_{y, \lambda}|
+|\min_{M\le\bar M}E_{y, -\lambda}-\mathbb{E}\min_{M\le\bar M}E_{y, -\lambda}|\nonumber\\ 
&+|E(\lambda\phi_y)-\mathbb{E}E(\lambda\phi_y)|
+|E(-\lambda\phi_y)-\mathbb{E}E(-\lambda\phi_y)|
+\frac{\lambda^2L}{8}\nonumber\\
&\mbox{}\quad\;\;\mbox{provided}\quad M(h_*)+\frac{\lambda L^2}{8}\le\bar M.
\end{align}
Summing over $y=1,\cdots,L-1$, this implies by definition (\ref{ao79})
\begin{align*}
\lefteqn{\lambda M(h_*)\le \frac{\lambda^2L^2}{8}}\nonumber\\
&+ L|\min_{M\le\bar M}E-\mathbb{E}\min_{M\le\bar M}E|
+\sum_{y=1}^{L-1}\big(|\min_{M\le\bar M}E_{y, \lambda}-\mathbb{E}\min_{M\le\bar
	M}E_{y, \lambda}|
+|\min_{M\le\bar M}E_{y, -\lambda}-\mathbb{E}\min_{M\le\bar
	M}E_{y, -\lambda}|\big)\nonumber\\
&+\sum_{y=1}^{L-1}\big(|E(\lambda\phi_y)-\mathbb{E}E(\lambda\phi_y)|
+|E(-\lambda \phi_y)-\mathbb{E}E(-\lambda\phi_y)|\big)
\mbox{}\quad\;\;\mbox{provided}\quad M(h_*)+\frac{\lambda L^2}{8}\le\bar M.
\end{align*}
This allows us to control the event that 
$M(h^*)\in[\frac{\bar M}{4},\frac{\bar M}{2}]$; yielding at first
\begin{align*}
\lefteqn{\frac{\lambda\bar M}{4}I\big(M(h^*)\in[\frac{\bar M}{4},\frac{\bar M}{2}]\big)
\le \frac{\lambda^2L^2}{8}+L|\min_{M\le\bar M}E-\mathbb{E}\min_{M\le\bar M}E|}\nonumber\\
&+\sum_{y=1}^{L-1}\big(|\min_{M\le\bar M}E_{y, \lambda}-\mathbb{E}\min_{M\le\bar
	M}E_{y, \lambda}| + |\min_{M\le\bar M}E_{y, -\lambda}-\mathbb{E}\min_{M\le\bar
	M}E_{y, -\lambda}|\big) \nonumber\\
& +\sum_{y=1}^{L-1}\big(|E(\lambda\phi_y)-\mathbb{E}E(\lambda\phi_y)|
+|E(-\lambda\phi_y)-\mathbb{E}E(-\lambda\phi_y)|\big)\mbox{}\;\;\quad\mbox{provided}\quad\bar
M\ge\frac{\lambda L^2}{4}\nonumber 
\end{align*}
and then
\begin{align*}
\lefteqn{\frac{\lambda\bar M}{8}I\big(M(h^*)\in[\frac{\bar M}{4},\frac{\bar M}{2}]\big)
\le L|\min_{M\le\bar M}E-\mathbb{E}\min_{M\le\bar M}E|}\nonumber\\
&+\sum_{y=1}^{L-1}\big(|\min_{M\le\bar M}E_{y, \lambda}-\mathbb{E}\min_{M\le\bar
	M}E_{y, \lambda}| + |\min_{M\le\bar M}E_{y, -\lambda}-\mathbb{E}\min_{M\le\bar
	M}E_{y, -\lambda}|\big) \nonumber\\
& +\sum_{y=1}^{L-1}\big(|E(\lambda\phi_y)-\mathbb{E}E(\lambda\phi_y)|
+|E(-\lambda\phi_y)-\mathbb{E}E(-\lambda\phi_y)|\big)\mbox{}\;\;\quad\mbox{provided}\quad\bar
M\ge\lambda L^2.\nonumber 
\end{align*}
By the triangle inequality for $\|\cdot\|_2$ and its monotonicity,
cf.~(\ref{ao24}), this yields
by (\ref{ao117}) and (\ref{c1}) in Lemma~\ref{L:4} in conjunction with (\ref{ao16})
\begin{align*}
\lambda\bar M\|I(M(h^*)\in[\frac{\bar M}{4},\frac{\bar M}{2}])\|_2
&\lesssim 
L\sqrt{\bar M}+
\sum_{y=1}^{L-1}\sqrt{\lambda M(\phi_y)}\quad\mbox{provided}\quad\bar
M\ge\lambda L^2.
\end{align*}
In view of (\ref{ao10}) this implies
\begin{align}\label{ao17}
\|I(M(h^*)\in[\frac{\bar M}{4},\frac{\bar M}{2}])\|_2
&\lesssim
\frac{L}{\lambda\sqrt{\bar M}}
+\frac{L^2}{\sqrt{\lambda}\bar M}
\quad\mbox{provided}\quad\bar M\ge\lambda L^2.
\end{align}
%
We choose $\lambda$ such that the two r.~h.~s.~terms in (\ref{ao17}) 
balance, that is, $\lambda=\frac{\bar M}{L^2}$. Then (\ref{ao17}) turns into
\begin{align}\label{ao18}
\|I(M(h^*)\in[\frac{\bar M}{4},\frac{\bar M}{2}])\|_2
&\lesssim(\frac{L^2}{\bar M})^\frac{3}{2}.
\end{align}
It follows from definition (\ref{ao24}) that for an event $A$ we have
\begin{align*}
\mathbb{P}(A)\exp\Big(\frac{1}{\|I(A)\|_2^2}\Big)\leq
\mathbb{E}\exp\big((\frac{I(A)}{\|I(A)\|_2})^2\big)\leq e
\end{align*}
which implies, for $\|I(A)\|_2\ll 1$,
\begin{align*}
-\ln\mathbb{P}(A)\gtrsim\frac{1}{\|I(A)\|_2^2}.
\end{align*}
We then learn from (\ref{ao18}) that
\begin{align*}
-\ln\mathbb{P}\big(M(h^*)\in[\frac{\bar M}{4},\frac{\bar M}{2}]\big)
&\gtrsim(\frac{\bar M}{L^2})^3\quad\mbox{provided}\quad\bar M\gg L^2.
\end{align*}
By the elementary estimate $\ln\sum_{k\ge 0}\exp(-2^k\nu)$ $\lesssim-\nu$
for $\nu\ge 1$
this in turn implies
\begin{align*}
-\ln\mathbb{P}\big(M(h^*)\ge\bar M\big)
&\gtrsim (\frac{\bar M}{L^2})^3\quad\mbox{provided}\quad\bar M\gg L^2.
\end{align*}
%
It remains to apply Lemma~\ref{L:8} to $X = M(h_*)/(CL^2)$, where $C$ is universal and
large enough. 

\qed

{\sc Proof of Corollary~\ref{C:1}}. We start from \eqref{ao89} with $\bar M = \lambda L^2$,
which implies
\begin{align*}
	& \mathbb{P}\big(|h_*(y)| \geq \lambda L\;\mbox{and}\;M(h_*) \leq \frac{\lambda L^2}{2}\big)
	\leq \mathbb{P}\Big(|\min_{M\le
	\lambda L^2}E-\mathbb{E}\min_{M\le\lambda L^2}E| \\ & +
	|\min_{M\le\lambda L^2}E_{y, \lambda} -\mathbb{E}\min_{M\le\lambda L^2}E_{y,
	\lambda}| + |\min_{M\le\lambda L^2}E_{y, -\lambda}
		-\mathbb{E}\min_{M\le\lambda L^2}E_{y, -\lambda}| \\ & 
	+|E(\lambda\phi_y)-\mathbb{E}E(\lambda\phi_y)| +
	|E(-\lambda\phi_y)-\mathbb{E}E(-\lambda\phi_y)| \geq \frac{\lambda^2 L}{2}\Big).
\end{align*}
This inequality, together with the equivalence in law \eqref{ao16},
the concentration properties (\ref{ao117}) and 
\eqref{c1}, and the fact that 
\[
	M(\pm \lambda \phi_y) = \lambda M(\phi_y) 
        \stackrel{(\ref{ao10})}{\le}
	\frac{\lambda L^2}{8},
\]
yields by Chebyshev's inequality (see \eqref{ao124} in Lemma \ref{L:8})
\begin{align*}
\ln\mathbb{P}(|h_*(y)| \geq \lambda L\;\mbox{and}\;M(h_*) \leq \frac{\lambda}{2}L^2)
\lesssim-\frac{(\lambda^2L)^2}{\lambda L^2}=-\lambda^3\qquad\mbox{for
all}\;\lambda\gg1.
\end{align*}
Combined with Proposition~\ref{P:1} 
\[
	\ln\mathbb{P}\big(M(h_*) > \frac{\lambda}{2}L^2\big)
	\overset{\eqref{ao124}}{\lesssim}-\lambda^3 \qquad
	\mbox{for all}\;\lambda\gg1,
\]
we obtain
\begin{align*}
\ln\mathbb{P}(|h_*(y)| \geq \lambda L)
\lesssim-\lambda^3\quad\mbox{provided}\;\lambda\gg 1.
\end{align*}
We conclude by appealing to Lemma~\ref{L:8}.
%
\qed

\medskip

{\sc Proof of Lemma~\ref{L:MP}.} At the core is the inequality 
\begin{align}\label{e40}
D(\min\{h,h'\}) + D(\max\{h,h'\}) \leq D(h) + D(h')\quad\mbox{for a pair of configurations}
\;(h,h'),
\end{align}
which would be obvious as an equality in the continuum case, 
and for which we give the argument below. 
Since the same obviously holds with equality for the local
$W$ replacing $D$, the inequality (\ref{e40}) also holds with $E=D-W$ replacing $D$:
\begin{align}\label{ao126}
E(\min\{h,h'\})+E(\max\{h,h'\}) \leq E(h) + E(h')\quad\mbox{for a pair of configurations}
\;(h,h').
\end{align}
We now turn to (\ref{ao125}) and treat the case $h_0,h_1\ge 0$.
Applying (\ref{ao126}) to $h=h_{*}$ and $h'=h_{h_0,h_1,*}$
we learn that $\min\{h,h'\}$ and $\max\{h,h'\}$ are minimizers
given their boundary condition $(0,0)$ and $(h_0,h_1)$, respectively.
By the uniqueness statement of Lemma \ref{L:7} we thus must have
$\min\{h_*,h_{h_0,h_1,*}\}$ $=h_*$ almost surely, which amounts
to (\ref{ao125}). Since Lemma~\ref{L:7} only refers to $h_*$ and not to
$h_{h_0,h_1,*}$, the null set is independent on $(h_0,h_1)$.

\medskip

We now turn to the argument for the inequality (\ref{e40}), 
which we establish on every interval $\{x-1,x\}$ separately where it assumes the form
%
\begin{align*}
&\frac{1}{2}(\min\{h_1,h_1'\}-\min\{h_0,h_0'\})^2
+\frac{1}{2}(\max\{h_1,h_1'\}-\max\{h_0,h_0'\})^2\\
&\leq \frac{1}{2}(h_1-h_0)^2+\frac{1}{2}(h_1'-h_0')^2
\end{align*}
with $(h_0,h_0',h_1,h_1')$ $=(h(x-1),h'(x-1),h(x),h'(x))$.
Expanding the squares and using the obvious $(\min\{h_i,h_i'\})^2+(\max\{h_i,h_i'\})^2$
$=h_i^2+{h_i'}^2$, we see that the above is equivalent to 
\[
h_0h_1+h_0'h_1'\leq \min\{h_0,h_0'\}\min\{h_1,h_1'\} + \max\{h_0,h_0'\}\max\{h_1,h_1'\}.
\]
We may w.~l.~o.~g.~assume $h_0\le h_0'$
so that the above turns into
$h_0h_1+h_0'h_1'$ $\leq h_0\min\{h_1,h_1'\}$ $+h_0'\max\{h_1,h_1'\}$, 
which can be reformulated as 
$h_0\max\{h_1-h_1',0\}$ $\leq h_0'\max\{h_1-h_1',0\}$ and thus holds true.
\qed 

\medskip

{\sc Proof of Corollary \ref{C:2}}.
The passage from Corollary \ref{C:1} (with $x=\frac{L}{2}$)
to Corollary \ref{C:2} relies on two ingredients
involving the affine interpolation of the boundary values
\begin{align*}
\bar h_{h_0,h_1}(x):=h_0(1-\frac{x}{L})+h_1\frac{x}{L},
\end{align*}
namely shear invariance, see also \cite[Corollary 2.2]{Peled}, and the comparison principle.

\medskip

We start with the shear-invariance,
which we use in a more restricted form than (\ref{ao16}) 
in the proof of Proposition \ref{P:1}:
\begin{align*}
h\mapsto W(h-\bar h_{h_0,h_1})-W(h_{h_0,h_1})
\quad\mbox{has same law as}\quad h\mapsto W(h)
\end{align*}
in conjunction with the obvious orthogonality
\begin{align*}
D(h-\bar h_{h_0,h_1})-D(\bar h_{h_0,h_1})=D(h),
\end{align*}
which implies
\begin{align}\label{ao134}
	h_{h_0,h_1,*}-\bar h_{h_0,h_1}\quad\mbox{has same law as}\quad h_*.
\end{align}
Since non-uniqueness may hold, \eqref{ao134} must be interpreted as equality in
law for sets of minimizers. From \eqref{ao134} we have, in particular,
\begin{align}\label{ao23}
	h_{h_0, h_1, *}(\frac{L}{2})-\frac{1}{2}(h_1+h_0)\quad\mbox{has same law as}\quad h_*(\frac{L}{2}).
\end{align}
With the same reasoning, \eqref{ao23} can be upgraded as constancy for the joint
law of minimizers within the same bin of boundary conditions, namely
\begin{align*}
\mbox{the joint law
of}\;\big\{h_{h_0,h_1,*}-\frac{1}{2}(h_0+h_1)\big\}_{h_0\in[\bar
h_0-\frac{L}{2},\bar h_0+\frac{L}{2}]}^{h_1\in[\bar
h_1-\frac{L}{2},\bar h_1+\frac{L}{2}]}
\;\mbox{does not depend on}\;(\bar
h_0,\bar h_1).
\end{align*}
This implies (\ref{ao31}). It also implies that we may extend the comparison principle
(\ref{ao125}) to
\begin{align}\label{ao127}
h_{h_0,h_1,*}\le h_{h_0',h_1', *}\quad\mbox{for}\;h_0\le h_0'\;\mbox{and}\;h_1\le h_1'
\quad\mbox{almost surely},
\end{align}
where the null set can be chosen to be independent of either $(h_0,h_1)$ or
$(h_0',h_1')$.

\medskip

Turning now to (\ref{ao28}) 
we may restrict ourselves to $\bar h_0 = \bar h_1 = 0$ by (\ref{ao31}).
By the comparison principle in form of (\ref{ao127}) we have
\[
	h_{-\frac{L}{2},-\frac{L}{2},*}(\frac{L}{2}) 
\leq h_{h_0, h_1, *}(\frac{L}{2}) \leq
	h_{\frac{L}{2},\frac{L}{2},*}(\frac{L}{2})
\quad\mbox{for all}\;\;(h_0,h_1)\in[-\frac{L}{2},\frac{L}{2}]^2
\;\;\mbox{almost surely}.
\]
With an eye on (\ref{ao23}) we rewrite this as
\begin{align*}
\big(h_{-\frac{L}{2},-\frac{L}{2},*}(\frac{L}{2})+\frac{L}{2}\big)-L
\leq h_{h_0, h_1, *}(\frac{L}{2})-\frac{1}{2}(h_0+h_1)
\leq \big(h_{\frac{L}{2},\frac{L}{2},*}(\frac{L}{2})-\frac{L}{2}\big)+L
\end{align*}
for all $(h_0,h_1)\in[-\frac{L}{2},\frac{L}{2}]^2$ almost surely,
so that by definition (\ref{ao30})
\begin{align*}
LX(0,0)\le\max\big\{-\big(h_{-\frac{L}{2},-\frac{L}{2},*}(\frac{L}{2})+\frac{L}{2}\big),
\big(h_{\frac{L}{2},\frac{L}{2},*}(\frac{L}{2})-\frac{L}{2}\big)\big\}+L\quad\mbox{almost surely},
\end{align*}
to which we apply (the monotone norm) $\|\cdot\|_3$, and then appeal to the identity in law
(\ref{ao23}):
\begin{align*}
L\|X(0,0)\|_3\le 2\|h_*(\frac{L}{2})\|_3+L,
\end{align*}
so that (\ref{ao28}) now follows from (\ref{ao21}).\qed
\ignore{
By definition (\ref{ao30}) this implies
\[
	X(0,0) \leq \frac{1}{L}(|h_{-\frac{L}{2},-\frac{L}{2},*}(\frac{L}{2})+L| +
	|h_{\frac{L}{2},\frac{L}{2},*}(\frac{L}{2})-L|) + 1.
\]
Applying the norm $\|\cdot\|_3$, using Corollary~\ref{C:1} and \eqref{ao23}, one
concludes.
\qed
}

\medskip 

{\sc Proof of Lemma \ref{L:1}}. The passage from Corollary \ref{C:2} to Lemma \ref{L:1}
relies on (\ref{ao46}), which we rewrite as 
\begin{align*}
\big|\frac{h_*((2\hat x-1)l)-\frac{1}{2}\big(h_*(2(\hat x-1)l)
+h_*(2\hat x l)\big)}{l}\big|\le 2X_{\hat x},
\end{align*}
where the random variable $X_{\hat x}$ abbreviates 
$X_{[2(\hat x-1)l,2\hat xl]}(\bar h(2(\hat x-1)l),\bar h(2\hat xl))$. 
Hence taking the cubic power and summing over $\hat x$ we obtain from 
(\ref{ao27}) for $p=3$
\begin{align}\label{ao111}
\frac{D_3(h_{*,l})}{L}\le
8\frac{2l}{L}\sum_{\hat x=1}^{\frac{L}{2l}-1} X_{\hat x}^3.
\end{align}
Since $X_{\hat x}$ depends on $W$ only through 
$\{W(x,\cdot)\}_{x\in(2l(\hat x-1),2l\hat x)\cap\mathbb{Z}}$ and by our assumption of independence, 
we have
\begin{align*}
\{X_{\hat x}\}_{\hat x=1,\cdots,\frac{L}{2l}-1}
\quad\mbox{are independent and identically distributed},
\end{align*}
where the latter follows from (\ref{ao31}).
In order to capitalize on this independence, we monitor the characteristic
function:
\begin{align}\label{ao97}
\mathbb{E}\exp(\lambda\frac{D_3(h_{*,l})}{L})
\stackrel{(\ref{ao111})}{\le}\big(\mathbb{E}
\exp(8\lambda\frac{2l}{L}X_1^3)\big)^{\frac{L}{2l}-1}\qquad\text{for}\;\lambda\geq0.
\end{align}
We choose
\begin{align}\label{ao98}
\lambda=\frac{1}{8}\frac{L}{2l}\frac{1}{\|X_1\|_3^3}\quad\mbox{so that}\quad
8\lambda\frac{2l}{L}=\frac{1}{\|X_1\|_3^3},
\end{align}
to the effect that by definition (\ref{ao24}) of $\|\cdot\|_s$, (\ref{ao97}) turns into
\begin{align*}
\mathbb{E}\exp(\lambda\frac{D_3(h_{*,l})}{L})
\le e^{\frac{L}{2l}-1}\le e^{\frac{L}{2l}}.
\end{align*}
By 
\begin{align*}
\mathbb{P}(\frac{D_3(h_{*,l})}{L}\ge\nu)
\le \exp(\frac{L}{2l}-\lambda\nu)\qquad\text{for all}\;\nu\geq0,
\end{align*}
which by the choice (\ref{ao98}) of $\lambda$ turns into
\begin{align*}
\mathbb{P}(\frac{D_3(h_{*,l})}{L}\ge\nu)
\le \exp\big(\frac{L}{2l}(1-\frac{1}{8}\frac{1}{\|X_1\|_3^3}\nu)\big).
\end{align*}
This yields (\ref{ao34}) using (\ref{ao28}).\qed

\medskip

{\sc Proof of Proposition~\ref{P:2}}. From \eqref{ao46} and \eqref{ao29}, we learn that if
\[
	|h_*(2(\hat x -1)l)|, |h_*(2\hat xl)| \leq \hat H L
\]
then
\[
	\Big|\frac{h_*((2\hat x-1)l)-\frac{1}{2}\big(h_*(2(\hat x-1)l)+h_*(2\hat xl)\big)}{2l}\Big|
	\le \sup_{h_0, h_1 \in [-\hat HL, \hat HL] \cap 2l\mathbb{Z}}X_{[2(\hat x-1)l,2\hat
	xl]}(h_0,h_1).
\]
From \eqref{ao28} in Corollary~\ref{C:2}, $\|X_{[2(\hat x-1)l,2\hat xl]}(h_0,h_1)\|_3 \lesssim
1$ and thus by Chebyshev's inequality (see \eqref{ao124} in Lemma \ref{L:8})
\begin{align}\label{ao118}
\ln \mathbb{P}(X_{[2(\hat x-1)l,2\hat xl]}(h_0,h_1) \geq \nu
\ln^{\frac{1}{3}}\frac{\hat HL}{2l}) \lesssim -\nu^3
\ln\frac{\hat HL}{2l}\qquad\text{for}\;\nu,\hat H\gg1,
\end{align}
so that, by the union bound, 
\begin{align*}
	& \ln \mathbb{P}\Big(\sup_{h_0, h_1 \in [-\hat HL, \hat HL] \cap 2l\mathbb{Z}}X_{[2(\hat x-1)l,2\hat
xl]}(h_0,h_1) \geq \nu \ln^{\frac{1}{3}}\frac{\hat HL}{2l}\Big) \\ & \leq
\ln \Big(\sum_{h_0, h_1 \in
	[-\hat HL, \hat HL] \cap 2l\mathbb{Z}} \mathbb{P}(X_{[2(\hat x-1)l,2\hat
xl]}(h_0,h_1) \geq \nu\ln^{\frac{1}{3}}\frac{\hat HL}{2l})\Big) \\ &
\overset{\eqref{ao118}}{\lesssim}
\ln\frac{\hat HL}{2l}-\nu^3\ln\frac{\hat
HL}{2l}\lesssim
-\nu^3\qquad\text{for}\;\nu,\hat H\gg1.
\end{align*}
With the choice $\hat H = \nu$ and by Corollary~\ref{C:1}, we obtain
\begin{align}\label{ao122}
&\ln\mathbb{P}\Big(\Big|\frac{h_*((2\hat x-1)l)-\frac{1}{2}\big(h_*(2(\hat
x-1)l)+h_*(2\hat xl)\big)}{2l}\Big| \geq \nu\ln^{\frac{1}{3}}\frac{\nu
L}{2l}\Big)
\nonumber\\ & \leq \ln\mathbb{P}\Big(|h_*(2(\hat x -1)l)|>\nu L\; \text{or}\; |h_*(2\hat
xl)| > \nu L\nonumber\\ &\qquad\text{or} \;\sup_{h_0, h_1 \in [-\nu L, \nu L] \cap
2l\mathbb{Z}}X_{[2(\hat x-1)l,2\hat
xl]}(h_0,h_1) \geq \nu \ln^{\frac{1}{3}}\frac{\nu L}{2l}\Big) \lesssim -\nu^3
\end{align}

\medskip

Since for $\epsilon > 0$ and $\nu \gg_\epsilon 1$
\begin{align*}
\ln^{\frac{1}{3}}\frac{\nu L}{2l}=\Big(\ln\nu+\ln
\frac{L}{2l}\Big)^{\frac{1}{3}}\lesssim\ln^{\frac{1}{3}}\nu+
\ln^{\frac{1}{3}}\frac{L}{2l}\leq(\ln^{\frac{1}{3}}\frac{L}{2l})\nu^\epsilon,
\end{align*}
then \eqref{ao122} implies the weaker
\begin{align*}
\ln\mathbb{P}\Big(\Big|\frac{h_*((2\hat x-1)l)-\frac{1}{2}\big(h_*(2(\hat
x-1)l)+h_*(2\hat xl)\big)}{2l}\Big| \geq
(\ln^{\frac{1}{3}}\frac{L}{2l})\nu^{1+\epsilon}\Big)\lesssim_\epsilon-\nu^{3}.
\end{align*}
An application of \eqref{ao39} in Lemma~\ref{L:8} thus yields for any $1\leq s<3$,
\begin{equation}\label{e22}
\Big\|\frac{h_*((2\hat x-1)l)-\frac{1}{2}\big(h_*(2(\hat
x-1)l)+h_*(2\hat xl)\big)}{2l}\Big\|_s \lesssim_s
\ln^{\frac{1}{3}}\frac{L}{2l}.
\end{equation}

\medskip

We proceed with an iteration argument on dyadic scales. Let $0 \leq x_0 = x < x_1 = y \leq L$. We have that either $x_0 \leq
L/2$ or $x_1 \geq L/2$. Without loss of generality assume that $x_0 \leq
L/2$. Define $x_n = x_0 + 2^{n-1}(x_1 - x_0)$. The sequence $(x_n)_{n \geq
0}$ is increasing and we have that
\begin{align}\label{ao119}
	x_n - x_0 = 2^{n-1}(x_1 - x_0).
\end{align}
Let $N$ be the largest positive integer such that $x_N \leq L$. Since $0 \leq x_0 \leq L/2$, then
\begin{align}\label{ao120}
	\log_2\frac{L}{|x_1 - x_0|} - 1 \leq N \leq \log_2\frac{L}{|x_1 - x_0|}
	+ 1.
\end{align}
We can estimate the modulus of continuity in the following way
\begin{equation}\label{e11}
	|\frac{h_*(x_1) - h_*(x_0)}{x_1 - x_0}| \leq \sum_{n =
	1}^{N-1} |\frac{h_*(x_{n+1}) - h_*(x_0)}{x_{n+1} - x_0} -
	\frac{h_*(x_n) - h_*(x_0)}{x_n - x_0}| + |\frac{h_*(x_N) -
	h_*(x_0)}{x_N - x_0}|.
\end{equation}
Due to the dyadic structure of the $ x_n $'s, each of the terms in the sum is of the form
\begin{gather*}
	\frac{h_*(x_{n+1}) - h_*(x_0)}{x_{n+1} - x_0} - \frac{h_*(x_n) -
		h_*(x_0)}{x_n - x_0} = -\frac{2}{x_{n+1} -
	x_0}\Big(h_*\big(\frac{x_{n+1} +
x_0}{2}\big) - \frac{h_*(x_{n+1}) + h_*(x_0)}{2}\Big)
\end{gather*}
and therefore can be estimated using \eqref{e22} with $l=x_n-x_0$ to get
\begin{equation}\label{e12}
\Big\|\frac{h_*(x_{n+1}) - h_*(x_0)}{x_{n+1} - x_0} - \frac{h_*(x_n) -
h_*(x_0)}{x_n - x_0}\Big\|_{s} \lesssim_s
\ln^{\frac{1}{3}}\frac{L}{2^{n-1}|x_0-x_1|}.
\end{equation}
Using \eqref{ao119} \& \eqref{ao120} in the form $x_N - x_0 \geq
\frac{L}{4}$, one obtains
\begin{equation}\label{e13}
	|\frac{h_*(x_N) - h_*(x_0)}{x_N - x_0}| \lesssim
	\frac{1}{L}(|h_*(x_N)| + |h_*(x_0)|),
\end{equation}
which by Corollary \ref{C:1} translates to
\begin{align}\label{ao123}
	\Big\|\frac{h_*(x_N) - h_*(x_0)}{x_N - x_0}\Big\|_s \lesssim_s 1.
\end{align}
Finally \eqref{e12} \& \eqref{ao123} yield
\[
\Big\|\frac{h_*(x_1) - h_*(x_0)}{x_1 - x_0}\Big\|_{s}
\lesssim_s
1+\sum_{n=1}^{N-1}\ln^{\frac{1}{3}}\frac{L}{2^{n-1}|x_0-x_1|}
\lesssim1+\ln^{\frac{4}{3}}\frac{L}{|x_0-x_1|}.
\]
\qed

\medskip

{\sc Proof of Lemma \ref{L:2}}. As in (\ref{ao67}) we have
\begin{align}\label{ao72}
\bar h=\sum_{\rho}\bar h_\rho,
\end{align}
where the sum over $\rho$ restricts to (\ref{ao71}) since $\bar h_{\ge 2l}=\bar
h$, cf.~\eqref{ao32}.
According to (\ref{ao112}), each component $\bar h_\rho$ is parameterized by its 
values in $(2\hat x-1)\rho$
for $\hat x\in\{1,\cdots,\frac{L}{2\rho}-1\}$. Taking also
\eqref{ao112} \& \eqref{ao27} into account 
we see
\begin{align}\label{ao35}
\frac{D(\bar h_\rho)}{L}=\frac{\rho}{L}
\sum_{\hat x=1}^{\frac{L}{2\rho}-1}
\big(\frac{\bar h_\rho((2\hat x-1)\rho)}{\rho}\big)^2.
\end{align}
We learn from the representation (\ref{ao112}) that (\ref{ao32}) implies
%
\begin{align*}
\bar h_{\rho}((2\hat x-1)\rho) 
\in l\mathbb{Z}\quad\mbox{for}\;\hat x=1,\cdots,\frac{L}{2\rho}-1.
\end{align*}
Together with (\ref{ao72}) \& (\ref{ao35}) we obtain
\begin{align}\label{ao73}
\ln\#\big\{\,\bar h\;\mbox{of form}\;(\ref{ao32})\,|\,
\max_\rho(\frac{l}{\rho})^2\frac{D(\bar h_\rho)}{L}<\hat D\,\big\}
\le\sum_\rho \ln Z\big(\frac{L}{2\rho},2(\frac{\rho}{l})^4\hat D\big),
\end{align}
where for $N\in\mathbb{N}$ and $D\in[0,\infty)$, $Z(N,D)$ is defined
via\footnote{The translation from $\bar h_\rho$ to the variables
$y_{\hat x}$ is done via $ly_{\hat x}:=\bar h_\rho((2\hat x-1)\rho)$.}
\begin{align}\label{ao75}
Z(N,D):=\#\big\{\,(y_1,\cdots,y_N)\in\mathbb{Z}^N\,|
\,\frac{1}{N}\sum_{\hat x=1}^Ny_{\hat x}^2<D\,\big\}.
\end{align}

\medskip

We now claim that (\ref{ao69}) follows from 
\begin{align}\label{ao74}
\ln Z(N,D)\lesssim N\ln D\quad\mbox{for}\;D\gg 1.
\end{align}
Indeed, inserting (\ref{ao74}) into (\ref{ao73}) we obtain
\begin{align*}
\lefteqn{\ln\#\big\{\,\bar h\;\mbox{of form}\;(\ref{ao32})\,|\,
\max_\rho(\frac{l}{\rho})^2\frac{D(\bar h_\rho)}{L}<\hat D\,\big\}}\nonumber\\
&\lesssim\sum_\rho \frac{L}{2\rho}\ln\big((\frac{\rho}{l})^4\hat
D\big)
=\frac{L}{2l}\sum_\rho \frac{l}{\rho}\big(4\ln\frac{\rho}{l}+\ln\hat D)
	\stackrel{\hat D\gg 1}{\lesssim}\frac{L}{l}\ln\hat D.
\end{align*}
%

\medskip

We finally turn to the proof of (\ref{ao74}) in the more precise form of
\begin{align}\label{ao70}
Z(N,D)
\le(C_0(D+1))^\frac{N}{2}\quad\mbox{for all}\;D\ge 0
\end{align}
which we reproduce for the convenience of the reader and establish by induction in $N$. 
Inequality (\ref{ao70})
is obviously satisfied for $N=1$ with $C_0=1$, cf.~(\ref{ao75}). The elementary
inclusion
\begin{align*}
\lefteqn{\{\,y\in\mathbb{Z}^{N+1}\,|
\,\sum_{\hat x=1}^{N+1}y_{\hat x}^2<(N+1) D\,\}}\nonumber\\
&\subset \bigcup_{y_{N+1}\in\mathbb{Z}\cap(-\sqrt{(N+1)D},\sqrt{(N+1)D})}
\{\,(y',y_{N+1})\,|\,y'\in\mathbb{Z}^{N}\;\mbox{and}\;
\sum_{\hat x=1}^{N}{y'}_{\hat x}^2< (N+1)D-y_{N+1}^2\,\}
\end{align*}
translates into
\begin{align*}
Z(N+1,D)\le\sum_{r\in\mathbb{Z}\cap(-\sqrt{(N+1)D},\sqrt{(N+1)D})}
Z(N,\frac{(N+1)D-r^2}{N}).
\end{align*}
Hence in order to propagate (\ref{ao70}), we need to show that there exists a $1\le C_0<\infty$
such that
\begin{align*}
(C_0(D+1))^{\frac{N+1}{2}}\ge\sum_{r\in\mathbb{Z}\cap(-\sqrt{(N+1)D},\sqrt{(N+1)D})}
\big(C_0\frac{(N+1)D-r^2+N}{N}\big)^\frac{N}{2},
\end{align*}
which we rewrite as
\begin{align*}
\sum_{r\in\mathbb{Z}\cap(-\sqrt{(N+1)D},\sqrt{(N+1)D})}
\frac{1}{\sqrt{D+1}}\big(1+\frac{D-r^2}{N(D+1)}\big)^\frac{N}{2} \le\sqrt{C_0}.
\end{align*}
Appealing to the elementary inequality on the summand
\begin{align*}
\big(1+\frac{D-r^2}{N(D+1)}\big)^\frac{N}{2}
\le\big(\exp(\frac{D-r^2}{N(D+1)})\big)^\frac{N}{2}
=\exp(\frac{D-r^2}{2(D+1)})\le\sqrt{e}\exp(-\frac{r^2}{2(D+1)}),
\end{align*}
we see that it is enough to show
\begin{align*}
\sum_{r\in\mathbb{Z}}
\frac{1}{\sqrt{D+1}}\exp(-\frac{r^2}{2(D+1)})\le\sqrt{\frac{C_0}{e}}.
\end{align*}
Appealing to the monotonicity of the integrand to control the sum by an integral, it suffices to show
\begin{align*}
\frac{1}{\sqrt{D+1}}+2\int_0^{\infty}\frac{dr}{\sqrt{D+1}}
\exp(-\frac{r^2}{2(D+1)})\le\sqrt{\frac{C_0}{e}}.
\end{align*}
Under the change of variables $r=\sqrt{D+1}\hat r$, this turns into
\begin{align*}
\frac{1}{\sqrt{D+1}}+\int_{-\infty}^{\infty}d\hat r
\exp(-\frac{\hat r^2}{2})\le\sqrt{\frac{C_0}{e}},
\end{align*}
which clearly is satisfied for $C_0$ large enough.
\qed

\medskip

{\sc Proof of Corollary \ref{C:3}}. Note that (\ref{ao38}) implies by the union
bound
\begin{align*}
\lefteqn{\mathbb{P}\big(\max_\rho(\frac{l}{\rho})^2\frac{D(h_{*,\rho})}{L}\le\hat D
\quad\mbox{and}\quad\frac{D_3(h_{*,l})}{L\ln\hat D}\ge\nu\big)}\nonumber\\
&\le
\sum_{\bar h\;\mbox{of the form (\ref{ao32}) with}
\;\max_\rho(\frac{l}{\rho})^2\frac{D(\bar h_\rho)}{L}\le 4\hat D}
\mathbb{P}\big(\mbox{(\ref{ao29}) holds and}
\quad\frac{D_3(h_{*,l})}{L}\ge\nu\ln\hat D\big).
\end{align*}
According to Lemmas \ref{L:1} and \ref{L:2}, where we make the implicit constants
explicit as $c_1>0$ and $C_2<\infty$, this implies for $\nu\ge 1$
(and thus $\nu\ln\hat D\gg 1$)
\begin{align*}
\ln\mathbb{P}\big(\max(\frac{l}{\rho})^2\frac{D(h_{*,\rho})}{L}\le\hat D
\quad\mbox{and}\quad\frac{D_3(h_{*,l})}{L\ln\hat D}\ge\nu\big)
\le C_2\frac{L}{l}\ln\hat D
-c_1\frac{L}{l}\nu\ln\hat D.
\end{align*}
Obviously, for $\nu\gg 1$, the second term dominates the first; 
in addition, since $\hat D\gg 1$, the sum is $\le-\nu$.
Hence we may apply Lemma~\ref{L:8} to $s=1$ and $X$ 
$=I(\max_\rho(\frac{l}{\rho})^2\frac{D(h_{*,\rho})}{L}$ 
$\le\hat D)\frac{D_3(h_{*,l})}{L\ln\hat D}$. The output
$\|X\|_1\lesssim 1$ yields the desired (\ref{ao37}). \qed

\medskip

{\sc Proof of Lemma \ref{L:3}}. We fix a 
\begin{align}\label{ao76}
\hat D\sim e+\max_\rho\big\|\frac{D(h_{*,\rho})}{L}\big\|_1
\quad\mbox{and such that Corollary \ref{C:3} holds},
\end{align}
dyadically bin the random variable $\max_\rho(\frac{l}{\rho})^2\frac{D(h_{*,\rho})}{L}$
in terms of units of $\hat D$, and use the triangle inequality and a union bound to obtain
\begin{align*}
\big\|\frac{D_p(h_{*,l})}{L}\big\|_1
&\le\big\|I(\max_\rho(\frac{l}{\rho})^2\frac{D(h_{*,\rho})}{L}\le\hat D)
\frac{D_p(h_{*,l})}{L}\big\|_1\nonumber\\
&+\sum_{k=1}^\infty
\big\|I(2^{k-1}\hat D<\max_\rho(\frac{l}{\rho})^2
\frac{D(h_{*,\rho})}{L}\le 2^k\hat D)\frac{D_p(h_{*,l})}{L}\big\|_1\\
&\le\big\|I(\max_\rho(\frac{l}{\rho})^2\frac{D(h_{*,\rho})}{L}\le\hat D)
\frac{D_p(h_{*,l})}{L}\big\|_1\nonumber\\
&+\sum_{k=1}^\infty\sum_{\rho'}
\big\|I\big(\frac{D(h_{*,\rho'})}{L}> 2^{k-1}(\frac{\rho'}{l})^2\hat D\big)\,
I\big(\max_\rho(\frac{l}{\rho})^2\frac{D(h_{*,\rho})}{L}\le 2^k\hat D\big)\,
\frac{D_p(h_{*,l})}{L}\big\|_1,
\end{align*}
where also $\rho'$ ranges over (\ref{ao71}). In order to capitalize on
the unlikeliness of $\frac{D(h_{*,\rho'})}{L}$ $>2^{k-1}(\frac{\rho'}{l})^2\hat D$
in the last summand, we need to invest stochastic integrability:
Equipped with the following auxiliary statements on our norms
(\ref{ao24})
\begin{align}
\|X\|_1&\lesssim\|X\|_{\frac{3}{p}}\qquad\mbox{for}~1\le p\le3,\label{ao41}
\end{align}
which is a direct consequence of Lemma~\ref{L:8} and that we use
on the first r.~h.~s.~term, and (see Lemma~\ref{L:15}.(iii) in the appendix)
\begin{align}
\|I(Y\ge\bar Y)X\|_1&\lesssim(\frac{\|Y\|_1}{\bar Y})^{1-\frac{p}{3}}
\|X\|_{\frac{3}{p}}\label{ao42}
\end{align}
for two random variables $X,Y\ge 0$ and a constant $\bar Y\ge 1$ and use for 
$X=I(\max_\rho(\frac{l}{\rho})^2\frac{D(h_{*,\rho})}{L}\le 2^k\hat D)\frac{D_p(h_{*,l})}{L}$,
$Y=\frac{D(h_{*,\rho'})}{L}$, and $\bar Y=2^{k-1}(\frac{\rho'}{l})^2\hat D$, 
this yields by our choice (\ref{ao76}) of $\hat D$ that ensures
$\|\frac{D(h_{*,\rho'})}{L}\|_s$ $\lesssim\hat D$
\begin{align*}
\big\|\frac{D_p(h_{*,l})}{L}\big\|_1
&\lesssim
\big\|I(\max_\rho(\frac{l}{\rho})^2\frac{D(h_{*,\rho})}{L}\le\hat D)
\frac{D_p(h_{*,l})}{L}\big\|_{\frac{3}{p}}\nonumber\\
&+\sum_{k=1}^\infty\sum_{\rho'}(\frac{1}{2^{k-1}}(\frac{l}{\rho'})^2)^{1-\frac{p}{3}}
\big\|I(\max_\rho(\frac{l}{\rho})^2\frac{D(h_{*,\rho})}{L}\le 2^k\hat D)
\frac{D_p(h_{*,l})}{L}\big\|_{\frac{3}{p}}.
\end{align*}
The sum in $\rho'$ is geometric and thus dominated by $\rho'=l$ so that
\begin{align*}
\big\|\frac{D(h_{*,l})}{L}\big\|_1
&\lesssim_p
\sum_{k=0}^\infty(\frac{1}{2^k})^{1-\frac{p}{3}}
\big\|I(\max_\rho(\frac{l}{\rho})^2\frac{D(h_{*,\rho})}{L}\le 2^k\hat D)
\frac{D_p(h_{*,l})}{L}\big\|_{\frac{3}{p}}\nonumber\\
&\stackrel{(\ref{ao99})}{\le}
\sum_{k=0}^\infty(\frac{1}{2^k})^{1-\frac{p}{3}}
\big\|I(\max_\rho(\frac{l}{\rho})^2\frac{D(h_{*,\rho})}{L}\le 2^k\hat D)
\frac{D_3(h_{*,l})}{L}\big\|_{1}^\frac{p}{3}.
\end{align*}
We now appeal to Corollary \ref{C:3} and our choice (\ref{ao76}) of $\hat D$ to obtain
\begin{align*}
\big\|\frac{D(h_{*,l})}{L}\big\|_p
&\lesssim_p\sum_{k=0}^\infty(\frac{1}{2^k})^{1-\frac{p}{3}}
\ln^\frac{p}{3}2^k(e+\max_\rho\|\frac{D(h_{*,\rho})}{L}\big\|_1)\nonumber\\
&=\sum_{k=0}^\infty(\frac{1}{2^k})^{1-\frac{p}{3}}
\big(k\ln 2+\ln(e+\max_\rho\big\|\frac{D(h_{*,\rho})}{L}\big\|_1)\big)^\frac{p}{3}.
\end{align*}
Again, in the sum in $k$ the geometric decay dominates the algebraic growth so
that the dominant contribution comes form $k=0$, whence (\ref{ao40}) holds.
\qed

\medskip

{\sc Proof of Theorem \ref{T}}. Since for $l=\frac{L}{2}$ the set (\ref{ao71}) 
of $\rho$'s is empty, (\ref{ao40}) turns into $\|\frac{D(h_{*,\frac{L}{2}})}{L}\|_1$
$\lesssim 1$. Hence we learn from iterating (\ref{ao40}) for $p=2$ that $\|\frac{D(h_{*,l})}{L}\|_1$
and thus $\max_{l\le\rho<L}\|\frac{D(h_{*,\rho})}{L}\|_1$ is (qualitatively) finite. 
Post-processing (\ref{ao40}) to
\begin{align}\label{ao100}
\max_{l\le\rho<L}\|\frac{D_p(h_{*,\rho})}{L}\|_1
\le C_0\ln^\frac{p}{3}(e+\max_{l\le\rho<L}\|\frac{D(h_{*,\rho})}{L}\|_1),
\end{align}
and using (\ref{ao100}) for $p=2$ we thus learn from the qualitative finiteness 
that it is quantitatively finite: $\max_{l\le\rho<L}\|\frac{D(h_{*,\rho})}{L}\|_1$ $\lesssim 1$.
Using this information once more in (\ref{ao100}) we obtain (\ref{ao20}).

\medskip

We now turn to the passage from (\ref{ao20}) to (\ref{ao80}).
It follows from the deterministic estimate
\begin{align}\label{ao105}
\sum_{x=1}^L|h_{\ge l}(x)-h_{\ge l}(x-1)|^p
&\lesssim
\sum_{x=1}^L\big(\sum_\rho(h_\rho(x)-h_\rho(x-1))^2\big)^\frac{p}{2}
\quad\mbox{for}~p\in(1,\infty),
\end{align}
which we shall establish below (see \eqref{ao107}), combined with H\"older's inequality in $\rho$ in form of
\begin{align}\label{ao104}
\big(\sum_\rho(h_\rho(x)-h_\rho(x-1))^2\big)^\frac{p}{2}&
\le(\log_2\frac{L}{2l})^{\frac{p}{2}-1}\sum_\rho|h_\rho(x)-h_\rho(x-1)|^p,
\end{align}
where $\rho$ runs over $l,2l,\cdots,\frac{L}{2}$. Indeed, inserting (\ref{ao104})
into (\ref{ao105}) and dividing by $L$ we obtain by definition (\ref{ao106})
\begin{align*}
\frac{D_p(h_{\ge l})}{L}\lesssim(\ln\frac{L}{2l})^{\frac{p}{2}-1}\sum_\rho\frac{D_p(h_\rho)}{L}.
\end{align*}
Using this for $h=h_*$, applying $\|\cdot\|_1$, 
and appealing to the triangle inequality we obtain
\begin{align*}
\|\frac{D_p(h_{*,\ge l})}{L}\|_1\lesssim (\ln\frac{L}{2l})^{\frac{p}{2}-1}
\sum_\rho\|\frac{D_p(h_{*,\rho})}{L}\|_1.
\end{align*}
Since the r.~h.~s.~sum has $\log_2\frac{L}{2l}$ summands, (\ref{ao20}) yields (\ref{ao80}).

\medskip

We finally turn to the proof of (\ref{ao105}), which we (exactly) rewrite in the continuum form of
\begin{align}\label{ao107}
\int_0^Ldx|\frac{dh_{\ge l}}{dx}|^p
\lesssim_p\int_0^Ldx\big(\sum_\rho(\frac{dh_\rho}{dx})^2\big)^\frac{p}{2}\quad\mbox{for}\;
p\in(1,\infty);
\end{align}
note that the definitions (\ref{ao66}) \& (\ref{ao67}) make also sense in the continuum case. 
Recalling \eqref{ao116}, $\frac{dh_\rho}{dx}$ can be interpreted
in terms of Haar wavelets. In particular,
$(\sum_\rho(\frac{dh_\rho}{dx})^2)^\frac{1}{2}$ acts what is known as the square function
in the case of Littlewood-Paley decomposition, as that (\ref{ao107}) can be interpreted
as an instance of the Littlewood-Paley theorem. 

\medskip

An efficient proof in our case of wavelets runs as follows:
According to \cite{Gaposhkin}, 
the Haar wavelets form an unconditional basis in $L^q(0,L)$ for all $q\in(1,\infty)$, 
which means
\begin{align*}
\int_0^Ldx|\sum_\rho\sigma_\rho\frac{dh_\rho}{dx}|^{q}
\lesssim_q \int_0^Ldx|\frac{dh_{\ge l}}{dx}|^{q}\quad\mbox{for all}
\;\{\sigma_\rho\}_\rho\subset\{-1,+1\}.
\end{align*}
Choosing the $\sigma_\rho$'s to be i.~i.~d.~and taking the expectation, 
this implies by Khintchine's inequality the 
$L^q(0,L)$-boundedness of the square function:
\begin{align*}
\int_0^Ldx\big(\sum_\rho(\frac{dh_\rho}{dx})^2\big)^\frac{q}{2}
\lesssim \int_0^Ldx|\frac{dh_{\ge l}}{dx}|^q.
\end{align*}
This in turn implies (\ref{ao107}) by duality.
\qed

\medskip

{\sc Proof of Lemma~\ref{L:9}.}
We start by specifying the two-scale construction of a competitor, which we will call $\bar h_* $, for
the problem at scale $L$: Given at intermediate scale $1\le
l\in2^{-\mathbb{N}}L$, let us define the set of Brownian motions
\begin{align}\label{f1}
\hat W_{\hat x}(\hat h):=\frac{1}{l}\sum_{x=l(\hat x-1)+1}^{l\hat x}W_x(l\hat h)
\qquad\mbox{for $\hat x\in\{0,\ldots,\frac{L}{l}\}$}
\end{align}
and denote by $\hat h_*$ be the minimizer of the energy
\begin{align*}
\hat E(\hat h)=\frac{1}{2}\sum_{\hat x=1}^{\frac{L}{l}}(\hat h(\hat x)-\hat h(\hat x-1))^2
-\sum_{\hat x=1}^{\frac{L}{l}-1}\hat W_{\hat x}(\hat h(x))
\end{align*}
on the rescaled lattice $\hat
x\in\{0,\ldots,\frac{L}{l}\}$ with zero boundary conditions. Given $\hat h_*$ we introduce the bin
$\bar{\hat h}_*$ as the unique function on the lattice $\{0,\ldots,\frac{L}{l}\}$ such that
\begin{align}\label{f10}
\bar{\hat h}_*(\hat x)\in\mathbb{Z}\quad\mbox{and}\quad \hat h_*(\hat x)\in(\bar{\hat
h}_*(\hat x)-\frac{1}{2},\bar{\hat h}_*(\hat x)+\frac{1}{2}].
\end{align}
Equipped with $\bar{\hat h}_*$, we also define its rescaled version $\bar h$ on the original lattice
$\{0,\ldots,L\}$ with the property that
\begin{align}\label{f9}
\bar h(l\hat x)=l\bar{\hat h}_*(\hat x)\quad\mbox{and}\quad
\mbox{$\bar h$ is affine in $\{l(\hat x-1),\ldots,l\hat x\}$}\qquad\mbox{for $\hat
x\in\{1,\ldots,\frac{L}{l}\}$}.
\end{align}
Finally, let us optimize in the smaller scales and construct the competitor
$\bar h_*$ as the minimizer of $E$ among all functions $h$ on the lattice $\{0,\ldots,L\}$ that satisfy
\begin{align}\label{f11}
h(l\hat x)=\bar h(l\hat x)\qquad\mbox{for all $\hat
x\in\{0,\ldots,\frac{L}{l}\}$.}
\end{align}

\medskip

Furthermore, we define
\begin{align}\label{defn-cl}
c_L:=-\frac{\mathbb{E}\min E}{L\ln L},\quad\mbox{where the minimum
is taken over the class \eqref{ao49}} .
\end{align}
Using this notation and the previously defined competitor $ \bar h_* $ will show in the sequel that
\begin{align}\label{ao163}
\mathbb{E}E(\bar h_*)+c_{\frac{L}{l}} L\ln\frac{L}{l}+c_lL\ln l\lesssim L\ln^{\frac{1}{2}}\frac{L}{l},
\end{align}
which implies the superadditivity (up to lower order terms) property
\begin{align*}
c_L \ln L-c_l\ln l-c_\frac{L}{l}\ln\frac{L}{l}\gtrsim-\ln^{\frac{1}{2}}\frac{L}{l} .
\end{align*}
Before establishing \eqref{ao163} we will argue that this is enough to prove the existence of the limit $ \in [ 0 , \infty ] $ of $c_L$ as $L\uparrow\infty$.

\medskip

For convenience of notation, we will exchange the role of $l$ and $\frac{L}{l}$, and thus work with
\begin{align}\label{ao169}
c_L\ln L-c_l\ln l-c_{\frac{L}{l}}\ln\frac{L}{l}\gtrsim-\ln^{\frac{1}{2}}l.
\end{align}
Let $q\in\mathbb{N}$ and $1\le l'<l$ be such that $L=l^ql'$. Iterating
\eqref{ao169}, we obtain
\begin{align*}
c_L\ln L-qc_l\ln l-c_{l'}\ln l'\gtrsim -q\ln^\frac{1}{2}l.
\end{align*}
Since $\ln L\in(\ln l)[q, q+1)$, dividing by $q$ and sending $L\uparrow\infty$, we get
\begin{align*}
(\liminf_{L\uparrow\infty}c_L)\ln l-c_l\ln l\gtrsim-\ln^{\frac{1}{2}}l.
\end{align*}
Dividing further by $\ln l$ and sending now $l\uparrow\infty$, we conclude that
\begin{align*}
\liminf_{L\uparrow\infty}c_L\ge\limsup_{l\uparrow\infty}c_l ,
\end{align*}
as desired.

\medskip

Finally, here comes the argument for the core estimate \eqref{ao163}. We will deduced it from the sum of the following:
\begin{align}
&l\mathbb{E}\hat E(\hat h_*)=-c_{\frac{L}{l}}L\ln\frac{L}{l},&\mbox{(large-scales energy)}\label{f2}\\
&l\big(\mathbb{E}\hat E(\bar{\hat h}_*)-\mathbb{E}\hat E(\hat h_*)\big)\lesssim 
l\|\hat E(\bar{\hat h}_*)-\hat E(\hat h_*)\|_{\frac{3}{2}}\lesssim
L\ln^{\frac{1}{2}}\frac{L}{l},&\mbox{(error from binning)}\label{f3}\\
&\mathbb{E}E(\bar h)-l\mathbb{E}\hat E(\bar{\hat h}_*)\lesssim\|E(\bar h)-l\hat E(\bar{\hat h}_*)\|_{\frac{3}{2}}\lesssim
L\ln^{\frac{1}{4}}\frac{L}{l},&\mbox{(error from scaling)}\label{f6}\\
&\mathbb{E}E(\bar h_*)-\mathbb{E}E(\bar h)+c_lL\ln l\lesssim
\|E(\bar h_*)-E(\bar h)+c_lL\ln l\|_{\frac{3}{2}}\lesssim L.&\mbox{(small-scales energy)}\label{f5}
\end{align}

\medskip

Equality \eqref{f2} follows directly from the definitions of $ \hat h_* $ and $c_{\frac{L}{l}}$
after noticing that $\{\hat W_{\hat x}(\cdot)\}_{\hat x}$ is a collection of
independent Brownian motions (cf.~\eqref{f1}).

\medskip

Here comes the proof of \eqref{f3}: The estimate will be obtained as the sum of
the error on the noise
\begin{align}\label{f7}
l\|\sum_{\hat x=1}^{\frac{L}{l}}\hat W_{\hat x}(\bar{\hat h}_*(\hat x))-\hat W_{\hat
x}(\hat h_*(\hat x))\|_2\lesssim L
\end{align}
and the deterministic estimate on the Dirichlet energy
\begin{align}\label{f8}
l|\hat D(\bar{\hat h}_*)-\hat D(\hat h_*)|\lesssim L+(Ll)^{\frac{1}{2}}\hat D^{\frac{1}{2}}(\hat
h_*),
\end{align}
together with the bound given by Corollary~\ref{C:4}
\begin{align*}
\|\hat D(\hat h_*)\|_{\frac{3}{2}}\lesssim \frac{L}{l}\ln\frac{L}{l}.
\end{align*}

\medskip

The estimate \eqref{f7} is an instance of Lemma~\ref{L:12}. Indeed, defining
\begin{align*}
X_{\hat x}(\bar{\hat h}):=\sup\Big\{|\hat W_{\hat x}(\bar{\hat h}(\hat x))-
\hat W_{\hat x}(\hat h(\hat x))|:|\bar{\hat h}(\hat x)-\hat h(\hat x)|
\le\frac{1}{2}\;\mbox{for all}\;\hat x\in\{0,\ldots,\frac{L}{l}\}\Big\},
\end{align*}
we notice that
\begin{align*}
\big|\sum_{\hat x=1}^{\frac{L}{l}}\hat W_{\hat x}(\bar{\hat h}_*(x))-\hat W_{\hat x}(\hat h_*(x))\big|\le\sum_{\hat x=1}^{\frac{L}{l}}X_{\hat x}(\bar{\hat h}_*).
\end{align*}
For $\bar{\hat h}$ fixed, the random variables $X_{\hat x}(\bar{\hat h})$ are
independent, so that Lemma~\ref{L:12} (applied with $s=2,p=\infty,\alpha=0$ and
$L$ replace by $\frac{L}{l}$) allows to reduce the bound \eqref{f7} to the
simpler
\begin{align*}
\mbox{for any $\bar{\hat h}$ and any $\hat x$}\quad\|X_{\hat x}(\bar{\hat h})\|_2\lesssim 1,
\end{align*}
which is true by standard estimates on the maximum of Brownian motion.

\medskip

The error on the Dirichlet energy in \eqref{f8}, can be obtained
recall that is a norm $\hat D^\frac{1}{2}$, i.e.~satisfies the triangle inequality) as follows
\begin{align*}
|\hat D(\bar{\hat h}_*)-\hat D(\hat h_*)|&
=(\hat D^{\frac{1}{2}}(\bar{\hat h}_*)+\hat D^{\frac{1}{2}}(\hat h_*))
|\hat D^{\frac{1}{2}}(\bar{\hat h}_*)-\hat D^{\frac{1}{2}}(\hat h_*)|
\\&\leq (2\hat D^{\frac{1}{2}}(\hat h_*)+
|\hat D(\bar{\hat h}_*)-\hat D(\hat h_*)|^{\frac{1}{2}})\hat D^{\frac{1}{2}}(\bar{\hat h}_*-\hat h_*)\\
&\overset{\eqref{f10}}{\le}(2\hat D^{\frac{1}{2}}(\hat h_*)+|D(\bar{\hat h}_*)-D(\hat
h_*)|^{\frac{1}{2}})(\frac{L}{2l})^{\frac{1}{2}}.
\end{align*}
Applying Young's inequality on the second term $|D(\bar{\hat h})-D(\hat
h_*)|^{\frac{1}{2}}(\frac{L}{l})^{\frac{1}{2}}$ and multiplying by $l$, we obtain \eqref{f8}.

\medskip

Here comes the argument for \eqref{f6}: Because of the choice in \eqref{f9}, in the difference $E(\bar h)-l\hat E(\bar{\hat h}_*)$, the Dirichlet energies simplify and we get
\begin{align*}
|E(\bar h)-l\hat E(\bar{\hat h}_*)|=\Big|\sum_{\hat x=1}^{\frac{L}{l}}
\Big(\sum_{x=l(\hat x-1)+1}^{l\hat x}W_x(\bar h(x))-W_x(\bar h(l\hat x))\Big)\Big|
\le\sum_{\hat x=1}^{\frac{L}{l}}lX_{\hat x}(\bar{\hat h}_*) ,
\end{align*}
where (for every $\bar{\hat h}$ and $\hat x$) we defined the errors according to
\begin{align*}
X_{\hat x}(\bar{\hat h}):=\frac{1}{l}\Big|\sum_{x=l(\hat x-1)+1}^{l\hat
x}W_x(\bar h(x))-W_x(\bar h(l\hat x))\Big| ,
\end{align*}
with $\bar h$ constructed from every $\bar{\hat h}$ as in \eqref{f9}. For
every (fixed, deterministic) $\bar{\hat h}$, we have
\begin{align*}
&\|X_{\hat x}(\bar{\hat h})\|_2^2
\lesssim\frac{1}{l^2}\mathbb{E}\Big|\sum_{x=l(\hat
x-1)+1}^{l\hat x}W_x(\bar h(x))-W_x(\bar h(l\hat x))\Big|^2
&\mbox{(using Gaussianity)}\\
&=\frac{1}{l^2}\sum_{x=l(\hat x-1)+1}^{l\hat x}|\bar h(x)-\bar h(l\hat x)|
&\mbox{(independence)}\\
&\overset{\eqref{f9}}{=}|\bar{\hat h}(\hat x-1)-\bar{\hat h}(\hat x)|\sum_{x=l(\hat x-1)+1}^{l\hat x}
\frac{|x-l\hat x|}{l^2}\le\frac{1}{2}|\bar{\hat h}(\hat x-1)-\bar{\hat h}(\hat x)|,
\end{align*}
so that
\begin{align*}
\Big(\frac{l}{L}\sum_{\hat x=1}^{\frac{L}{l}}\|X_{\hat x}(\bar{\hat
h})\|_2^2\Big)^\frac{1}{2}\lesssim
\Big(\frac{l}{L}\sum_{\hat x=1}^{\frac{L}{l}}|\bar{\hat h}(\hat x-1)-\bar{\hat
h}(\hat x)|\Big)^\frac{1}{2}
\le\Big(\frac{2D(\bar{\hat h})}{L/l}\Big)^\frac{1}{4}.
\end{align*}
An application of Lemma~\ref{L:12} with $s=p=2,\alpha=\frac{1}{4}$ and
$\frac{L}{l}$ in place of $L$ then gives
\begin{align*}
\|\sum_{\hat x=1}^{\frac{L}{l}}X_{\hat x}(\bar{\hat h}_*)\|_{\frac{3}{2}}\lesssim
\frac{L}{l}\ln^{\frac{1}{4}}\frac{L}{l}.	
\end{align*}

\medskip

Finally, here comes the proof of \eqref{f5}: Denoting by 
\begin{align*}
E_{\hat x}(h)=\frac{1}{2}\sum_{x=l(\hat x-1)+1}^{l\hat x}(h(x)-h(x-1))^2-\sum_{x=l(\hat x-1)+1}^{l\hat x} W_x(h(x))
\end{align*}
the restriction of the energy to the interval $\{l(\hat x-1)+1,\ldots,l\hat x\}$, one can rewrite \eqref{f5} as
\begin{align}
\|\sum_{\hat x=1}^{\frac{L}{l}}(E_{\hat x}(\bar
h_*)-E_{\hat x}(\bar h)+c_ll\ln
l)\|_{\frac{3}{2}}\lesssim L.
\end{align}
This can be proven by a further application of Lemma~\ref{L:12} with the choice
\begin{align*}
	X_{\hat x}(\bar{\hat h})&=|E_{\hat x}(\bar h_*)-E_{\hat x}(\bar
	h)+c_ll\ln l|.
\end{align*}
with $\bar h$ and $\bar h_*$ constructed from every $\bar{\hat h}$ using
\eqref{f9} and \eqref{f11}. By the shear invariance of the energy $E$ (see
\eqref{ao134}) and \eqref{defn-cl} we can
rewrite
\begin{align*}
X_{\hat x}(\bar{\hat h})&=\big|E_{\hat x}(\bar h_*)-
E_{\hat x}(\bar h)-\mathbb{E}\big(E_{\hat x}(\bar h_*)
-E_{\hat x}(\bar h)\big)\big|,
\end{align*}
and observe that the law of $X_{\hat x}(\bar{\hat h})$ does not depend on
$\bar{\hat h}$. Because of the concentration result of Lemma~\ref{L:10} applied
to each subsystem of size $l$, we know that for each deterministic $\bar{\hat h}$
\begin{align*}
\mbox{$\|X_{\hat x}(\bar{\hat h})\|_{\frac{3}{2}}\lesssim l$ and $\{X_{\hat x}(\bar{\hat h})\}_{\hat x}$ are independent.}
\end{align*}
This is enough to apply Lemma~\ref{L:12} with $s=\frac{3}{2},p=\infty,\alpha=0$ and conclude.
\qed

\medskip

{\sc Proof of Lemma~\ref{L:10}.}
Because of \eqref{ao39} in Lemma~\ref{L:8}, to prove \eqref{ao138} it is enough to show that
\begin{align*}
\ln\mathbb{P}(|\min \frac{E}{L}-\mathbb{E}\min \frac{E}{L}|\geq
\nu)\lesssim-\nu^{\frac{3}{2}}\qquad\mbox{for all}\;\nu\gg1.
\end{align*}
We start by claiming that
\begin{align}\label{ao143}
\ln\mathbb{P}(|\min \frac{E(h)}{L}-\mathbb{E}\min_{\frac{M(h)}{L^2}\leq
\sqrt{\nu}}\frac{E(h)}{L}|\geq \nu)\lesssim-\nu^{\frac{3}{2}}.
\end{align}
Indeed, by the union bound, \eqref{ao143} amounts to estimating the probability of the two events
\begin{align*}
\Big\{\frac{M(h_*)}{L^2}\geq\sqrt{\nu}\Big\}\qquad\mbox{and}\qquad\Big\{|\min_{\frac{M(h)}{L^2}\leq\sqrt{\nu}}
\frac{E(h)}{L}-\mathbb{E}\min_{\frac{M(h)}{L^2}\leq\sqrt{\nu}}\frac{E(h)}{L}|\geq
\nu\Big\}.
\end{align*}
This is done by using Proposition~\ref{P:1} (together with \eqref{ao124} in
Lemma~\ref{L:8}) and the constrained concentration \eqref{c1} in Lemma~\ref{L:4}.

\medskip

Having established \eqref{ao143}, it is enough to apply the following fact
(see \cite[Claim 5.6]{Peled}): If for a random variable $X$ and $s\ge1$
\begin{align}\label{ao141}
\inf_a\ln\mathbb{P}(|X-a|\geq \nu)\leq-\nu^s\qquad \nu\gg1,
\end{align}
then
\begin{align}\label{ao142}
\ln\mathbb{P}(|X-\mathbb{E}X|\geq \nu)\lesssim-\nu^s\qquad \nu\gg1.
\end{align}
For the convenience of the reader, we give a proof in the appendix.\qed

\medskip

{\sc Proof of Lemma~\ref{L:12}.} By monotonicity in $p$, one can assume $p=s'$
the dual exponent of $s$. Let us show that for a fixed $\bar h$, we have
\begin{align}\label{f13}
\ln\mathbb{P}\big(\frac{1}{L}\sum_{x=1}^LX_x(\bar h)\ge\nu\ln^\alpha L\big)\le
L-L(\frac{D(\bar h)}{L\ln L})^{-s\alpha}\nu^s.
\end{align}
It follows from applying H\"older's inequality
\begin{align*}
\frac{1}{L}\sum_{x=1}^LX_x(\bar h)
&\le\Big(\frac{1}{L}\sum_{x=1}^L(\frac{X_x(\bar h)}{\|X_x(\bar h)\|_s})^s\Big)^{\frac{1}{s}}
\Big(\frac{1}{L}\sum_{x=1}^L\|X_x(\bar h)\|_s^{s'}\Big)^\frac{1}{s'}\\
&\overset{\eqref{f12}}{\le}\Big(\frac{1}{L}\sum_{x=1}^L(\frac{X_x(\bar h)}{\|X_x(\bar h)\|_s})^s\Big)^{\frac{1}{s}}
\big(\frac{D(\bar h)}{L}\big)^\alpha,
\end{align*}
together with $ s $-exponential Chebyshev's inequality and the independence in $\eqref{ao189}$
\begin{align*}
&\ln\mathbb{P}\big(\frac{1}{L}\sum_{x=1}^LX_x(\bar h)\ge\nu\ln^\alpha L\big)
\le\ln\mathbb{P}\Big(\big(\frac{1}{L}\sum_{x=1}^L(\frac{X_x(\bar h)}{\|X_x(\bar h)\|_s})^s\big)^{\frac{1}{s}}
\big(\frac{D(\bar h)}{L}\big)^\alpha\ge\nu\ln^\alpha L\Big)\\
&=\ln\mathbb{P}\Big(\exp\big(\sum_{x=1}^L(\frac{X_x(\bar h)}{\|X_x(\bar h)\|_s})^s\big)\ge
\exp\big(L(\frac{D(\bar h)}{L\ln L})^{-s\alpha}\nu^s\big)\Big)\\
&\overset{\text{(Chebyshev})}{\le}-L(\frac{D(\bar h)}{L\ln L})^{-s\alpha}\nu^s+
\ln\mathbb{E}\exp\big(\sum_{x=1}^L(\frac{X_x(\bar h)}{\|X_x(\bar h)\|_s})^s\big) \\
&\overset{\text{(independence})}{\le}-L(\frac{D(\bar h)}{L\ln L})^{-s\alpha}\nu^s+
\sum_{x=1}^L\ln\mathbb{E}\exp (\frac{X_x(\bar h)}{\|X_x(\bar h)\|_s})^s 
\overset{\eqref{ao189}}{\le}-L(\frac{D(\bar h)}{L\ln L})^{-s\alpha}\nu^s+L.
\end{align*}

\medskip

For constants $\hat D_1, \hat D_2>0$, let us denote by $\mathcal{H}_{\hat
D_1,\hat D_2}$ the set of functions $\bar h$ which satisfies
\begin{align*}
D(\bar h)\le\hat D_2L\ln L\qquad\mbox{and}\qquad D(\bar h_\rho)\le
\hat D_1\rho^2L\;\mbox{for every dyadic scale $1\le\rho<L$.}
\end{align*}
Distinguishing whether $\bar h_*$ belongs to $\mathcal{H}_{\hat D_1,\hat D_2}$ 
or not, we get the following inclusion
\begin{align*}
\Big\{\frac{1}{L}\sum_{x=1}^LX_x(\bar h_*)\ge\nu\ln^\alpha L\Big\}
&\subset\Big\{\sup_{\bar h\in\mathcal{H}_{\hat D_1,\hat D_2}}
\frac{1}{L}\sum_{x=1}^LX_x(\bar h)\ge\nu\ln^\alpha L\Big\}\\
&\cup\big\{D(\bar h_{*,\rho})>\hat D_1\rho^2L\;\mbox{for a dyadic $1\le\rho<L$}\big\}\\
&\cup\big\{D(\bar h_*)>\hat D_2L\ln L\big\},
\end{align*}
so that it is enough to estimate the probability of each of the three events on the r.h.s..

\medskip

Let us estimate the first one using \eqref{f13} and a union bound:
\begin{align*}
\ln\mathbb{P}\big(\sup_{\bar h\in\mathcal{H}_{\hat D_1,\hat D_2}}
\frac{1}{L}\sum_{x=1}^LX_x(\bar h)\ge\nu\big)
&\le\ln(\#\mathcal{H}_{\hat D_1,\hat D_2})+L-L
\big(\sup_{\bar h\in\mathcal{H}_{\hat D_1,\hat D_2}}\frac{D(\bar h)}{L\ln L}\big)^{-s\alpha}\nu^s\\
&\overset{\text{Lemma~\ref{L:2}}}{\le} cL\ln\hat D_1+L-L\hat D_2^{-s\alpha}\nu^s,
\end{align*}
for some universal constant $c$. The probability of the second term
can be estimated using Theorem~\ref{T} (with $p=2$) that in view of \eqref{ao124} and the union bound implies
\begin{align*}
&\ln\mathbb{P}\big(D(\bar h_{*,\rho})>\hat D_1\rho^2L\;
\mbox{for a dyadic $1\le\rho<L$}\big)
\lesssim \ln \sum_{ \substack{1\le\rho<L\\\text{dyadic}} } e^{ - c \rho^2 \hat D_1 }
\lesssim - \hat D_1 \qquad\mbox{for}\;\hat D_1\gg1 .
\end{align*}
Lastly, by the same arguments, Corollary~\ref{C:4} shows that the third probability is
\begin{align*}
&\ln\mathbb{P}\big(D(\bar h_*)>\hat D_2L\ln L\big)\lesssim-\hat D_2^{\frac{3}{2}}
\qquad\mbox{for}\;\hat D_2\gg 1.
\end{align*}

\medskip

With the final choice
\begin{align*}
\frac{1}{t}:=\frac{1}{s}+\frac{2}{3}\alpha,\qquad \hat D_1:=\nu^t\qquad
\mbox{and $\hat D_2$ such that}\quad \hat D_2^{ \frac{3}{2} }= \hat D_2^{ - s \alpha } \nu^{ s} ,
\end{align*}
i.e.~$ \hat D_2^{ \frac{3}{2} } = \nu^{ t } $,
the three quantities above are $\lesssim-\nu^t$ for $\nu\gg 1$ and we get
\begin{align*}
\ln\mathbb{P}\big(\frac{1}{L}\sum_{x=1}^LX_x(\bar h_*)\ge\nu\ln^\alpha L\big)
\lesssim-\nu^t\qquad\mbox{for}\;\nu\gg1.
\end{align*}
Because of Lemma~\ref{L:8}, this is enough to conclude.
\qed

\medskip

{\sc Proof of Lemma~\ref{L:13}.}
Let us use the notation from the proof of Lemma~\ref{L:5}, namely
\begin{align*}
	\hat W_1:=\sup_{\frac{D(h)}{L}\leq 1}\frac{W(h)}{L}.
\end{align*}
Recall the formulation \eqref{ao56} of the problem
\begin{align*}
	\frac{D(h_*)}{L}=\mbox{argmin}_{\hat D}\{\hat D-\sup_{\frac{D(h)}{L}\le\hat
	D}\frac{W(h)}{L}\}.
\end{align*}
If we could substitute the field term with its expectation, using the scaling
\eqref{ao63} we would have the convex problem
\begin{align*}
	\text{argmin}_{\hat D}\{\hat D-\hat D^{\frac{1}{4}}\mathbb{E}\hat W_1\}
\end{align*}
which admits the unique solution
\begin{align}\label{ao201}
\hat D_*:=(\frac{1}{4})^{\frac{4}{3}}\mathbb{E}^{\frac{4}{3}}\hat W_1.
\end{align}
It is then convenient by adding and subtracting the expectation to rewrite the
problem as
\begin{align}\label{ao192}
\mbox{argmin}_{\hat D}\{\hat D-\hat D^{\frac{1}{4}}\mathbb{E}\hat
W_1+\frac{3}{4}(\frac{1}{4})^{\frac{1}{3}}\mathbb{E}^{\frac{4}{3}}\hat
W_1-X_{\hat D}\} ,
\end{align}
where $ X_{ \hat D } $ is defined by
\begin{align}\label{fo02}
X_{\hat D}:=\sup_{\frac{D(h)}{L}\le\hat D}\frac{W(h)}{L}
-\mathbb{E}\sup_{\frac{D(h)}{L}\le\hat D}\frac{W(h)}{L}
\end{align}
for which we have the concentration bound \eqref{c2} in Lemma
\ref{L:4}, which together with the scaling \eqref{ao63} implies
\begin{align}\label{ao190}
\|X_{\hat D}\|_2\lesssim \hat D^{\frac{1}{4}}.
\end{align}
The constant
$\frac{3}{4}(\frac{1}{4})^{\frac{1}{3}}\mathbb{E}^{\frac{4}{3}}\hat
W_1$ has been added so that the minimum of the deterministic part
is zero. 

\medskip

To control the fluctuations $X_{\hat D}$, we use the strong convexity
of the deterministic problem which can be rewritten as
\begin{align*}
	\hat D-\hat D^{\frac{1}{4}}\mathbb{E}\hat
	W_1+\frac{3}{4}(\frac{1}{4})^{\frac{1}{3}} \mathbb{E}^{\frac{4}{3}} \hat W_1 = \mathbb{E}^{\frac{4}{3}} \hat
W_1 \, p\Big((\frac{\hat D}{\mathbb{E}^{\frac{4}{3}}\hat
W_1})^{\frac{1}{4}}\Big)
\end{align*}
where we introduced
\begin{align*}
p(t)=t^4-t+\frac{3}{4}(\frac{1}{4})^{\frac{1}{3}}\
=(t-(\frac{1}{4})^{\frac{1}{3}})^2\times(\mbox{positive polynomial of
degree}~2).
\end{align*}
As a consequence
\begin{align}
\hat D-\hat D^{\frac{1}{4}}\mathbb{E}\hat
W_1+\frac{3}{4}(\frac{1}{4})^{\frac{1}{3}} \mathbb{E}^{\frac{4}{3}} \hat
W_1&\overset{\eqref{ao201}}{\gtrsim}\mathbb{E}^{\frac{4}{3}} \hat W_1 \Big(
(\frac{\hat D}{ \mathbb{E}^{\frac{4}{3}} \hat
W_1 } )^{\frac{1}{4}}-(\frac{\hat D_*}{\mathbb{E}^{\frac{4}{3}}\hat
W_1})^{\frac{1}{4}}\Big)^2\nonumber\\
   &\gtrsim\frac{\mathbb{E}^{\frac{2}{3}}\hat W_1}{\hat D^{\frac{3}{2}}+\hat
		   D_*^{\frac{3}{2}}}(\hat D-\hat D_*)^2\overset{\eqref{ao108}}{\gtrsim}\frac{\ln^{\frac{1}{2}}L}{\hat
   D^{\frac{3}{2}}+\hat D_*^{\frac{3}{2}}}(\hat D-\hat D_*)^2\label{ao191}
\end{align}
where in the last line we used that
\begin{align*}
|x-y|\gtrsim\frac{|x^4-y^4|}{x^3+y^3}\quad\mbox{and}\quad(x+y)^2\lesssim x^2+y^2
\end{align*}
in order to let the error $\hat D-\hat D_*$ explicitly appear.

\medskip

Next, we will argue that for $\mu\gg1$
\begin{align}\label{ao208}
(\frac{D(h_*)}{L}-\hat D_*)^2\lesssim \mu^{\frac{3}{2}}(\ln L)(-X_{\hat
D_*}+\sup_{|\hat D-\hat D_*|\le\mu\ln L}X_{\hat D})
\quad\mbox{if} \quad |\frac{D(h_*)}{L}-\hat D_*|\le\mu\ln L.
\end{align}
To this end, we use that $\hat D_*$ is a competitor for the problem \eqref{ao192}, so that
\begin{align*}
\frac{D(h_*)}{L}-(\frac{D(h_*)}{L})^{\frac{1}{4}}\mathbb{E}\hat
W_1+\frac{3}{4}(\frac{1}{4})^{\frac{1}{3}}\mathbb{E}^{\frac{4}{3}}\hat W_1
\le-X_{\hat D_*}+X_{\frac{D(h_*)}{L}}.
\end{align*}
To deduce \eqref{ao208}, we assume $|\frac{D(h_*)}{L}-\hat D_*|\le \mu\ln L$ and estimate the r.~h.~s.~from above by
\begin{align*}
-X_{\hat D_*}+X_{\frac{D(h_*)}{L}}\le-X_{\hat D_*}+\sup_{|\hat D-\hat
D_*|\le\mu\ln L}X_{\hat D}
\end{align*}
and the l.~h.~s.~from below using \eqref{ao191}
\begin{align*}
&\frac{D(h_*)}{L}-(\frac{D(h_*)}{L})^{\frac{1}{4}}\mathbb{E}\hat
W_1+\frac{3}{4}(\frac{1}{4})^{\frac{1}{3}}\mathbb{E}^{\frac{4}{3}}\hat W_1\\
&\gtrsim\frac{\ln^{\frac{1}{2}}L}{(\frac{D(h_*)}{L})^{\frac{3}{2}}+\hat D_*^{\frac{3}{2}}}(\frac{D(h_*)}{L}-\hat D_*)^2
\gtrsim\frac{1}{\mu^{\frac{3}{2}}\ln L}(\frac{D(h_*)}{L}-\hat D_*)^2,
\end{align*}
where in the last inequality we used that the Orlicz norm bounds the
expectation so that
\begin{align}\label{ao201cont}
\hat D_*\overset{\eqref{ao201}}{\lesssim}\mathbb{E}^{\frac{4}{3}}\hat W_1
\overset{\eqref{ao82}}{\lesssim}\ln L\lesssim\mu\ln L\qquad\mbox{provided}~\mu\gg1.
\end{align}
The combination of the last two statements implies \eqref{ao208}.

\medskip

On the basis of \eqref{ao208}, we will prove the following uniform bound on the fluctuations
\begin{align}\label{ao204}
\|\sup_{|\hat D-\hat D_*|\le\mu\ln L}X_{\hat
D}\|_2\lesssim(\mu\ln L)^{\frac{1}{4}}\ln^{\frac{1}{2}}\ln
L\qquad\mbox{for}\;\mu\gg1
\end{align}
below.
Assuming this, let us conclude. Since from \eqref{ao190} and
\eqref{ao201cont}, $\|X_{\hat
D_*}\|_2\lesssim\ln^{\frac{1}{4}}L$, taking the Orlicz norm
of \eqref{ao208} and using that $\|X\|_4=\|X^2\|_2^{\frac{1}{2}}$ for any random
variable $X$ together with the bound \eqref{ao204}, we obtain
\begin{align}
&\|(\frac{D(h_*)}{L}-\hat D_*)I\big(|\frac{D(h_*)}{L}-\hat D_*|\le\mu\ln
L\big)\|_4\nonumber\\
&=\|(\frac{D(h_*)}{L}-\hat D_*)^2I\big(|\frac{D(h_*)}{L}-\hat D_*|\le\mu\ln
L\big)\|_2^{\frac{1}{2}}
\lesssim
\mu^{\frac{7}{8}}(\ln^{\frac{5}{8}}L)\ln^{\frac{1}{4}}\ln
L. \label{fo01}
\end{align}
We now express \eqref{fo01} this in term of
\begin{align*}
X:= \frac{D(h_*)}{L}-\hat D_* ,
\end{align*}
for which we have
\begin{align*}
\|X \|_{\frac{3}{2}}\le\|\frac{D(h_*)}{L}\|_{\frac{3}{2}}+|\hat D_*|
\overset{\eqref{ao82},\eqref{ao201cont}}{\lesssim}\ln L,
\end{align*}
so that \eqref{fo01} implies
\begin{align*}
\|XI\big(|X|\le
\mu\| X \|_{\frac{3}{2}}\big)\|_4
\lesssim\mu^{\frac{7}{8}}(\ln^{\frac{5}{8}}L)\ln^{\frac{1}{4}}\ln
L\qquad\mbox{for}\;\mu\gg1.
\end{align*}
We now appeal to the following fact on the Orlicz norms (see
Lemma~\ref{L:15}.(iv) in the appendix)
\begin{align}
&\mbox{if}\qquad\|X~I(X\le\mu\|X\|_s)\|_t\le\mu^\alpha A
\quad\mbox{for}\quad\mu\gg1,\label{ao218}\\
&\mbox{then}\qquad\|X\|_1\lesssim A\qquad
\mbox{provided}\quad \frac{st}{s+\alpha t} \ge 1.\label{ao212}
\end{align}
With the choice $ X = \frac{D(h_*)}{L}-\hat D_*,
s=\frac{3}{2},t=4,\alpha=\frac{7}{8}$ and $A\sim(\ln^{\frac{5}{8}}L)\ln^{\frac{1}{4}}\ln
L$, it implies \eqref{ao171}.

\medskip

It remains to prove \eqref{ao204}, which due to \eqref{ao201cont}, is equivalent
to
\begin{align}\label{ao211}
\|\sup_{0\le\hat D\le\mu\ln L}X_{\hat
D}\|_2\lesssim(\mu\ln L)^{\frac{1}{4}}\ln^{\frac{1}{2}}\ln
L\qquad\mbox{for}\;\mu\gg1.
\end{align}
We need the following fact (see Lemma~\ref{L:15}.(ii) in the appendix): For a finite set of random variables $\{Z_i\}_{i \in I}$
\begin{align}\label{ao206}
\|\sup_{i \in I}Z_i\|_2\lesssim\ln^{\frac{1}{2}}(\#I)\sup_{i\in I}\|Z_i\|_2 .
\end{align}
To infer \eqref{ao211} from \eqref{ao206}, we estimate a bin via
\begin{align}\label{ao205}
\|\sup_{\hat D'\in[\hat D-\Delta,\hat D]} X_{\hat D'} \|_2\lesssim(\mu\ln
L)^{\frac{1}{4}}\qquad\mbox{provided}\;\Delta\leq\mu^{\frac{1}{4}}(\ln^{-\frac{1}{2}}L)\hat
D^{\frac{3}{4}},
\end{align}
which we will establish below. We then construct a sequence of intervals of the form
$[\hat D_{n-1}, \hat D_n]$ with $(\hat D_n)_{n\ge 0}$ satisfying the
implicit relation
\begin{align}\label{ao210}
\hat D_n-\hat D_{n-1}=\mu^{\frac{1}{4}}(\ln^{-\frac{1}{2}}L)\hat D_n^{\frac{3}{4}}
\qquad\hat D_0=0
\end{align}
such that \eqref{ao205} applies for $\hat D'\in[\hat D_{n-1},\hat D_n]$. We will
argue that
\begin{align}\label{ao209}
\hat D_n\geq\mu\ln L\qquad\mbox{for}\;n\gg\ln^{\frac{3}{4}}L.
\end{align}
This implies that with $\lesssim\ln^{\frac{3}{4}}L$ intervals it is possible to cover the set $[ 0, \mu\ln L ]$, so that from \eqref{ao205} and \eqref{ao206} we deduce
\eqref{ao211}.

\medskip

Here comes the argument for \eqref{ao205}. It is sufficient to control the supremum from above,
since it is bounded from below by $X_{\hat D}$, for which we can apply
\eqref{ao190}. If $\hat D'\in[\hat D-\Delta,\hat D]$, then using
the monotonicity of the map
\begin{align}\label{ao202}
\hat D\mapsto\sup_{\frac{D(h)}{L}\leq\hat D}\frac{W(h)}{L}
\end{align}
we have
\begin{align*}
X_{\hat D'}
&\le\sup_{\frac{D(h)}{L}\le\hat
D}\frac{W(h)}{L}-\mathbb{E}\sup_{\frac{D(h)}{L}\le\hat
D-\Delta}\frac{W(h)}{L}\\
&=X_{\hat D}+\Big(\mathbb{E}\sup_{\frac{D(h)}{L}\le\hat
D}\frac{W(h)}{L}-\mathbb{E}\sup_{\frac{D(h)}{L}\le\hat
D-\Delta}\frac{W(h)}{L}\Big).
\end{align*}
From \eqref{ao190} we know that
\begin{align*}
\|X_{\hat D}\|_2\lesssim \hat D^{\frac{1}{4}}\le(\mu\ln
L)^{\frac{1}{4}} . 
\end{align*}
Therefore, it is enough to show that the second term is $\lesssim(\mu\ln
L)^{\frac{1}{4}}$.
Because of the scaling in \eqref{ao63} it satisfies
\begin{align*}
&\mathbb{E}\sup_{\frac{D(h)}{L}\le\hat
D}\frac{W(h)}{L}-\mathbb{E}\sup_{\frac{D(h)}{L}\le\hat
D-\Delta}\frac{W(h)}{L}=\big(\hat D^{\frac{1}{4}}-(\hat D-\Delta)^{\frac{1}{4}}\big)\mathbb{E}\hat W_1
\\&
=\big(1-(1-\frac{\Delta}{\hat D})^{\frac{1}{4}}\big)\hat
D^{\frac{1}{4}}\mathbb{E}\hat W_1\lesssim\frac{\Delta}{\hat D}\hat
D^{\frac{1}{4}}\ln^{\frac{3}{4}}L ,
\end{align*}
With the condition $\Delta\le\mu^{\frac{1}{4}}(\ln^{-\frac{1}{2}}L)\hat
D^{\frac{3}{4}}$, the last term is exactly bounded by $ (\mu\ln L)^{\frac{1}{4}} $.
We used that $\hat D\ge\Delta$ and that $ 1-(1-x)^{\frac{1}{4}} \leq x$ for $x\in[0,1]$.

\medskip

Here comes the argument for the lower bound in \eqref{ao209}. To this end, we like to think of the iteration formula \eqref{ao210} as discrete version of the ODE
\begin{align*}
f'(x)=\mu^{\frac{1}{4}}(\ln^{-\frac{1}{2}}L)f^{\frac{3}{4}}(x)\qquad f(0)=0 ,
\end{align*}
which is solved by the function $ f(x)=( \mu \ln^{-2}L) (\frac{x}{4})^4$. By comparison this implies a lower bound on $ \hat D_n $. Indeed, the sequence $(f(n))_{n\in\mathbb{N}}$ is a subsolution for
the problem in \eqref{ao210}, since
\begin{align*}
f(n)-f(n-1)=\int_{n-1}^ndx~f'(x)=
\int_{n-1}^n dx ~ \mu^{\frac{1}{4}}(\ln^{-\frac{1}{2}}L)f^{\frac{3}{4}}(x)
\le\mu^{\frac{1}{4}}(\ln^{-\frac{1}{2}}L)f^{\frac{3}{4}}(n) ,
\end{align*}
and thus by comparison principle
\begin{align*}
\hat D_n\ge f(n) \geq \mu\ln L\qquad\mbox{for}\;n\gg\ln^{\frac{3}{4}}L.
\end{align*}

\medskip

We can now derive \eqref{ao200} from \eqref{ao171}. Starting from
\begin{align*}
\frac{W(h_*)}{L}=\sup_{\frac{D(h)}{L}\le\frac{D(h_*)}{L}}\frac{W(h)}{L},
\end{align*}
if we could substitute $\frac{D(h_*)}{L}$ with $\hat D_*$ (cf.~\eqref{ao171}), the
r.~h.~s.~would be
\begin{align*}
\sup_{\frac{D(h)}{L}\le\hat D_*}\frac{W(h)}{L}
\end{align*}
the expectation of which is exactly
\begin{align}\label{ao201bis}
	\hat W_*:=\mathbb{E}\sup_{\frac{D(h)}{L}\le\hat
	D_*}\frac{W(h)}{L}\overset{\eqref{ao63}}{=}
	\hat D_*^{\frac{1}{4}}\mathbb{E}\hat
	W_1\overset{\eqref{ao201}}{=}(\frac{1}{4})^{\frac{1}{3}}\mathbb{E}^{\frac{4}{3}}\hat
	W_1.
\end{align}

\medskip

To quantify the error made in these two steps, we will show 
that
\begin{align*}
\ln\mathbb{P}\Big(|\frac{W(h_*)}{L}-\hat W_*|\ge\nu\Big)\lesssim-\frac{\nu}{(\ln^{\frac{5}{8}}L)\ln^{\frac{1}{4}}\ln
L}\qquad\mbox{for}\;\nu\gg(\ln^{\frac{5}{8}}L)\ln^{\frac{1}{4}}\ln L.
\end{align*}
We will assume that $ \frac{W(h_*)}{L} \geq \hat W_* $, since the other case is similar. Using the
monotonicity of the sup in \eqref{ao202}, we have the inclusion
\begin{align*}
\Big\{\frac{W(h_*)}{L}-\hat W_*\ge\nu\Big\}\subset
\Big\{\frac{D(h_*)}{L}-\hat D_*\ge\frac{\nu}{2}\Big\}\cup
\Big\{\sup_{\frac{D(h)}{L}\le\hat D_*+\frac{\nu}{2}}\frac{W(h)}{L}-\hat
W_*\ge\nu\Big\}.
\end{align*}
Because of \eqref{ao171} (and \eqref{ao124} in Lemma \ref{L:8}), the probability
of the first event on the r.~h.~s.~satisfies the bound
\begin{align*}
	\ln\mathbb{P}\Big(\frac{D(h_*)}{L}-\hat
		D_*\ge\frac{\nu}{2}\Big)\lesssim-\frac{\nu}{(\ln^{\frac{5}{8}}L)\ln^{\frac{1}{4}}\ln
L}\qquad\mbox{for}\;\nu\gg(\ln^{\frac{5}{8}}L)\ln^{\frac{1}{4}}\ln L.
\end{align*}
Regarding the second event, we will prove below that its probability is of lower
order since
\begin{align}\label{ao203}
\ln\mathbb{P}\Big(\sup_{\frac{D(h)}{L}\le\hat D_*+\frac{\nu}{2}}\frac{W(h)}{L}-\hat
W_*\ge\nu\Big)\lesssim-\frac{\nu}{\ln^{\frac{1}{2}}L}\qquad
\mbox{for}\;\nu\gg\ln^{\frac{1}{2}}L.
\end{align}
Therefore, together with \eqref{ao39}, this proves \eqref{ao200}. 

\medskip

Here comes the proof of \eqref{ao203}. Adding and subtracting the expected value and using
the scaling in \eqref{ao63}, we have
\begin{align*}
\sup_{\frac{D(h)}{L}\le\hat D_*+\frac{\nu}{2}}\frac{W(h)}{L}-\hat W_*=
X_{\hat D_*+\frac{\nu}{2}}+\Big(\big(\frac{\hat D_*+\frac{\nu}{2}}{\hat
D_*}\big)^{\frac{1}{4}}\hat W_*-\hat W_*\Big) ,
\end{align*}
(recall the definition of $ X_{ \hat D } $ in \eqref{fo02}).
By the concavity of the function $x\mapsto(1+x)^{\frac{1}{4}}-1$ and the relation
$4\hat D_*=\hat W_*$ (cf.~\eqref{ao201},\eqref{ao201bis}), we can estimate the term in parenthesis as
\begin{align*}
\big(\frac{\hat D_*+\frac{\nu}{2}}{\hat D_*}\big)^{\frac{1}{4}}\hat W_*-\hat W_*
\le\frac{\nu}{2}\frac{\hat W_*}{4 \hat D_*}=\frac{\nu}{2}.
\end{align*}
This finally leads to
\begin{align*}
\ln\mathbb{P}\Big(\sup_{\frac{D(h)}{L}\le\hat
D_*+\frac{\nu}{2}}\frac{W(h)}{L}-\hat W_*\ge\nu\Big)
\le\ln\mathbb{P}\Big(X_{\hat D_*+\frac{\nu}{2}}\ge\frac{\nu}{2}\Big)
\overset{\eqref{ao190},\eqref{ao124}}{\lesssim}-\frac{\nu^2}{\hat
D_*^{\frac{1}{2}}+\nu^{\frac{1}{2}}}\lesssim-\frac{\nu}{\ln^{\frac{1}{2}}L}
\end{align*}
for $\nu\gg\ln^{\frac{1}{2}}L$.
\qed

\medskip

{\sc Proof of Theorem~\ref{T:2}.} We first prove the convergence in the Orlicz
norm. Because of Lemmas~\ref{L:9} and~\ref{L:10},
we have already established the existence of the limit
\begin{align*}
\lim_{L\to\infty}\frac{\min E}{L\ln L}
\end{align*}
with respect to $\|\cdot\|_{\frac{3}{2}}$ (and therefore in $\|\cdot\|_1$). Since
\begin{align*}
\frac{\min E}{L\ln L}=\frac{D(h_*)}{L\ln L}-\frac{W(h_*)}{L\ln L}
\end{align*}
from the concentrations in \eqref{ao171} and \eqref{ao200} we have (with respect to $ \| \cdot \|_1 $)
\begin{align*}
\lim_{L\to\infty}\frac{\min E}{L\ln
L}=\lim_{L\to\infty}\big((\frac{1}{4})^{\frac{4}{3}}-(\frac{1}{4})^{\frac{1}{3}}\big)
(\frac{\mathbb{E}\hat W_1}{\ln^{\frac{3}{4}}L})^{\frac{4}{3}}.
\end{align*}
Together with the concentration \eqref{c2}, this proves the existence of the
limit in \eqref{ao137}. Equipped with this, thanks to \eqref{ao171}, we
learn that the limit in \eqref{ao136} exists with respect to $\|\cdot\|_1$.
Finally, due to the bound \eqref{ao82} of $D(h_*)$ with respect to
$\|\cdot\|_{\frac{3}{2}}$ and the Hölder inequality for Orlicz norms
(see Lemma~\ref{L:15}.(i) in the appendix)
\begin{align*}
\|X\|_s\le\|X\|_r^{\lambda\frac{r}{s}}\|X\|_t^{(1-\lambda)\frac{t}{s}}
\qquad\mbox{for}\;s=\lambda r+(1-\lambda)t,\lambda\in(0,1)
\end{align*}
this limit also holds in any $\|\cdot\|_s$ with $1\le s<\frac{3}{2}$.

\medskip

Let us now deduce the almost sure convergence for the minimal energy
\eqref{ao135}: the convergence in Orlicz norm implies by Jensen that
\begin{align*}
|\frac{\mathbb{E}\min E}{L\ln L}+\alpha|\lesssim\|\frac{\min E}{L\ln
L}+\alpha\|_{\frac{3}{2}}\to0.
\end{align*}
Hence it is enough to prove that
\begin{align*}
\frac{\min E}{L\ln L}-\frac{\mathbb{E}\min E}{L\ln L}\to 0\quad\mbox{almost
surely}.
\end{align*}
This follows from this general fact: If a collection of random variables
$\{X_k\}_{k\ge 1}$ satisfies
\begin{align*}
\|X_k\|_s\lesssim k^{-\gamma}~\mbox{for some}~s\ge 1,
\gamma>0\quad\mbox{then}\quad X_k\to0~\mbox{almost surely}.
\end{align*}
By Lemma~\ref{L:10}, this applies to $X_k:=(\min E)/(L\ln L)$ where $k$ is such
that $L=2^k$. The almost sure limits in \eqref{ao136} and \eqref{ao137} follow similarly by
applying \eqref{ao171} and \eqref{c2}.

\medskip

Here comes the proof of the above fact. Notice that
\begin{align*}
\ln\mathbb{P}\big(|X_k|\ge k^{-\frac{\gamma}{2}}\big)\lesssim-k^{\frac{\gamma
s}{2}},\quad\mbox{which implies}\quad\sum_k\mathbb{P}(|X_k|\ge
k^{-\frac{\gamma}{2}})<\infty.
\end{align*}
By the Borel-Cantelli lemma we deduce that almost surely
\begin{align*}
X_k\le k^{-\frac{\gamma}{2}}\quad\mbox{for all large enough $k$},
\end{align*}
which is enough to conclude.

\qed

\section{Acknowledgements}

FO and CW thank Ron Peled for insightful discussions on the white-noise multi-dimensional
case in the Fall of 2023. CW thanks Barbara Dembin for the discussion during a
workshop in Spring 2025.

\medskip

The work was done while the authors were affiliated with the Max Planck Institute for
Mathematics in the Sciences; CW thanks the MPI for the support and warm hospitality.

\medskip

This version of the article has been accepted for publication, after peer review but is
not the Version of Record and does not reflect post-acceptance improvements, or any corrections. The Version of Record is
available online at: https://doi.org/10.1007/s00440-026-01468-y


\newpage

\section{Appendix}

{\sc Proof of Lemma~\ref{L:7}}. To show the existence, it is enough
to notice that almost surely $E$ is a continuous function and admits compact sublevel sets;
indeed, the law of the iterated logarithm in form of
\begin{align*}
W(h) & \lesssim_{\omega, L} \sum_{x = 1}^L \sqrt{|h(x)|\,\ln\ln(e+|h(x)|)}
\lesssim_{L} D(h)^{1/4} \ln \ln (e+D(h)),
\end{align*}
where the subscript $\omega$ indicates that
the constant depends on the realization of the $W$,
yields coercivity of $E$ in the sense of
\begin{align*}
E(h) = D(h) - W(h) \geq \frac{1}{2}D(h)\quad\mbox{whenever}\;D(h) \gg_{\omega, L} 1.
\end{align*}
Since $D(h)$ is a norm on the finite-dimensional and thus locally compact
configuration space (\ref{ao49}), this suffices to show that $E$
has compact sublevel sets.

\medskip

We now turn to the uniqueness statement. Following \cite[Lemma 5.4]{Peled} we introduce
for every site $y\in\{1,\cdots,L-1\}$ and height $\bar h\in\mathbb{R}$
the two constrained minimal energies
\begin{align}\label{ao121}
E_\pm(y,\bar h): = \min\{E(h): h(0) = h(L) = 0\;\mbox{and}\;\pm h(y) \geq \pm \bar h\}.
\end{align}
These are defined such that the event of non-uniqueness in (\ref{ao48}) is contained in 
the union of the events 
\begin{align}\label{u1}
E_+(y,\bar h_+)=E_-(y,\bar h_-)\quad\mbox{while}\quad\bar h_+>\bar h_-,
\end{align}
indexed by $(y, \bar h_+,\bar h_-) \in \{1, \ldots, L-1\} \times
	\mathbb{Q}\times\mathbb{Q}$.
Since this index set is countable,
it is enough to prove that an event of the form (\ref{u1}) has vanishing probability.

\medskip

We thus fix $y\in\{1,\cdots,L-1\}$ and $\bar h_+>\bar h_-$.
We note that by definitions (\ref{ao53}), (\ref{ao48}), and (\ref{ao121}),
\begin{align*}
E_+(y,\bar h_+)+W(y,\bar h_+)&\quad\mbox{depends on}\;W(y,\cdot)\;\mbox{only through}\;
\{W(y,h)-W(y,\bar h_+)\}_{h\ge\bar h_+},\\
E_-(y,\bar h_-)+W(y,\bar h_-)&\quad\mbox{depends on}\;W(y,\cdot)\;\mbox{only through}\;
\{W(y,h)-W(y,\bar h_-)\}_{h\le\bar h_-}.
\end{align*}
Hence writing 
\begin{align*}
\lefteqn{E_+(y,\bar h_+)-E_-(y,\bar h_-)}\nonumber\\
&=\Big(\big(E_+(y,\bar h_+)+W(y,\bar h_+)\big)-\big(E_-(y,\bar h_-)+W(y,\bar h_-)\big)\Big)
-\big(W(y,\bar h_+)-W(y,\bar h_-)\big)
\end{align*}
we learn from the independence of the increments of the (two-sided) Brownian motion $W(y,\cdot)$
that $E_+(y,\bar h_+)$ $-E_-(y,\bar h_-)$ can be expressed as 
the sum of two independent random variables, where one of them is the 
non-degenerate Gaussian $W(y,\bar h_+)-W(y,\bar h_-)$. In particular, the distribution
of the random variable $E_+(y,\bar h_+)$ $-E_-(y,\bar h_-)$ contains no atoms,
and thus the event (\ref{u1}) has vanishing probability.
\ignore{
By shifting the energy $ W $ by $ W ( y, q_1 ) $, we can assume w.~l.~o.~g.~that $ q_1 = 0 $.
The energy can be decomposed as
\begin{align*}
	E(h) & = \Big(D(h) - \sum_{x \neq y} W(x, h(x)) - (W(y,
	h(y)) - W(y, q_2))\Big) + ( W(y, 0)  - W ( y, q_2 ) )
\end{align*}
so that
\begin{align*}
	& E_+(y, q_2) - E_-(y, q_1) \\
	& = \Big(\min\{  E(h) + W(y, q_2)  : h(y) \geq
	q_2\} - E_-(y, 0)\Big) 
	+ (W(y,
	0) - W(y, q_2)) ,
\end{align*}
where the first summand only depends on the increments of $W(y, \cdot)$ below $ q_1 = 0 $
and above $ q_2 $, and the second summand only depends on increments
between $ q_1 $ and $ q_2 $. Thus, since the Brownian motion has independent
increments, these random variables
are independent. The second one is a non-degenerate Gaussian. The law of their
sum, being the convolution of a positive measure with a Gaussian distribution,
has smooth density and the event in \eqref{u1} has null probability.
}
\qed

\medskip

{\sc Proof of Lemma \ref{L:5}}.
We start with the $\gtrsim$ direction in (\ref{ao83}), and argue
\begin{align}\label{ao86}
\mathbb{P}(\hat W_{1}\ge\nu)
\le\mathbb{P}(\hat W_{1}\ge2\nu)
+\mathbb{P}(\frac{D(h_*)}{L}\ge\frac{\nu^\frac{4}{3}}{C})\quad\mbox{for all}\;\nu\ge 0
\end{align}
for some generic\footnote{the value of which may change from line to line} universal constant $C$.
We are using the notation
\begin{align}\label{ao85}
\hat W_{\hat D}:=\sup_{h:\frac{D(h)}{L}\le\hat D}\frac{W(h)}{L}.
\end{align}
Here comes the argument for \eqref{ao86}. By (\ref{ao84}) and (\ref{ao56}) we have
\begin{align*}
\frac{D(h_*)}{L}-\frac{W(h_*)}{L}\le\hat D-\hat W_{\hat D}\quad\mbox{for all}\;\hat D\ge 0.
\end{align*}
By definition (\ref{ao85}) and neglecting $\frac{D(h_*)}{L}\ge 0$ this implies
\begin{align*}
	\hat W_{\hat D}\le \hat W_{\bar D}+\hat D\quad\mbox{for
	all}\;\bar D\ge \frac{D(h_*)}{L}.
\end{align*}
Hence we gather
\begin{align*}
	\mathbb{P}(\hat W_{\hat D}\ge\nu)\le\mathbb{P}(\hat W_{\bar D}\ge\nu-\hat D)
	+\mathbb{P}(\frac{D(h_*)}{L}\ge\bar D)
	\quad\mbox{for all}\;\nu, \bar D\ge0,
\end{align*}
which by the invariance in law (\ref{ao63}) yields
\begin{align*}
\mathbb{P}(\hat W_{1}\ge\frac{\nu}{\hat D^\frac{1}{4}})
\le\mathbb{P}(\hat W_{1}\ge\frac{\nu-\hat D}{\bar D^\frac{1}{4}})
+\mathbb{P}(\frac{D(h_*)}{L}\ge\bar D).
\end{align*}
Redefining $\nu$ we rewrite this as
\begin{align*}
\mathbb{P}(\hat W_{1}\ge\nu)
\le\mathbb{P}(\hat W_{1}\ge\frac{\hat D^\frac{1}{4}\nu-\hat D}{\bar D^\frac{1}{4}})
+\mathbb{P}(\frac{D(h_*)}{L}\ge\bar D),
\end{align*}
which allows us to optimize in $\hat D$ 
by choosing $\hat D$ $=(\frac{\nu}{4})^\frac{4}{3}$ to the effect of 
\begin{align*}
\mathbb{P}(\hat W_{1}\ge\nu)
\le\mathbb{P}(\hat W_{1}\ge\frac{\nu^\frac{4}{3}}{C\bar
D^\frac{1}{4}})
+\mathbb{P}(\frac{D(h_*)}{L}\ge\bar D),
\end{align*}
where we recall that $\bar D> 0$ can be freely chosen.
Hence we pick $\bar D$ such that $\frac{\nu^\frac{4}{3}}{C\bar D^\frac{1}{4}}=2\nu$
which gives (\ref{ao86}).

\medskip

Redefining $\nu$, we may rewrite (\ref{ao86}) as 
\begin{align*}
\mathbb{P}(\frac{\hat W_{1}}{\|\hat D_*\|_s^\frac{3}{4}}\ge\nu)
\le\mathbb{P}(\frac{\hat W_{1}}{\|\hat D_*\|_s^\frac{3}{4}}\ge2\nu)
+\mathbb{P}(\frac{\hat D_*}{\|\hat D_*\|_s}\ge\frac{\nu^{\frac{4}{3}}}{C})
\quad\mbox{where}\quad \hat D_*:=\frac{D(h_*)}{L},
\end{align*}
which, based on the qualitative property that
$\mathbb{P}(\frac{\hat W_{1}}{\|\hat D_*\|_s^\frac{3}{4}}\ge\nu)$ vanishes as $\nu\uparrow\infty$,
implies by iteration
\begin{align*}
\mathbb{P}(\frac{\hat W_{1}}{\|\hat D_*\|_s^\frac{3}{4}}\ge\nu)
\le
\sum_{k=0}^\infty\mathbb{P}(\frac{\hat D_*}{\|\hat D_*\|_s}\ge\frac{(2^k\nu)^{\frac{4}{3}}}{C}).
\end{align*}
By (\ref{ao124}) this implies
\begin{align*}
\mathbb{P}(\frac{\hat W_{1}}{\|\hat D_*\|_s^\frac{3}{4}}\ge\nu)
\le\sum_{k=0}^\infty\exp\big(-(\frac{(2^k\nu)^{\frac{4}{3}}}{C})^s\big)
\quad\mbox{for}\;\nu\ge 1.
\end{align*}
Since the r.~h.~s.~is $\sim \exp(-(\frac{\nu}{C})^{\frac{4}{3}s})$, we obtain
by redefining $\nu$
\begin{align}\label{ao133}
\mathbb{P}(\frac{\hat W_{1}}{C\|\hat D_*\|_s^\frac{3}{4}}\ge\nu)
\le\exp(-\nu^{\frac{4}{3}s})\quad\mbox{for}\;\nu\gg 1,
\end{align}
which by Lemma~\ref{L:8} yields the desired 
$\|\frac{\hat W_{1}}{C\|\hat D_*\|_s^\frac{3}{4}}\|_{\frac{4}{3}s}$ $\lesssim 1$.

\medskip

We now turn to the $\lesssim$ direction in (\ref{ao83}), 
and first note that
\begin{align}\label{ao87}
\Big(\frac{D(h_*)}{L}\in[\hat D, 2\hat D)\quad\Longrightarrow\quad
\sup_{h:\frac{D(h)}{L}\le 2\hat D}\frac{W(h)}{L}\ge \hat D\Big)\quad\mbox{for all}\;\hat D\ge 0.
\end{align}
Indeed, the l.~h.~s.~statement implies that
$\sup_{h:\frac{D(h)}{L}\le 2\hat D}\frac{W(h)}{L}$
$\ge\frac{W(h_*)}{L}$. Since by definition (\ref{ao84}), $(D-W)(h_*)$ $\le(D-W)(0)$ $=0$ holds,
we obtain $\frac{W(h_*)}{L}$ $\ge\frac{D(h_*)}{L}$. Using once more the l.~h.~s.~statement,
we obtain the r.~h.~s.~statement of (\ref{ao87}).

\medskip

By the invariance in law (\ref{ao63}), we obtain from (\ref{ao87})
\begin{align}\label{ao132}
\mathbb{P}(\hat D_*\ge\hat D)\le\mathbb{P}(\hat D_*\ge2\hat D)+
\mathbb{P}(\hat W_1\ge 2^{-\frac{1}{4}}\hat D^\frac{3}{4})
\quad\mbox{for all}\;\hat D\ge 0.
\end{align}
where again
\begin{align*}
\hat D_*:=\frac{D(h_*)}{L}\;\;\mbox{and}\;\;
\hat W_{1}:=\sup_{\hat h:\frac{D(\hat h)}{L}\le1}\frac{W(\hat h)}{L}.
\end{align*}
Following the same steps the led from \eqref{ao86} to \eqref{ao133}, we learn
from \eqref{ao132} that
\begin{align}
\mathbb{P}(\frac{\hat D_*}{C\|\hat W_1\|_s^\frac{4}{3}}\ge\hat D)
\le\exp(-\hat D^{\frac{3}{4}s})\quad\mbox{for}\;\hat D\gg 1,
\end{align}
which again by Lemma~\ref{L:8} yields $\|\hat D_*\|_{\frac{3}{4}s}\lesssim_s\|\hat
W_1\|_s^{\frac{4}{3}}$.
\qed 

\medskip

{\sc Proof of Lemma \ref{L:6}}. Based on the decomposition \eqref{ao67}, 
we inductively construct $h_l$ for 
$l=\frac{L}{2},\frac{L}{4},\ldots,1$ such that
\begin{equation}\label{t1}
	D(h_l) \leq \frac{L}{2}\;\text{for every $l$}\quad\mbox{while}\quad 
	\mathbb{E} W\Big(\sum_l h_l\Big) \gtrsim L \ln{L}.
\end{equation}

\medskip

Assume that $h_{\ge 2l}$ has been constructed. The component $h_l$ is
determined by the values in $(2\hat x - 1) l$ for $\hat x \in \{1, \ldots,
\frac{L}{2l}-1\}$ and we will choose between
\begin{equation}\label{t0}
	h_l((2\hat x - 1)l) = 0 \qquad \text{and} \qquad h_l((2\hat x - 1)l) = l
\end{equation}
by an algorithm that we explain now.
The condition on the Dirichlet energy in~\eqref{t1} is
automatically satisfied since
\[
\frac{D(h_l)}{L} \overset{\eqref{ao112},\eqref{ao27}}{=} \frac{l}{L}\sum_{\hat x =
1}^{\frac{L}{2l}-1} \Big(\frac{h_l((2\hat x -
1)l)}{l}\Big)^2\overset{\eqref{t0}}{\leq}\frac{1}{2}.
\]
For each $\hat x$, define the triangle 
\begin{align*}
	T_{l, \hat x} = \big\{(x, y): & x \in (2(\hat x -1)l, 2\hat x l) \cap
		\mathbb{Z}\\ &\mbox{and}\;\;h_{\ge 2 l}(x) < y < h_{\ge 2 l}(x) + (l - |x - (2\hat x
-1)l|)\big\},
\end{align*}
which is the region between the graphs of $h_{\ge l}$ and $h_{\ge 2l}$ if
we make the choice $h_l((2\hat x - 1)l) = l$ in~\eqref{t0}. We also consider the
upper third part
\begin{align}\label{ao131}
	T'_{l, \hat x} =  \big\{(x, y): & x \in (2(\hat x -1)l +
		\frac{2}{3}l, 2\hat x l - \frac{2}{3}l) \cap
		\mathbb{Z}, \nonumber\\ & h_{\ge 2 l}(x) + \frac{2}{3}l< y < h_{\ge 2
		l}(x) + (l - |x - (2\hat x -1)l|)\big\}.
\end{align}
Following \cite{dw}, for a region $S\subset\mathbb{Z}\times\mathbb{R}$ we denote by
$W(S)$ the two-dimensional semi-discrete white noise integrated over 
$S$. Its law is characterized by this two properties:\footnote{By area of a set, we
mean its semi-discrete version $\text{Area}(S) = \sum_{x \in \mathbb{Z}}
\mathcal{L}^1\big(S \cap (\{x\}\times \mathbb{R})\big)$.}
\begin{align}
&\mbox{$W(S)$ is a centered Gaussian with variance Area$(S)$;}\label{ao128}\\
&\mbox{$W(S_1),\ldots,W(S_n)$ are independent whenever $S_1,\ldots,S_n$ are
disjoint.\label{ao129}}
\end{align}
We now may specify \eqref{t0}: we choose $h_l((2\hat x - 1)l) = l$ if
\[
	W(T'_{l, \hat x}) = \sum_{x \in (2(\hat x -1)l +
		\frac{2}{3}l, 2\hat x l - \frac{2}{3}l) \cap
		\mathbb{Z}}
		W\big(x, h_{\ge 2 l}(x) + (l - |x - (2\hat x
	-1)l|)\big) - W\big(x, h_{\ge 2 l}(x) + \frac{2}{3}l\big) \geq 0
\]
\begin{figure}
\centering
\includegraphics[width=0.45\linewidth]{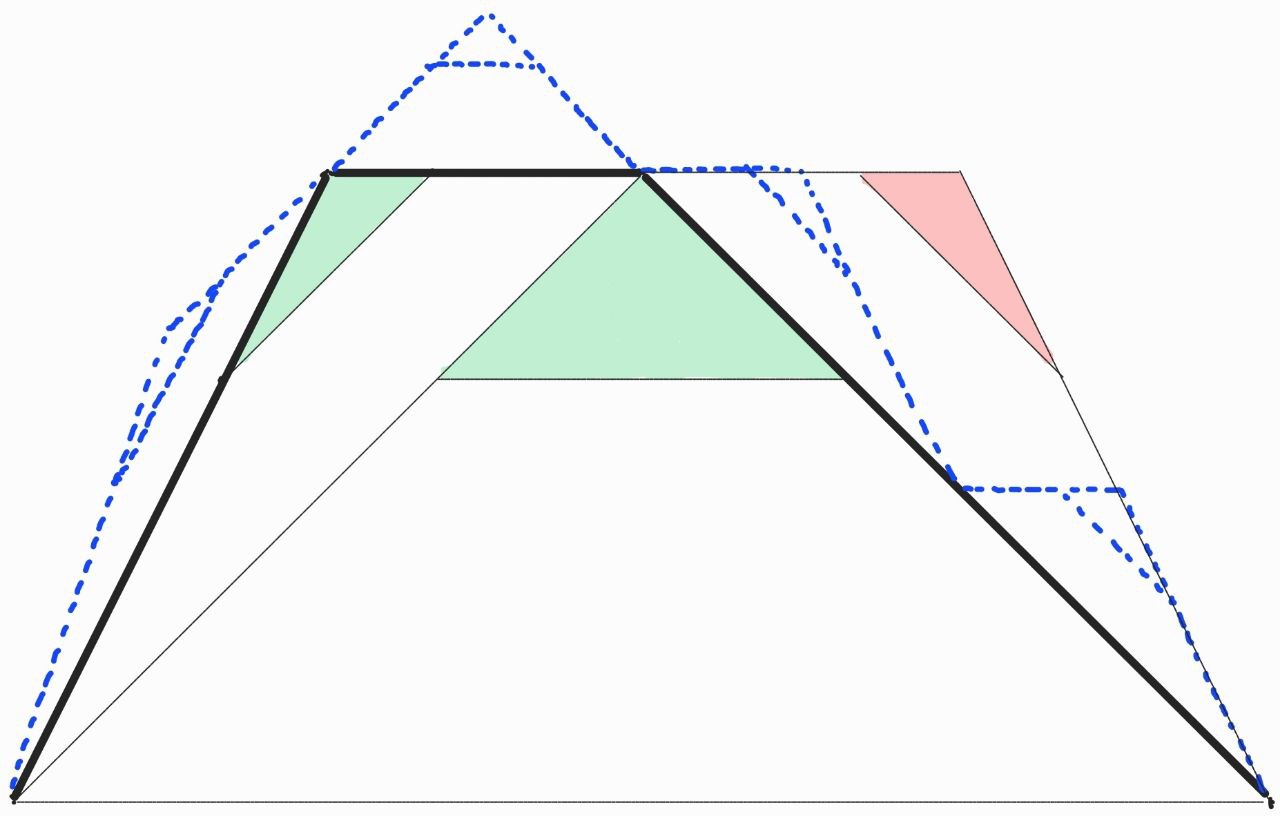}
\caption{Step 3 of the construction (dotted, in blue) in the case
	$W(T'_{\frac{L}{2}, 1}), W(T'_{\frac{L}{4},1})
> 0$, in green, and $W(T'_{\frac{L}{4},2}) < 0$, in red; see also \cite[Figure
1 and 2]{dw} for the more complicated two-dimensional geometric problem.}
\label{fig:f2}
\end{figure}
and $h_l((2\hat x - 1)l) = 0$ otherwise, so that
\begin{equation}\label{t3}
	W(h_{\geq l}) - W(h_{\geq 2 l}) = \sum_{\hat x =
	1}^{\frac{L}{2l}-1} W(T_{l, \hat x})I(W(T'_{l, \hat x}) \geq
	0).
\end{equation}

\medskip

In order to establish the second item in \eqref{t1}, we shall
prove below that
\begin{equation}\label{t2}
	\text{$T_{l, \hat x}$ is disjoint from all the possible
	choices of $\{T'_{\rho, \hat x'}\}_{\rho>l,\hat x'}$}.
\end{equation}
Because of the construction, the triangle $T_{l,\hat x}$ is a deterministic
function of $\{W(T_{\rho,\hat x'})\}_{\rho>l,\hat x'}$, while Area$(T_{l,\hat
x})$ is a deterministic constant. Using this and the two properties of
white noise \eqref{ao128} \& \eqref{ao129}, from \eqref{t2} we learn that
\begin{align*}
	\big(W(T_{l,\hat x}),W(T'_{l,\hat x})\big)\;\mbox{and}\;\{W(T'_{\rho,\hat
	x'})\}_{\rho>l,\hat x'}\;\mbox{are independent},
\end{align*}
which implies by construction of $h_{\ge 2l}$ that
\begin{align*}
	\big(W(T_{l,\hat x}),W(T'_{l,\hat x})\big)\;\mbox{and}\;h_{\ge 2l}\;\mbox{are independent}.
\end{align*}
As a further consequence, we obtain
\begin{align}\label{ao130}
	\big(W(T_{l,\hat x}\backslash T'_{l, \hat x}), W(T'_{l, \hat x})\big)
	\overset{\text{law}}{=}  
	(\sqrt{\text{Area}(T_{l, \hat x}\backslash T'_{l, \hat x})}Z_1,
	\sqrt{\text{Area}(T'_{l, \hat x})}Z_2),
\end{align}
where $(Z_1, Z_2)$ is a couple of independent standard Gaussians. This implies
\begin{align*}
&\mathbb{E}[W(T_{l, \hat x})I(W(T'_{l, \hat x})\geq 0)]\\
&=\mathbb{E}[W(T'_{l, \hat x})I(W(T'_{l, \hat x})\geq 0)]+
\mathbb{E}[W(T_{l, \hat x}\backslash T'_{l,\hat x})I(W(T'_{l, \hat x})\geq 0)]\\
&\overset{\eqref{ao130}}{=}\mathbb{E}[W(T'_{l, \hat x})I(W(T'_{l, \hat x})\geq 0)]\overset{\eqref{ao130}}{\gtrsim}
\sqrt{\text{Area}(T'_{l, \hat x})} \overset{\eqref{ao131}}{\gtrsim} l.
\end{align*}
And summing we obtain
\begin{align*}
	\mathbb{E}[W(h_{\ge l}) - W(h_{\ge 2l})]\overset{\eqref{t3}}{=} \sum_{\hat{x} =
	1}^{\frac{L}{2 l} - 1} \mathbb{E}\big[W(T_{l, \hat x})
	I(W(T'_{l, \hat x})\geq 0)\big] \gtrsim L.
\end{align*}
Summing over all scales we get the lower bound in \eqref{t1}.

\medskip

Here comes the argument for \eqref{t2}. We prove the equivalent
\begin{equation}\label{t4}
	\text{$T'_{l, \hat x}$ is disjoint from all the possible
	choices of $\{T_{\rho, \hat x'}\}_{\rho<l,\hat x'}$}.
\end{equation}
where $l$ must be thought as fixed and $\rho$ ranges over scales
smaller than $l$ (as opposite with respect to \eqref{t2}).
We have two cases, depending on \eqref{t0}. If $h_l((2 \hat x -1)l) =
l$, then all subsequent triangles (corresponding
to $\rho < l$) will be constructed above $T_{l,
\hat{x}}$, so there will be no intersection. If $h_l((2 \hat x -1)l) =
0$, denoting by
\begin{align*}
	h_{\le \frac{l}{2}}:=\sum_{\substack{\rho\;\text{dyadic} \\ 1\le\rho<l}}
	h_\rho
\end{align*}
the part in $h$ below scale $l$, then \eqref{t4} is implied by
\begin{align*}
\sup_x h_{\le \frac{l}{2}}(x) \leq \frac{2}{3}l,
\end{align*}
which we are now going to establish. To do so, it is important to couple
adjacent scales and use the triangular shape of $h_\rho$ together with the
definition in \eqref{t0} to get the finer 
\[
	\sup_x (h_{\rho}(x) + h_{\frac{\rho}{2}}(x)) \leq \rho,
\]
so that we gain a faster decay in the following geometric series
and we conclude
\[
	\sup_x h_{\le \frac{l}{2}}(x) \leq \sum_{k = 0}^{\log_2(\frac{l}{2})} \sup_x
	(h_{l 2^{-2k+1}}(x) + h_{l 2^{-2k}}(x)) \leq 
	\sum_{k = 0}^\infty \frac{1}{4^k}\frac{l}{2} = \frac{2}{3}l.
\]

\medskip

{\sc Proof of Lemma~\ref{L:8}}. Replacing $X$ by $\max\{X-\nu_0,0\}$,
and appealing to the monotonicity and normalization of $\|\cdot\|_s$, cf.~(\ref{ao24}),
it is enough to consider the case of $\nu_0=0$. The argument is based on the identity
\begin{align*}
\mathbb{E}\exp((\frac{X}{C_0})^s)
=1+\int_0^\infty d\nu\mathbb{P}(\frac{X}{C_0}\ge\nu)\frac{d}{d\nu}\exp(\nu^s)
\end{align*}
for some constant $1<C_0<\infty$ to be determined,
which by the l.~h.~s.~of (\ref{ao39}) turns into
\begin{align*}
	\mathbb{E}\exp((\frac{X}{C_0})^s)
&\leq 1+\int_0^\infty d\nu\exp(-(C_0\nu)^s)\frac{d}{d\nu}\exp(\nu^s)\nonumber\\
&=1-\frac{1}{C_0^s-1}\int_0^\infty d\nu\frac{d}{d\nu}\exp(-(C_0^s-1)\nu^s)
=1+\frac{1}{C_0^s-1}.
\end{align*}
Fixing $C_0$ so large that the l.~h.~s.~is $\le e$, we obtain $\|X\|_s\le C_0$
by definition (\ref{ao24}). 

\medskip

Regarding \eqref{ao124}, it is enough to apply Chebyshev's inequality and obtain
\begin{align*}
	\ln\mathbb{P}(X\geq\nu)=\ln\mathbb{P}\Big(\exp\big((\frac{X}{\|X\|_s})^s\big)
	\geq\exp\big((\frac{\nu}{\|X\|_s})^s\big)\Big)
	\leq 1-(\frac{\nu}{\|X\|_s})^s\lesssim-(\frac{\nu}{\|X\|_s})^s
\end{align*}
for $\nu\gg\|X\|_s$.
\qed

\medskip

\begin{lemma}\label{L:15}
The following estimates hold for the Orlicz norms:
\begin{enumerate}
\item[(i)] The following Hölder type inequality holds:
\begin{align*}
\|X\|_s\le\|X\|_r^{\lambda\frac{r}{s}}\|X\|_t^{(1-\lambda)\frac{t}{s}}
\qquad \mbox{for} ~ s=\lambda r+(1-\lambda)t,\lambda\in(0,1) .
\end{align*}
\item[(ii)] As a consequence of the union bound we have
\begin{align*}
	\|\sup_{i\in I}Z_i\|_2\lesssim(\ln\#I)^{\frac{1}{2}}\sup_{i\in I}\|Z_i\|_2 .
\end{align*}
\item[(iii)] For $X,Y$ random variables, $\bar Y\ge 1$ constant and $1\le p\le3$
\begin{align*}
\|I(Y\ge\bar Y)X\|_1&\lesssim(\frac{\|Y\|_1}{\bar Y})^{1-\frac{p}{3}}
\|X\|_{\frac{3}{p}}.
\end{align*}
\item[(iv)]If
\begin{align*}
&\|X~I(X\le\mu\|X\|_s)\|_t\le\mu^\alpha A
\quad\mbox{for}\quad\mu\gg1
\qquad\mbox{with}\quad\frac{st}{s+\alpha t} \ge 1
\end{align*}
then $\|X\|_1\lesssim A$.
\end{enumerate}
\end{lemma}

{\sc Proof of Lemma~\ref{L:15}.} To show (i) we apply Young's inequality
\begin{align*}
|\frac{X}{\|X\|_r^{\lambda\frac{r}{s}}\|X\|_t^{(1-\lambda)\frac{t}{s}}}|^s
=|\frac{X}{\|X\|_r}|^{\lambda r}|\frac{X}{\|X\|_t}|^{(1-\lambda)t}
\le\lambda |\frac{X}{\|X\|_r}|^r+(1-\lambda)|\frac{X}{\|X\|_t}|^t.
\end{align*}
Taking the expectation of the exponential and applying Hölder's inequality, we
conclude
\begin{align*}
&\mathbb{E}\exp(|\frac{X}{\|X\|_r^{\lambda\frac{r}{s}}\|X\|_t^{(1-\lambda)\frac{t}{s}}}|^s)\le
\mathbb{E} \exp(\lambda|\frac{X}{\|X\|_r}|^r)\exp((1-\lambda)|\frac{X}{\|X\|_t}|^t) \\
&\le\big(\mathbb{E}\exp(|\frac{X}{\|X\|_r}|^r)\big)^\lambda
\big(\mathbb{E}\exp(|\frac{X}{\|X\|_t}|^t)\big)^{1-\lambda}\le e.
\end{align*}

\medskip

Let us now prove (ii). It is sufficient to bound the $\sup_iZ_i$
from above, since a bound from below is given by any $Z_j$. By the union bound
\begin{align*}
&\ln\mathbb{P}(\sup_{i\in I}Z_i\ge\nu)
\le\ln\Big(\sum_{i\in I}\mathbb{P}(Z_i\ge\nu)\Big)\\
&\overset{\eqref{ao124}}{\lesssim}\ln(\# I)-(\frac{\nu}{\sup_{i\in I}\|Z_i\|_2})^2
	\qquad&&\mbox{for}\;\nu\gg\sup_{i\in I}\|Z_i\|_2\\
&\lesssim-(\frac{\nu}{\sup_{i\in
I}\|Z_i\|_2})^2\lesssim-
(\frac{\nu}{(\ln\#I)^{\frac{1}{2}}\sup_{i\in I}\|Z_i\|_2})^2\qquad&&\mbox{for}\;\nu\gg( \ln \#
I)^{\frac{1}{2}} \sup_{i\in I}\|Z_i\|_2.
\end{align*}
Using \eqref{ao39} with the change of variable $\nu'=\nu(\ln \# I)^{\frac{1}{2}} \sup_{i\in
I}\|Z_i\|_2$, we conclude.

\medskip

To prove (iii), one can work (up to scaling) with the assumption that
$\|Y\|_1=\|X\|_{\frac{3}{p}}=1$. Let us start from 
\begin{align*}
\mathbb{P}(I(Y\ge\bar Y)X\ge\nu)\le\min\{\mathbb{P}(Y\ge\bar Y),\mathbb{P}(X\ge\nu)\},
\end{align*}
which by (\ref{ao124}) implies that $\ln\mathbb{P}(I(Y\ge\bar Y)X\ge\nu)$ 
$\lesssim-\max\{\bar Y,\nu^\frac{3}{p}\}$ for $\nu\ge 1$.
Noting that $\max\{\bar Y,\nu^\frac{3}{p}\}$ $\gtrsim\bar Y^{1-\frac{p}{3}}\nu$
this yields $\ln\mathbb{P}(\theta I(Y\ge\bar Y)X\ge\nu)$
$\lesssim -\bar Y^{1-\frac{p}{3}}\frac{\nu}{\theta}$ for $\theta\le 1$.
By Lemma~\ref{L:8}, this implies $\|\theta I(Y\ge\bar Y)X\|_1$ $\lesssim 1$
provided $\theta\ll\bar Y^{1-\frac{p}{3}}$, as desired.

\medskip

Here comes the proof of (iv). By scaling, it is enough to treat the case $A=1$. For a threshold $\nu$, we have
\begin{align*}
\mathbb{P}(|X|\ge\nu)\le\mathbb{P}(|X|~I(|X|\le\mu\|X\|_s)\ge\nu)+
\mathbb{P}(|X|\ge\mu\|X\|_s).
\end{align*}
For the second event, we know from \eqref{ao124}
\begin{align*}
\ln\mathbb{P}(|X|\ge\mu\|X\|_s)
\lesssim-\mu^s
\qquad\mbox{for}\;\mu\gg1 ,
\end{align*}
while for the first one \eqref{ao124} implies
\begin{align*}
\ln\mathbb{P}(|X|~I(|X|\le\mu\|X\|_s)\ge\nu)
\stackrel{(iv)}{\lesssim} -(\frac{\nu}{\mu^\alpha})^t\qquad\mbox{for}\;\nu\gg\mu^\alpha.
\end{align*}
Choosing $\mu$ so that
\begin{align*}
\mu^s=(\frac{\nu}{\mu^\alpha})^t\quad\Longleftrightarrow\quad
\mu=\nu^{\frac{t}{s+\alpha t}}
\end{align*}
the two events have probabilities of the same order, so that we obtain
\begin{align*}
	\ln\mathbb{P}(|X|\ge\nu)\lesssim-\nu^{\frac{st}{s+\alpha
	t}}\qquad\mbox{for}~ \nu\gg1.
\end{align*}
Using \eqref{ao39}, this is enough to conclude.
\qed

\medskip

{\sc Proof that \eqref{ao141} implies \eqref{ao142}.} Fix $\bar \nu\gg1$ and $\bar a$ such that
\begin{align}
	&\mathbb{P}(|X-\bar a|\geq \bar \nu)<\frac{1}{2},\qquad\mbox{so that}
	\qquad\mathbb{P}(X\in[\bar a-\bar \nu,\bar a+\bar
	\nu])>\frac{1}{2}\label{ao144}.
\end{align}
For $\nu \geq\bar \nu$, let $a_\nu$ be such that
\begin{align}
\ln\mathbb{P}(|X-a_\nu|\geq \nu)\leq-\nu^s\qquad\mbox{and in
particular}\qquad\mathbb{P}(X\in[a_\nu-\nu,a_\nu+\nu])>\frac{1}{2}. \label{ao146}
\end{align}
From the second items in \eqref{ao144} and \eqref{ao146}, we learn that the intervals $ [ a_{ \nu } + \nu , a_{ \nu } - \nu ] $ and $ [ \bar a + \bar \nu , \bar a - \bar \nu ] $ have non empty intersection. Thus
\begin{align*}
	|a_\nu-\bar a|\leq \bar \nu+\nu\leq2\nu ,
\end{align*}
which, together with \eqref{ao146}, gives
\begin{align}\label{ao147}
	\ln\mathbb{P}(|X-\bar a|\geq3\nu)\leq\ln\mathbb{P}(|X-a_\nu|\geq
	\nu)\lesssim-\nu^s\qquad \nu \geq\bar \nu.
\end{align}
This allows to apply \eqref{ao39} in Lemma \ref{L:8} to $X-\bar a$ to get
$\|X-\bar a\|_s\lesssim_s 1$. By Jensen's inequality
\begin{align*}
|\mathbb{E}X-\bar a|\le\mathbb{E}|X-\bar a|\le\|X-\bar a\|_s\lesssim_s1.
\end{align*}
Together with \eqref{ao147} again, this gives \eqref{ao142}.\qed


\end{document}